\pdfoutput=1
\date{} 
\title{Pivotal, cluster and interface measures for critical planar percolation}
\author{Christophe Garban \and G\'abor Pete \and Oded Schramm}

\documentclass[12pt,preprint,letterpaper]{article}

\newif\iffigures\figurestrue
\figuresfalse

\newif\ifhyper\IfFileExists{hyperref.sty}{\hypertrue}{\hyperfalse}

\ifhyper\usepackage{hyperref}
\def\hitem#1#2{\item[\hypertarget{#1}{#2}]\expandafter\gdef\csname LBL#1ITM\endcsname{#2}}
\def\iref#1{\hyperlink{#1}{\csname LBL#1ITM\endcsname}}
\else
\def\hitem#1#2{\item[{#2}]\expandafter\gdef\csname LBL#1ITM\endcsname{#2}}
\def\iref#1{{\csname LBL#1ITM\endcsname}}
\fi

\newif\ifdraft
\drafttrue 
\long\def\comment#1{}
\long\def\old#1{}

\usepackage{amsmath}
\usepackage{amssymb}
\usepackage{amsthm}
\usepackage{amsfonts}
\usepackage{graphicx}

\usepackage[vmargin={3cm, 3.5cm},hmargin=3cm]{geometry}
\input epsf.sty
\input labelfig.tex

\numberwithin{equation}{section}
\numberwithin{figure}{section}
\newtheorem{theorem}{Theorem}
\numberwithin{theorem}{section}
\newtheorem{corollary}[theorem]{Corollary}
\newtheorem{lemma}[theorem]{Lemma}
\newtheorem{proposition}[theorem]{Proposition}
\newtheorem{conjecture}[theorem]{Conjecture}

\newtheorem{question}[theorem]{Question}

\theoremstyle{remark}\newtheorem{definition}[theorem]{Definition}   
\theoremstyle{remark}\newtheorem{remark}[theorem]{Remark}
\def\eqref#1{(\ref{#1})}
\let\qqed=\qed

\def\nn{\nonumber}

\def\QED{\qqed\medskip}
\let\qed=\QED

\newcommand{\R}{\mathbb{R}}

\newcommand{\C}{\mathbb{C}}
\newcommand{\Z}{\mathbb{Z}}
\newcommand{\N}{\mathbb{N}}
\def\H{\mathbb{H}}

\def\Quad{Q}
\def\CC{\mathscr C}

\def\SS{\mathcal{S}}
\def\diam{\mathrm{diam}}

\def\dist{\mathrm{dist}}
\def\length{\mathop{\mathrm{length}}}

\def\floor#1{\lfloor{#1}\rfloor}
\def\Im{{\rm Im}\,}

\def\SLEkk#1/{$\mathrm{SLE}(#1)$}
\def\SLEr#1/{$\mathrm{SLE(\kappa;#1)}$}
\def\SLEkr#1;#2/{$\mathrm{SLE(#1;#2)}$}
\def\SLEk/{\SLEkk{\kappa}/}
\def\SLEtwo/{\SLEkk2/}
\def\SLE/{$\mathrm{SLE}$}
\def\SLEab/{\SLEkr 4; {a/\hco-1}, {b/\hco-1}/}

\def\Ito/{It\^o}
\def \eps {\epsilon}
\def \P {{\bf P}}
\def\md{\mid}
\def\Bb#1#2{{\def\md{\bigm| }#1\bigl[#2\bigr]}}
\def\BB#1#2{{\def\md{\Bigm| }#1\Bigl[#2\Bigr]}}
\def\Bs#1#2{{\def\md{\mid}#1[#2]}}
\def\Pb{\Bb\P}
\def\Eb{\Bb\E}
\def\PB{\BB\P}
\def\EB{\BB\E}
\def\Ps{\Bs\P}
\def\Es{\Bs\E}

\def \p {{\partial}}
\def \E {{\bf E}}

\def\closure{\overline}
\def\ev#1{{\mathcal{#1}}}
\def \proof {{ \medbreak \noindent {\bf Proof.} }}
\def\proofof#1{{ \medbreak \noindent {\bf Proof of #1.} }}
\def\proofcont#1{{ \medbreak \noindent {\bf Proof of #1, continued.} }}

\def\bl{\bigl}

\def\noopsort#1{}

\def\bl{\begin{lemma}}
\def\el{\end{lemma}}
\def\bth{\begin{theorem}}
\def\eth{\end{theorem}}
\def\bc{\begin{corollary}}
\def\ec{\end{corollary}}
\def\bcj{\begin{conjecture}}
\def\ecj{\end{conjecture}}
\def\bpr{\begin{proposition}}
\def\epr{\end{proposition}}
\def\bde{\begin{definition}}
\def\ede{\end{definition}}
\newcommand{\be}{\begin{eqnarray}}
\newcommand{\ee}{\end{eqnarray}}
\newcommand{\bes}{\begin{eqnarray*}}
\newcommand{\ees}{\end{eqnarray*}}

\usepackage{mathrsfs}

\def\1{1\hspace{-2.55 mm}{1}}

\def\Qual{{\bf Q}}
\def\Sep{\mathcal{T}}
\def\Tg{\mathbb{T}} 
\def\RelQ{\Qual^*}
\def\llra{\longleftrightarrow}
\def\lra{\leftrightarrow}
\def\boxup{\boxdot}
\def\HH{\mathscr{H}}
\def\QUAD{\mathcal{Q}}
\def\T{\mathcal{T}}
\def\A{\mathcal{A}}
\def\CC{\mathscr{C}} 
\def\TC{\mathcal{T}_{\| \cdot \|_\infty}}
\def\ni{\noindent}

\def\param{\tau}
\def\corr{\mathsf{corr}}
\def\aream{\lambda}

\def\area{\mathrm{area}}
\def\length{\mathrm{length}}

\def\Cov{\mathcal{C}}
\def\O{\mathcal{O}}
\def\r{\mathfrak{r}} 

\def\Grimm{Grimmett:newbook}
\def\WWperc{arXiv:0710.0856}

\def\PommerenkeUniv{MR0507768}
\def\Mattila{MR1333890}

\def\KestenScaling{MR88k:60174}
\def\KestenIIC{MR88c:60196}
\def\NolinKesten{MR2438816}
\def\NolinWerner{MR2505301}

\def\SmirnovPerc{MR1851632}
\def\SchrammSmirnovNoise{MR2884873}
\def\SmirnovWerner{MR1879816}
\def\SmirnovICM{MR2275653}
\def\MakarovSmirnovICMP{arXiv:0909.5377}

\def\CFN{MR2284886}
\def\CamiaNewmanFull{MR2249794}
\def\CamiaNewmanConv{MR2322705}
\def\AizXiamen{MR1461346}
\def\SheffieldTrees{MR2494457}
\def\CamiaNewmanIsing{MR2504956}

\def\Hara{MR2393990}
\def\BeffaraDim6{MR2078552}
\def\BeffaraDimKappa{MR2435854}
\def\Sapozhnikov{MR2800910}
\def\Reimer{MR2001g:60017}
\def\Masson{MR2506124}
\def\BarlowMassonTail{MR2683633}
\def\BCKS{MR1868996}
\def\ZhanChordRev{MR2435856}
\def\DubedatCommutation{MR2358649}
\def\DubedatWatts{MR2226888}
\def\SchrammWatts{MR2884875}
\def\LawlerSheffield{MR2884877}
\def\LawlerZhou{arXiv:1006.4936}
\def\LawlerRezaBasic{arXiv:1203.3259}
\def\LawlerRezaMink{arXiv:1211.4146}
\def\LawlerVermesi{Vermesi}

\def\HaggPeresSteif{MR1465800}
\def\SteifSurvey{arXiv:0901.4760}
\def\BKS{MR2001m:60016}
\def\SchrammSteif{MR2630053}
\def\SchShDGFF{MR2486487}
\def\SchSLE{MR1776084}
\def\LawlerStrict{MR1645225}
\def\LSWan{MR1961197}
\def\LSWplane{MR2002m:60159b}

\def\GPS1{MR2736153}

\def\DPSL{DPSL}
\def\MST{MST}
\def\MSTprelim{arXiv:0909.3138}
\def\HPSIIC{HPS:IIC}
\def\OdedSurvey{MR2884872}
\def\MetricProp{MetricProp}

\begin{document}
\maketitle

\begin{abstract}
This work is the first in a series of papers devoted to the construction and study of scaling limits of dynamical and near-critical planar percolation and related objects like invasion percolation and the Minimal Spanning Tree.  We show here that the counting measure on the set of pivotal points of critical site percolation on the triangular grid, normalized appropriately, has a scaling limit, which is a function of the scaling limit of the percolation configuration. We also show that this limit measure is conformally covariant, with exponent 3/4. Similar results hold for the counting measure on macroscopic open clusters (the area measure), and for the counting measure on interfaces (length measure). 

Since the aforementioned processes are very much governed by pivotal sites, the construction and properties of the ``local time''-like pivotal measure are key results in this project. Another application is that the existence of the limit length measure on the interface is a key step towards constructing the so-called natural time-parametrization of the $\mathrm{SLE}_6$ curve.

The proofs make extensive use of coupling arguments, based on the separation of interfaces phenomenon. This is a very useful tool in planar statistical physics, on which we included a self-contained Appendix. Simple corollaries of our methods include ratio limit theorems for arm probabilities and the rotational invariance of the two-point function.
\end{abstract}


\begin{figure}[!htp]
\SetLabels
(0.07*0.12) {\footnotesize left-to-right }\\
(0.052*0.01) {\footnotesize pivotals}\\
(0.09*-0.1) {\footnotesize (dimension $3/4$)}\\
(0.42*0.01) {\footnotesize critical cluster}\\
(0.44*-0.1) {\footnotesize (dimension $91/48$)}\\
(0.83*-0.1) {\footnotesize $\mathrm{SLE}_6$ interface (dimension $7/4$)}\\
\endSetLabels
\begin{center}
\AffixLabels{
\includegraphics[width=0.32 \textwidth]{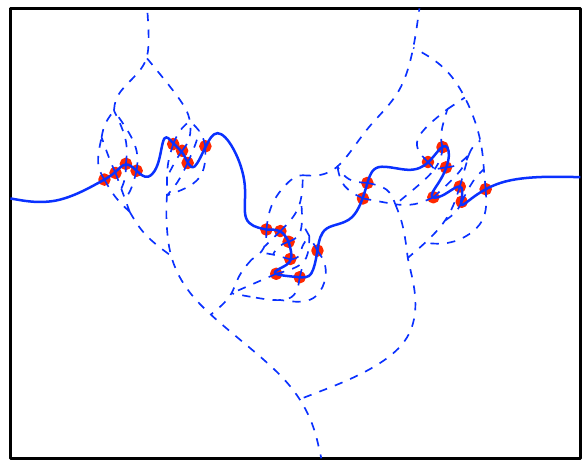}
\includegraphics[width=0.32 \textwidth]{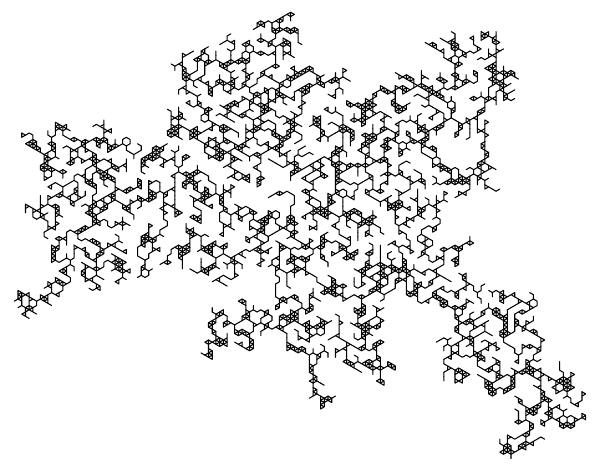}
\includegraphics[width=0.32 \textwidth]{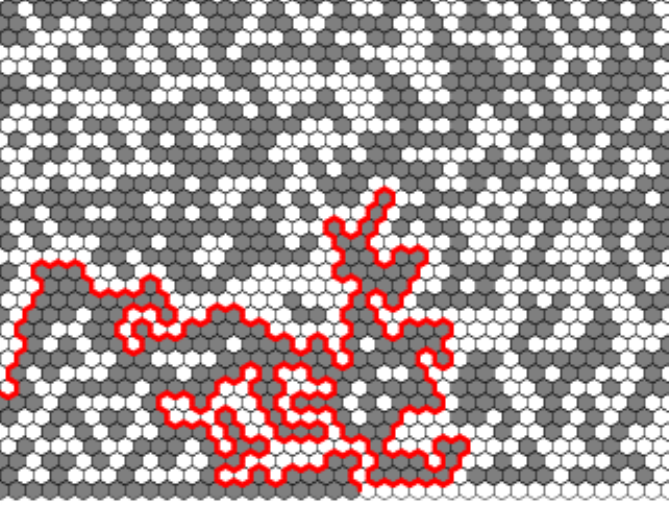}
}
\end{center}
\end{figure}


\tableofcontents

\section{Introduction}\label{s.introduction}

Critical planar percolation has become a central object of probability theory and statistical mechanics; see \cite{\Grimm,\WWperc} for background. A main reason for this is that the discrete process is known to have an interesting continuum scaling limit, in several different but closely related senses, and this scaling limit has turned out to be very useful in understanding also the discrete process, even more than Brownian motion helps in understanding random walks. However, in passing to any notion of a scaling limit, some microscopic information is always lost, and it is far from clear what natural discrete objects (functions of the discrete percolation configuration) will have limits that are meaningful in the scaling limit. This paper focuses on the scaling limits of counting measures on different types of special points of the discrete percolation configuration, normalized by a quantity comparable to their expected value, so that they have a chance to produce a finite limit. The main types of special points that one might think of are the following:
\vskip 0.2 cm
\noindent
1) {\bf macroscopic open clusters}: points that are contained in open clusters that have a macroscopic size (clusters that are still visible in the scaling limit);
\vskip 0.2 cm
\noindent
2) {\bf interfaces}: points on interfaces separating macroscopic open and closed clusters, or the points on a single exploration path that converges to $\mathrm{SLE}_6$;
\vskip 0.2 cm
\noindent
3) {\bf exterior boundaries of clusters}: points on interfaces that are not separated from infinity (or from the boundary of the domain) by the interface itself; for instance, the lowest open left-right crossing in a rectangular domain;
\vskip 0.2 cm
\noindent
4) {\bf pivotal points}: a bit (vertex in site percolation or edge in bond percolation) is called pivotal for some macroscopic crossing event in a given percolation configuration if flipping its state (open versus closed) changes the outcome of the event. 
\vskip 0.2 cm

Having scaling limits for these normalized counting measures is analogous to the construction of the local time measure of one-dimensional Brownian motion at zero via the counting measure on the zeroes of random walks. However, the proofs here are much harder, since we do not have the one-dimensional Markov structure that makes the number of zeroes in two disjoint time intervals independent once we condition on the value of the Brownian motion at any one moment between the two intervals. 

The reason for our paper is the applications of these limit measures. The quantity of pivotal points governs several fundamental dynamical processes related to critical percolation: {\bf dynamical percolation}, the {\bf near-critical ensemble} (percolation in the critical window), and related near-critical objects like {\bf invasion} and {\bf gradient percolation} and the {\bf Minimal Spanning Tree}. Our present work is the first paper in a series devoted to the construction and study of the scaling limits of these processes. 

Why are these processes so fundamental? The scaling limit of dynamical percolation is going to be the natural time evolution with stationary measure being the conformally invariant continuum percolation. The near-critical regime is of central interest in relation with any critical system; e.g., as critical planar systems are often described by conformal (massless) field theories in physics, their near-critical versions are described by massive field theories. Invasion percolation is a self-organized criticality version of percolation. Finally, the Minimal Spanning Tree, besides being a classical combinatorial object, is a natural candidate to possess a scaling limit that is invariant under translations, rotations and scaling, but not under general conformal maps.

Now, how is the scaling limit of the pivotal measure a key to the understanding of the scaling limits of these processes? Any description of continuum critical percolation is intimately related to macroscopic crossing events. Hence, in any reasonable dynamics, the effect of the dynamics on the pivotal points for macroscopic crossing events will probably determine the macroscopic evolution of the system: if, in a given configuration, there are more pivotals for a crossing event, then the dynamics will change that event more easily. However, the amount of pivotals for a macroscopic crossing event is, a priori, microscopic data, so it is not clear that the dynamics will make sense also in the scaling limit. The results of the present paper say that we can in fact tell the number of pivotals from macroscopic information only, and hence we can hope to build the scaling limits of these processes from continuum percolation plus extra randomness governing how pivotals are changing. Such a description also opens the way to understand the conformal properties of the dynamics. These are the goals of this project.

In another direction, the construction of these measures is closely related to the question of finding a so-called {\bf natural} (or {\bf physical}) {\bf time-parametrization} for SLE curves related to percolation, as opposed to the usual parametrization using conformal half-plane capacity. 

In this introduction, we will first discuss these applications in more detail (this time with references),
then explain the results of the present paper, together with the main tools and ideas used in the proofs.

\subsection{Applications to dynamical and near-critical percolation}\label{ss.DPSL}
 
Throughout the paper, we will restrict our attention to critical site percolation on the triangular grid $\eta\Tg$ with small mesh size $\eta>0$, where each vertex of the grid is chosen to be ``open'' or ``closed'' with probability $1/2$, independently of each other. This is the model where conformal invariance \cite{\SmirnovPerc} and its consequences are known: convergence of the exploration interface to $\mathrm{SLE}_6$ \cite{\SchSLE, \SmirnovPerc, \SmirnovICM, \CamiaNewmanConv} and existence of a full scaling limit as a collection of interface loops \cite{\CamiaNewmanFull}. See \cite{\WWperc} for a nice  textbook account of the subject.

Some partial results remain valid for $\Z^2$, for subsequential scaling limits of percolation, which we will 
discuss  at the end of Subsection~\ref{ss.results}.
 
In {\bf dynamical percolation} at criticality, a model introduced by \cite{\HaggPeresSteif} and Itai Benjamini, independently, each bit is changing its state back and forth between open and closed using independent Poisson clocks, with rates in the two directions that make the stationary measure of this process just critical percolation (i.e., for site percolation on $\eta\Tg$, where $p_c=1/2$, the two rates are equal). Two related questions on this process have been studied: {\rm 1) }{\it How does noise affect the macroscopic features of the system (e.g., how much time is needed for the event of having an open left-right crossing in a large square to decorrelate)? {\rm 2)} On the infinite lattice, are there random exceptional times at which the configuration satisfies events that almost surely do not happen in the stationary measure?} See \cite{\OdedSurvey,\SteifSurvey} for surveys of noise sensitivity and dynamical percolation, respectively. For critical planar percolation, both questions have been answered quite completely in a series of three papers over a decade \cite{\BKS, \SchrammSteif, \GPS1}: the time (the amount of noise) needed for decorrelation is given by the typical number of pivotal points, in the following way.
 
A ``quad''  $\Quad$ is a simply connected planar domain with piecewise smooth boundary and four marked points $\{a,b,c,d\}\in\p\Quad$. We think of $\Quad$ as given independently of the grid $\eta\Tg$, with roughly unit size. Site percolation on $\eta\Tg$ can be considered as a two-coloring of the faces of the hexagonal lattice, and hence we can talk about an open crossing in a percolation configuration $\omega_\eta$ from the ``left'' boundary arc $ab$ to the ``right'' boundary arc $cd$ using the hexagons intersecting $\Quad$. Since the probability of this left-right crossing converges to a conformally invariant quantity $\mathsf{cr}_\Quad\in(0,1)$ as $\eta\to 0$, it is very natural to focus on these macroscopic crossing events when considering scaling limits. Now notice that a site $x$ (i.e., a hexagon) is pivotal in a given configuration for left-right crossing if and only if there are open arms from it to $ab$ and $cd$ and closed arms to $bc$ and $da$, i.e., if it satisfies the alternating four-arm event  to $\p\Quad$ (with the arms ending at the respective boundary pieces). For most points $x\in \Quad$ (in the ``bulk''), the probability of this event is comparable  (up to multiplicative factors depending on $\Quad$ and on what we mean by ``bulk'') to $\alpha_4^\eta(\eta,1)$, the alternating four-arm probability from distance $\eta$ to 1 on $\eta\Tg$. See Subsection~\ref{ss.arms} for more on notation and the basics. In fact, it is not hard to show that the expected number of pivotal points in $\Quad$ is $\asymp_\Quad \eta^{-2}\alpha_4^\eta(\eta,1)$, which is known to be $\eta^{-3/4+o(1)}$ on $\eta\Tg$. (Note that it is going to infinity as $\eta\to 0$.) In dynamical percolation, the state of a quad (if it is crossed from left to right or not) changes exactly when the clock of a current pivotal point rings. Therefore, if we take the rate of the clocks in dynamical percolation to be 
$$r(\eta):=\eta^2 \alpha_4^\eta(\eta,1)^{-1}\ (=\eta^{3/4+o(1)} \text{ on }\eta\Tg)\,,$$ 
then, after time $t$, by the linearity of expectation and Fubini, the expected number of {\bf pivotal switches} (i.e., the number of changes of the crossing event) will be of order $t$, with factors depending, of course, on $\Quad$: for a larger quad there are more pivotals, hence more changes, in expectation. This dynamical percolation process will be denoted by $\{\omega_\eta(t)\}_{t\geq 0}$. The previous expectation result already tells us, by Markov's inequality, that in a short time interval $[0,t]$, the probability of seeing any change is small. However, for large $t$, the large expectation does not imply that there is a change with high probability, so it is not immediately clear that this $r(\eta)$ is really the right scaling if we want to follow the changes in macroscopic crossing events.

The possible issues are explained in \cite{\GPS1}, so we give only a brief summary. First of all, the expected number of pivotals for certain Boolean functions is much larger than their typical number. E.g., for the majority function on $2n+1$ bits, $n+1$ bits are pivotal with probability around $1/\sqrt{n}$, but there are no pivotals with probability close to 1. However, it is well-known that this issue does not happen with left-right crossing. This is good, but even if we assumed more, namely, that for a large time $t$ with high probability there are many pivotal switches, it could be the case that pivotal switches happen much faster in one direction than in the other, thus, at any given deterministic moment, the quad is much more likely to be in one of the states (probably the starting one). So, a priori, real macroscopic changes could start happening only at much larger time scales than the one we chose by setting our $r(\eta)$.

The solution to these matters is to use discrete Fourier analysis instead of trying to follow pivotals. With this tool, it is easy to show that at time $t=1$ the correlation 
$$
\corr_\eta(t):=\frac{\Pb{\Quad \text{ is crossed in both }\omega_\eta(0) \text{ and }\omega_\eta(t)}-\Pb{\Quad  \text{ is crossed in }\omega_\eta}^2}{\Pb{\Quad  \text{ is crossed in }\omega_\eta}-\Pb{\Quad  \text{ is crossed in }\omega_\eta}^2}
$$
will already be bounded away from 1, hence a pivotal switch has happened with positive probability, and macroscopic changes have visibly started. (This result follows from (1.2) of \cite{\GPS1}, using only the first and second moments of the size of the so-called Fourier spectral sample of the crossing event, which happen to coincide with the corresponding moments for the pivotal set.) It is much harder to prove, and this is the main result of \cite{\GPS1}, that $\corr_\eta(t)\to 0$ as $t\to\infty$, uniformly in $\eta$. (This needs precise estimates on the lower tail of the size of the spectral sample.) This was proved even for bond percolation on $\eta\Z^2$, while, for the left-right crossing of the square $[0,1]^2$ in site percolation on $\eta\Tg$, also the rate $\corr_\eta(t) \asymp t^{-2/3}$ was established. 

Our present paper is independent of \cite{\GPS1} and does not use Fourier analysis. The above discussion is only to argue that $r(\eta)=\eta^2 \alpha_4^\eta(\eta,1)^{-1}$ is the right scaling for the scaling limit of dynamical percolation. However, we still have to explain what kind of scaling limit we want to establish and how we will do this.

As introduced in \cite{\SchrammSmirnovNoise} and discussed in Subsection~\ref{ss.topology} below, a percolation configuration $\omega_\eta$ on some domain $\Omega$ can be considered as a (random) point in a compact metrizable space $(\HH_\Omega,\T_\Omega)$ that encodes all macroscopic crossing events (the states of all quads). The scaling limit $\omega$ of static percolation is the weak limit of the percolation measures as $\eta\to 0$, w.r.t.~the topology of this compact space. The uniqueness of this quad-crossing limit measure follows from \cite{\CamiaNewmanFull}, using the uniqueness of another notion of scaling limit (based on interface loops), as we will explain in Subsection~\ref{ss.black}. We chose to work in the quad-crossing space $(\HH_\Omega,\T_\Omega)$ because following the states of a countable dense family of quads in dynamical percolation is a very natural way of describing the dynamics; for instance, the correlation decay results of \cite{\GPS1} are formulated for crossing events. Also, this way, $\{\omega_\eta(t)\}_{t\geq 0}$ is a continuous time Markov process in the  {\it compact} space  $(\HH_\Omega,\T_\Omega)$, which certainly makes tightness arguments easier. One can easily describe a natural space-time topology in which the processes $\{\omega_\eta(t)\}_{t\geq 0}$ should converge to some $\{\omega(t)\}_{t\geq 0}$. 

How can we try to describe this scaling limit? For a fixed $\eta>0$, knowing the number of pivotal sites for each quad in the plane at any given moment means that we know the rates at which pivotal switches occur, hence we know the rate of transitions for our Markov process.   However, in the scaling limit there are no pivotal sites anymore that we could count. So, an obvious first task is to show that the scaling limit $\mu^\Quad$ of the normalized counting measure on the pivotals for a quad $\Quad$ in $\omega_\eta$ exists and is a function of the static scaling limit $\omega$. This is the main subject of the present paper. Then, one can hope to build the scaling limit of dynamical percolation from the starting static scaling limit configuration $\omega(0)$ plus independent Poisson clocks with rates given by $\mu^\Quad$ for each $\Quad$: in the dynamics, when a pivotal switch of a quad happens, we change all the crossing events appropriately, read off the new pivotal measure for each quad, and continue. This was the suggestion of \cite{\CFN}, and it has the virtue that it is built from critical percolation (plus some extra randomness), whose conformal properties can be used  to understand the dynamics, as well. Of course, there are some issues here to deal with. For the many small quads inside a bounded domain there are many pivotal sites, hence pivotals for small quads are switching almost all of the time, so there is no nice ordering of all the pivotal switches. One can try to introduce a cut-off $\eps>0$, consider only quads of size larger than $\eps$, then show that the dynamical picture stabilizes as $\eps\to 0$. A related question is whether $\omega(t)$, the states of all the macroscopic quads at any time $t>0$, is determined by the starting macroscopic configuration and the macroscopic pivotal switches in $[0,t]$, or rather there is a cascade of information during the dynamics from the microscopic to the macroscopic world, i.e., microscopic information we have forgotten in the scaling limit becomes macroscopically relevant at later times. Such a cascade would certainly exclude the existence of the macroscopic process  $\{\omega(t)\}_{t\geq 0}$ as a Markov process. These stability results and the existence and Markovianity of the dynamical percolation scaling limit are established in the sequel paper \cite{\DPSL}. In the present paper we also show how the conformal invariance properties of critical percolation translate to conformal covariance properties of the measures $\mu^\Quad$, which will imply similar properties for the dynamical percolation scaling limit.
\vskip 0.3 cm

Dynamical percolation also provides us with tools to understand the {\bf near-critical regime}. Namely, if we modify the dynamics to always switching to ``open'' when a clock rings, then, as time goes on, we are getting a more and more supercritical picture. As times flows in the negative direction, we always turn sites to ``closed'', getting a more and more subcritical picture. All these pictures of different densities are naturally coupled into one ensemble by the dynamics, where, as it is easy to check, the probability of a site being open at time $t\in(-\infty,\infty)$ is 1/2 plus asymptotically to $t \, r(\eta) /2$ as $\eta\to 0$. (Still using the rates $r(\eta)=\eta^2 \alpha_4^\eta(\eta,1)^{-1}\to 0$.)
An equivalent but more standard coupling of these configurations is done by assigning an independent Unif$[0,1]$ label to every site, and defining the configuration $\hat\omega_\eta(\lambda)$ by declaring a site open if its label is at most $1/2 + \lambda\,r(\eta)$. This $\{\hat\omega_\eta(\lambda)\}_{\lambda\in(-\infty,\infty)}$ is called the {\bf near-critical ensemble}, or
{\bf percolation in the critical window}. Why? From the dynamical percolation results of \cite{\GPS1} described above, it clearly follows that for large $\lambda>0$ the probability of having a left-right open crossing in a given quad $\Quad$ is close to 1, while for $\lambda<0$ very negative, the probability is close to 0, uniformly in $\eta$. A priori, this transition could happen even faster, in a smaller window around criticality than the dynamical percolation scale $r(\eta)$, since monotonicity helps us. Nevertheless, it was shown already by Kesten \cite{\KestenScaling} that the four-arm probabilities stay comparable to the critical $\alpha_4^\eta(\eta,1)$ in the entire regime $1/2+\lambda\, r(\eta)$, $\lambda\in(-\infty,\infty)$ (with the factors of comparability getting worse as $|\lambda|\to\infty$), and this implies (using Russo's formula and differential inequalities) that the probability of an open crossing cannot change so fast with $\lambda$ that it reaches 0 or 1 before $\lambda=\pm\infty$. That it does reach 0 and 1 at the endpoints $\lambda=\pm\infty$ of this $r(\eta)$-window was shown in \cite{\BCKS}. This is simpler than establishing the dynamical percolation window, and does not need Fourier analysis, somehow because this can now be understood as a question about a single configuration $\hat\omega_\eta(\lambda)$ instead of a long time interval. Also, the previously mentioned stability result for dynamical percolation that we prove in \cite{\DPSL} can be considered as a generalization of Kesten's stability result for the four-arm probabilities. 

In \cite{\DPSL}, in a way very similar to dynamical percolation, and still roughly following the suggestions of \cite{\CFN}, we prove that the $\eta\to 0$ limit of the near-critical ensemble $\{\hat\omega_\eta(\lambda)\}_{\lambda\in(-\infty,\infty)}$ exists and is a Markov process in $\lambda$. This refines \cite{\BCKS} by showing the convergence of the crossing probabilities, and suggests convergence of the distributions of the sizes of largest components for any fixed $\lambda \in (-\infty,\infty)$. It also refines \cite{\NolinWerner} by showing that the near-critical interface does have an actual scaling limit (the so-called {\bf massive $\mathrm{SLE}_6$}), not just subsequential limits. Massive SLE$_\kappa$ curves are discussed in \cite{\MakarovSmirnovICMP}, but progress is reported only for certain special values of $\kappa$: 2, 4, 8, and some partial results for 16/3 and 3. It is also explained there why these SLE-related curves are ``massive''. We hope that our techniques will give a differential equation for the Loewner driving function of the massive $\mathrm{SLE}_6$, at least on the conjectural level. This has the usual Brownian term $\sqrt{6} B_t$, plus a self-interacting drift term. One drawback is that we do not see how this equation could lead to exact computations like deducing Cardy's crossing formula in the ordinary case.

The near-critical ensemble is also closely related to other objects, like {\bf  gradient} and {\bf  invasion percolation} and the {\bf Minimal Spanning Tree}. The existence and properties of the gradient percolation scaling limit are proved in \cite{\DPSL}, while for the other two, in \cite{\MST}. Since \cite{\MSTprelim} already contains an informal description of the main ideas, we will not discuss these models here.

Finally, let us mention that a very different type of limit object for dynamical percolation, the construction of a local time on the set of exceptional percolation times, will appear in \cite{\HPSIIC}.

\subsection{The question of natural time-parametrization for SLE curves}\label{ss.time}

The item labels 1)-4) in the very first paragraph of this Introduction correspond to the number of macroscopic arms that a corresponding special point must have. So, it should not be surprising if analogs of the above 4-arm results hold also for other $k$-arm events. In the present subsection, we discuss the 2-arm case in a bit more detail, since it is interesting also from a different point of view.

SLE's have been extremely useful in understanding discrete processes with conformally invariant scaling limits. But there is also an opposite direction: if an SLE is known to be the scaling limit of a discrete process, then we often gain valuable extra information about the SLE. For instance, the reversibility of SLE$_\kappa$ for the special values $\kappa=2,4,6,8$ was known, using the discrete processes converging to them, much before the general $\kappa\leq 4$ result of Dapeng Zhan \cite{\ZhanChordRev}. The question of natural time-parametrizations of SLE curves has to do with both directions.

By definition, the standard time-parametrization for SLE$_\kappa$ is half-plane capacity for the chordal and disk capacity (conformal radius) for the radial processes. These parametrizations transform nicely under conformal maps, hence are often useful when working with SLE's --- except, for instance, when one wants to consider the path traversed in the reverse direction. On the other hand, the discrete processes converging to SLE's are naturally equipped with a time parametrization, given by the number of discrete lattice steps, which obviously has nothing to do with the capacity parametrization of the scaling limit. For instance, it has a locality property that capacity does not have, namely, congruent pieces in two curves have the same length, regardless of the rest of the curves. If the limit of the discrete time-parametrization existed, it should inherit all these natural properties, and, in fact, it should be the $d(\kappa)$-dimensional Minkowski content, where $d(\kappa)=1+\kappa/8$ is the dimension of the SLE$_\kappa$ \cite{\BeffaraDimKappa}. 

After weaker results by \cite{\LawlerSheffield}, it  was proved in \cite{\LawlerZhou} that SLE$_\kappa$ for all $\kappa < 8$ has a unique time-parametrization that satisfies a few natural conditions, such as the right $d(\kappa)$-dimensional transformation rule (i.e., $d(\kappa)$-dimensional conformal covariance). An improved proof, which also gives some information about regularity properties of the parametrization, appeared recently in \cite{\LawlerRezaBasic}. Even more recently, it has been proved in \cite{\LawlerRezaMink} that this time-parametrization equals the $d(\kappa)$-dimensional Minkowski content. However, it is not known for a single $\kappa$ that this parametrization is a limit of the natural discrete parametrization of counting steps.

In the present paper, the scaling limit of the normalized counting measure on a percolation interface (see~(\ref{e.Parameasure}) below) will be handled quite similarly to the case of pivotals, with two main differences: one, points on an interface between two boundary points $a,b\in\p\Quad$ are characterized by having one open and one closed arm to the respective boundary arcs $ab$ and $ba$, instead of four alternating arms; two, it turns out to be simpler and more natural to describe the scaling limit of the measure as a function of the $\mathrm{SLE}_6$ interface, rather than a function of the entire percolation configuration in the quad-crossing topology (see Subsection~\ref{ss.interface}). Now, having constructed this scaling limit and having proved that it has the right conformal covariance properties, it looks quite obvious that it is the natural time-parametrization we are looking for. However, there is an issue: we have the length measure of the entire interface from $a$ to $b$ and the measure restricted to small boxes, but any small box can contain several connected pieces of the interface, and the measure of the box does not differentiate between the pieces.  So, it is not clear that we can also get the length measure of all the starting segments of the interface.  This question, together with fine metric properties (e.g., H\"older-regularity) of the length measure will be addressed in \cite{\MetricProp}. 

\subsection{The main results and techniques}\label{ss.results}

Having been bombarded by two subsections of future applications, the reader hopefully has a good idea by now what the results themselves are. Nevertheless, here they are explicitly. For any missing definition, see Section~\ref{s.backgrounds}.

Let $\Omega\subset \hat\C=\C\cup\{\infty\}$ be open, and $(\HH_\Omega,\T_\Omega)$ be its quad-crossing space. For any quad $\Quad\subseteq \Omega$ with piecewise smooth boundary, consider the following measure on $\Quad$, a function of the discrete percolation configuration $\omega_\eta \in \HH_\Omega$:
\begin{equation}\label{e.Qmeasure}
\mu^\Quad_\eta := \mu_{\eta}^\Quad (\omega_{\eta})= \sum_{x\textrm{ is } \Quad\textrm{-pivotal}} \delta_x\, \frac {\eta^2}{ \alpha_4^\eta(\eta,1)}\,,
\end{equation}
where $\delta_x$ is a unit mass on the hexagon $x$ (say, a point pass at the center of $x$). This $\mu^\Quad_\eta$ is an element of the space $\mathfrak{M}(\Quad)$ of finite measures on the compact set $\Quad$, which is a complete, metrizable, separable space in the weak*-topology (also called the weak convergence of measures). 

Now, our main results concerning pivotals are the following:

\begin{theorem}[{\bf Measurable limit measure}]\label{th.4meas}  For any quad $\Quad\subseteq \Omega$ with piecewise smooth boundary, there is a limit $(\omega,\mu^\Quad)$ of the joint law $(\omega_\eta,\mu^\Quad_\eta)$ in the product topology of the quad-crossing space $(\HH_\Omega,\T_\Omega)$ and the weak*-topology on $\mathfrak{M}(\Quad)$, where $\mu^\Quad$ is a measurable function of the percolation scaling limit $\omega$.
\end{theorem}

\begin{theorem}[{\bf Conformal covariance}]\label{th.4cov} 
Let $f:\Omega\longrightarrow \widetilde\Omega$ be a conformal isomorphism, and assume that $\Quad\subset \closure\Quad \subset \Omega$. Let $\tilde\omega=f(\omega)$ be the image of the continuum percolation $\omega\in \HH_\Omega$, and  $f_*(\mu^\Quad(\omega))$ be the pushforward measure of $\mu^\Quad$. Then, for almost all $\omega$, the Borel measures $\mu^{f(\Quad)}(\tilde\omega)$ and $f_*(\mu^\Quad(\omega))$  on $f(\Quad)$ are absolutely continuous w.r.t.~each other, and their Radon-Nikodym
derivative satisfies, for any $w=f(z) \in \widetilde\Omega$,
\begin{equation*}
\frac {d \mu^{f(\Quad)}(\tilde\omega)}{ d f_*(\mu^\Quad(\omega))}\,(w) = |f'(z))|^{3/4},
\end{equation*}
or equivalently, for any Borel set $B\subset \Quad$,
\begin{equation*}
\mu^{f(\Quad)}(f(B))(\tilde\omega) = \int_{B} |f'|^{3/4} d\mu^\Quad(\omega) \,.
\end{equation*}
\end{theorem}

On a heuristic level, the scaling exponent $3/4$ comes from the fact that $\mu^\Quad$ is a ``natural'' measure supported on the set of pivotal points, which is known to be of Hausdorff-dimension $3/4$. (This dimension is not proved explicitly anywhere, but follows from a simple modification of \cite[Section 6.1]{\LSWplane} and the techniques of \cite{\BeffaraDim6}. See also \cite[Theorem 10.4]{\GPS1} for a very similar result. For more on the connections between the dimension of a set and the measures supported, see \cite[Chapter 8]{\Mattila}.) Of course, to make this explanation more grounded, one would also need to prove that the support of the measure $\mu^\Quad(\omega)$ is exactly the set of $\Quad$-pivotal points of $\omega$, not just a subset. This requires some work, which will be done in \cite{\MetricProp}.

Note that it might be more natural to 
use a different family of measures: instead of a measure for each quad $\Quad\subseteq \Omega$, there could be a measure $\bar \mu^\rho_\eta$ for each radius $\rho>0$, the normalized counting measure on the {\bf $\rho$-important sites}, i.e., sites $x\in \Omega$ that have the alternating four-arm event from $x$ to distance at least $\rho$. However, these $\rho$-important measures are slightly cumbersome to deal with, so we will introduce one more family of measures: the normalized counting measure $\mu^A_\eta$ on the so-called {\bf $A$-important sites}, where $A\subset \Omega$ is a topological annulus with piecewise smooth boundaries (called a proper annulus).  See Subsection~\ref{ss.A-imp} for the definitions, Theorems~\ref{th.measurable} and~\ref{pr.covariance} for the results, and Subsection~\ref{ss.filtering} for a comparison between the $\rho$-important and $A$-important families of measures.
\medskip

The case of the 1-arm event, i.e., the normalized counting measure on {\bf macroscopic open clusters}, is simpler than the case of pivotals, due to the monotonicity of the 1-arm event that is false for the alternating 4-arm event.  The actual formulation that we will prove concerns the following measure (using the 1-arm analogue of the $A$-important points mentioned above):
\begin{equation}\label{e.Areameasure}
\aream_\eta^A = \aream_\eta^A(\omega_\eta) := \sum_{x \in \Delta:\; \{x \leftrightarrow \p_2 A\} } \delta_x \, \frac {\eta^{2}} {\alpha_1^\eta(\eta,1)}\,,
\end{equation}
where $A$ is a proper annulus with outer boundary $\p_2 A$ and inner disk $\Delta$. See Subsections~\ref{ss.clusters} and~\ref{ss.cothers} for the results. Since $\alpha_1^\eta(\eta,1)=\eta^{5/48+o(1)}$, the covariance exponent is $91/48$ in this case.

Let us note that the cluster measure was also considered in
\cite{\CamiaNewmanIsing} for percolation and the FK-Ising model. 
In the case of FK percolation, as explained there,
proving a scaling limit result (similar to ours) for the area measure
on FK clusters would enable one to construct a magnetization field for
the scaling limit of the Ising spin model. In \cite{\CamiaNewmanIsing}, Camia and Newman 
sketched the proof for the tightness of the normalized discrete
measures, while the main step in this program, i.e., the question of convergence,  
is proved in \cite{CGN}.
\medskip

The 2-arm case, i.e., the counting measure on the {\bf interfaces}, will be slightly different. The setting is that $\Omega$ is a domain with smooth boundary and the so-called Dobrushin boundary conditions: two marked points $a,b\in\p\Omega$, and we are interested in the set of edges of the triangular grid that have one endpoint connected to the boundary arc $ab$ by an open path, and another endpoint connected to the arc $ba$ by a closed path, in the configuration $\omega_\eta$. Of course, one usually considers dual edges, i.e., edges on the hexagonal lattice that form the boundary of the open cluster of hexagons connected to $ab$, thus forming an interface $\gamma_\eta$ from $a$ to $b$. The normalized counting measure is now
\begin{equation}\label{e.Parameasure}
\param^{\Omega,a,b}_\eta = \param^{\Omega,a,b}_\eta(\omega_\eta) = \param^{\Omega,a,b}_\eta(\gamma_\eta):= \sum_{\mathrm{edges} \, e \in \gamma_\eta}  \delta_e \, \frac {\eta^2} {\alpha_2^\eta(\eta,1)}\,.
\end{equation}

As mentioned in Subsection~\ref{ss.time}, it will be easier and more natural to show that the limit measure $\param^{\Omega,a,b}$ exists as a function of not the continuum percolation $\omega$, an element of the quad-crossing space $\HH_\Omega$, but as a function of the $\mathrm{SLE}_6$ limit of the interface $\gamma_\eta$. See Subsection~\ref{ss.interface} for more details and the limit measure result, and Subsection~\ref{ss.cothers} for conformal covariance. The covariance exponent is $7/4$.
\medskip

It would also be natural to consider the 3-arm case, i.e., points on the {\bf exterior boundary} of large clusters or of a single interface. However, this has some additional technicalities that would make this paper even longer, while the results do not seem to have striking applications at present, hence we chose not to address this case in detail. See Subsection~\ref{ss.frontier} for what the issues are and for hints on how they could be solved.
\medskip

Let us now explain the main ideas in the proofs of Theorems~\ref{th.4meas} and~\ref{th.4cov}, together with some lemmas that are interesting also in their own rights. 

First of all, as we said in Subsection~\ref{ss.DPSL}, we want to prove the existence and conformal covariance of a measure that is measurable with respect to the Borel-$\sigma$-field of the {\bf quad-crossing space} $(\HH_\Omega,\T_\Omega)$, see Subsection~\ref{ss.topology}. Therefore, it is essential to know that the discrete percolation configurations $\omega_\eta$ themselves have a unique scaling limit $\omega$ in this space. (As it will be apparent from the proof, if we had only subsequential limits of percolation, then we could not show the conformal covariance of the subsequential limits of the measures.) However, \cite{\SchrammSmirnovNoise},  where this space was introduced, proves the existence and the ``black noise factorization property'' only along subsequences. Fortunately, the {\bf uniqueness of the quad-crossing limit} follows easily from the results of \cite{\CamiaNewmanFull}: Camia and Newman showed there that the percolation configuration described as the set of all {\bf interface loops} has a unique and conformally invariant scaling limit, and they also showed that the $\mathrm{SLE}_6$ interface in a subdomain with Dobrushin boundary conditions can be extracted from their scaling limit. This immediately implies that the quad-crossing events are also measurable in this scaling limit. So, the quad-crossing scaling limit is a function of a unique scaling limit object, hence is unique itself. This will be explained more in Subsection~\ref{ss.black}. From this we also get that the quad-crossing scaling limit is conformally invariant, which would not follow simply from Smirnov's theorem of the conformal invariance for a single quad \cite{\SmirnovPerc}: we need conformal invariance of the joint laws, not only the one-dimensional marginals provided by the Cardy-Smirnov formula. After uniqueness, we need to know that, e.g., the four-arm event in an annulus is actually a measurable function of the quad-crossing configuration $\omega$ and its probability is the $\eta\to 0$ limit of the probabilities in the $\omega_\eta$ configurations. Though not obvious, this is not hard, and is done in Subsection~\ref{ss.limitarms}. All these measurability results are parts of the general belief that the different notions of the full scaling limit for percolation are all equivalent. See \cite{\SchrammSmirnovNoise} for an overview of these matters.

Given the measurability results of the previous paragraph, one can start showing that the number of $\Quad$-pivotals in any ball $B\subseteq\Quad$ can be read off from macroscopic quad-crossing events only, so that the values $\mu^\Quad_\eta(B)$ will have a chance to have a limit that is a measurable function of $\omega$. The main idea is to take a grid of macroscopic but small mesh size $\eps$, count the number $Y^\eps_\eta$ of $\eps$-boxes in the grid that lie in $B$ and have the four-arm event to $\p\Quad$, and hope that, for small $\eps$, this number will give a good idea of the number $X_\eta$ of actual pivotal sites in $B$. To realize this hope, one needs to show that given that an $\eps$-box $Q_i$ has the four-arm event to $\p\Quad$, the distribution of the pivotal points inside $Q_i$ is more or less independent of the percolation configuration outside an enlarged box $kQ_i$, with some large $k$. With enough independence, one could try to prove a law of large numbers type result, saying that the contributions of the many $\eps$-boxes average out, so that $X_\eta$ will be well-approximated (up to a factor $1+o(1)$) by $Y^\eps_\eta$ times a deterministic factor $\beta(\eta,\eps)$, which is basically the expectation of the number of pivotals in an $\eps$-box given that the box has the four-arm event to $\p\Quad$. How do we fill in the details? One way to express the approximate independence is that given that $Q_i$ has the four-arm event to $\p\Quad$, and given two different copies of the entire percolation configuration outside $kQ_i$ that satisfy the conditioning, the configurations inside $Q_i$ can be coupled with large probability to agree in the two copies. Such {\bf coupling results} are proved in Section~\ref{s.coupling}, using the {\bf separation of interfaces phenomenon} that is discussed in Appendix~\ref{s.appendix}. This coupling property is closely related (for the one-arm event) to Kesten's construction of the Incipient Infinite Cluster \cite{\KestenIIC}, and (for the multi-arm events, where the separation of interfaces starts playing a role) to the construction of the multi-arm IIC measures by \cite{\Sapozhnikov}. In fact, the proof strategy in our Appendix~\ref{s.appendix} is based on \cite{\Sapozhnikov}. Then, we use our coupling results in Section~\ref{s.measurable} to gain enough independence to prove the approximation of $X_\eta$ by $\beta(\eta,\eps) Y^\eps_\eta$ in the $L^2$-sense.
Then there is some additional work in showing that the normalized counting measures $\mu^\Quad_\eta(B)$ actually have a limit. The existence of this limit measure, and also its scale invariance, use the following {\bf ratio limit theorems}, proved in Subsection~\ref{ss.ratio}: 
\begin{equation}\label{e.intratio}
\lim_{\eta\rightarrow 0}\frac{\alpha_4^\eta(\eta,r)}{ \alpha_4^\eta(\eta, 1)} = \lim_{\eps\to 0} \frac{\alpha_4(\eps,r)}{\alpha_4(\eps,1)} = r^{-5/4}\,,
\end{equation}
where $\alpha_4$ without the superscript $\eta$ is the continuum four-arm probability, given by $\mathrm{SLE}_6$. See Remark~\ref{r.ratioexplain} for an explanation why this is not at all obvious from the four-arm exponent being 5/4. Similar limit results hold for the $\alpha_1$ and $\alpha_2$ probabilities, and we prove the {\bf rotational invariance of the two-point function} $\Pb{x\mathop{\llra}\limits^{\omega_\eta} y}$ in the $\eta\to 0$ limit, see Subsection~\ref{ss.twopoint}.
\medskip

For the proof of conformal covariance in Section~\ref{s.covariance}, we have to deal with the fact that a conformal map distorts the $\eps$-grid that we use for the approximation, so we need to compare such a distorted approximation with the original. Locally, each $\eps$-box alone is not much distorted, almost just shifted, rotated and scaled, so we can start by proving the approximation results for all translated, rotated and scaled versions of the $\eps$-grid, which are still proper grids. However, if we tried to approximate locally the distorted grid everywhere at the same time, by covering the domain by small pieces of these proper grids, it would be impossible to glue these pieces so that we get a one-fold covering of the domain, counting every pivotal exactly once. The solution is to take an average of all the shifted and rotated approximations, w.r.t.~the uniform measure on the shifts $a\in [0,\eps)^2$ and angles $\theta\in[0,2\pi)$, and notice that if we take any conformal map $f$, then locally at each point $f(x)\in f(\Omega)$, in the push-forward of the average we will still see the uniform distribution on the shifted and rotated $|f'(x)|\eps$-boxes (or rather almost-boxes). (For any given pair $(a,\theta)$, the $a$-shifted $\theta$-rotated $|f'(x_1)|\eps$-box at some point $f(x_1)$ will come from a differently shifted and rotated $\eps$-grid than the $(a,\theta)$-moved $|f'(x_2)|\eps$-box at some other point $f(x_2)$, since the local rotation and shift given by $f$ is different at $f(x_1)$ and $f(x_2)$; however, the overall distributions are still uniform at each point.)  Therefore, the map $f$ changes this average by the integral of the right power of the local scaling $|f'|$. (To compare different scaled versions of a grid, (\ref{e.intratio}) is used again.) This will prove conformal covariance.
\medskip

Finally, what survives for other lattices, say, {\bf critical bond percolation on $\Z^2$}? Since the Russo-Seymour-Welsh technology works, all ingredients remain valid that do not use the uniqueness and conformal invariance of the percolation scaling limit: the basic topology and measurability results for the quad-crossing space $(\HH_\Omega,\T_\Omega)$, the tightness and moment estimates for the normalized counting measures, and the separation of interfaces phenomena and the coupling results. (See Remark~\ref{r.Z2} for a tiny issue to handle.) Therefore, we can still prove that whenever the scaling limit of static percolation exists along a subsequence $\{\eta_k\}_{k=1}^\infty$, then the scaling limit of the normalized measures also exists along that subsequence. Then we conclude in \cite{\DPSL} that the dynamical and near-critical percolation scaling limits also exist along the same subsequence and are Markov processes. Note that the mere existence of subsequential scaling limits for these processes would follow from the compactness of $(\HH_\Omega,\T_\Omega)$, but the possibility of using the same subsequence that already worked for the static scaling limit, and the Markovianity of the limit are highly non-trivial.

\medskip
\noindent {\bf Acknowledgments.} We are grateful to Art\"em Sapozhnikov for explaining to us the simple proof of the Strong Separation Lemma (Proposition~\ref{p.separation}) in emails, before it was written up in \cite{\Sapozhnikov}, to Vincent Beffara and Chuck Newman for explaining that the results of \cite{\CamiaNewmanFull} imply the uniqueness of the quad-crossing scaling limit, and to Asaf Nachmias for encouraging us to include the statement about the rotational invariance of the two-point connectivity function. We also thank Ren\'e Conijn, Demeter Kiss, and Alan Hammond for very useful comments on the manuscript, and Vincent Beffara for the simulation pictures. A large part of the work was done while all authors were at Microsoft Research.

\section{Arms, quads and loops}\label{s.backgrounds}

\subsection{Notation and basics}\label{ss.arms}

Recall that we can consider site percolation on $\eta\Tg$ as a two-coloring of the faces of the hexagonal lattice, and hence of the plane (with edges shared by two differently colored faces having both colors. In this context, we will use the words ``hexagon'', ``site'' and ``point'' interchangeably. We assume that the reader is familiar with the basics of discrete critical planar percolation: the FKG-inequality, the Russo-Seymour-Welsh (RSW) technology, and the mono- and polychromatic $k$-arm events. See \cite{\WWperc}. 

We will use the multi-arm events and probabilities in several different versions. First of all, $\alpha_k^\eta(r,R)$ will denote the probability that there are $k$ arms, not all of the same color when $k>1$, connecting the two boundary pieces of the annulus $A(r,R)=B_R\setminus B_r$, where $B_r$ will denote the square of Euclidean radius $r$ (hence sidelength $2r$). Of course, there is a small problem with this definition: the exact set of hexagons intersecting this annulus depends on where the center of the two squares is compared to the discrete lattice. For a while, the precise definition will not matter, as long as it is fixed; say, the center of the squares is the center of a hexagon. Later this will matter more, at which point we will refine the definition. The notation $\alpha_k^\eta(\eta,R)$, for $k\leq 6$, will always denote the $k$-arm probability from a single site to radius $R$.
Omitting the superscript $\eta$ will either mean (depending on the context) that we are talking about the continuum $\mathrm{SLE}_6$ quantity $\alpha_k(r,R)=\lim_{\eta\to 0}\alpha_k^\eta(r,R)$, or that it does not matter whether one uses the continuum or the discrete quantity. There will be several more versions of the four-arm event later on; see, for instance, Subsection~\ref{ss.ratio}.

A well-known important property of $k$-arm probabilities is {\bf quasi-multi\-pli\-cativity}:
$$
c_k \, \alpha_k(r_1,r_2) \, \alpha_k(r_2,r_3) \leq
\alpha_k(r_1,r_3) \leq 
\alpha_k(r_1,r_2)\, \alpha_k(r_2,r_3)\,,
$$
for some absolute constant $c_k>0$. For $k=1$ this is easy from RSW and FKG, for $k\geq 2$ it requires the (weak) separation of interfaces technology; see \cite{\KestenScaling, \SchrammSteif, \NolinKesten}.

\subsection{The quad-crossing topology}\label{ss.topology}

We now define what we mean by the scaling limit of critical percolation. We will work with the setup introduced in \cite{\SchrammSmirnovNoise}, which describes the scaling limit using ``left-right'' {\bf crossing events} in generalized quadrilaterals. For completeness, we review the basic definitions and results below. There are alternative setups that might work equally well, even if we felt that using crossing events is the most practical one for our purposes: Aizenman \cite{\AizXiamen} suggested using the set of {\bf all open paths}, Camia and Newman \cite{\CamiaNewmanFull} defined the scaling limit as the set of all {\bf interface loops}, and Sheffield defined a {\bf branching SLE tree} \cite{\SheffieldTrees}, which should be the limit of discrete exploration trees in all models where the interfaces converge to SLEs.  As pointed out in \cite{\SchrammSmirnovNoise}, where there is a discussion of several more possible definitions, it is a priori highly non-trivial to prove that all these different notions contain the same relevant information (i.e., that they are measurable with respect to each other). Nevertheless, as we will explain in Subsection~\ref{ss.black}, the results of \cite{\CamiaNewmanFull} imply that the quad-crossing scaling limit is a function of the Camia-Newman scaling limit, while Subsection~\ref{ss.limitarms} is concerned with what objects can be expressed in the quad-crossing limit. (We also have a proof that the quad-topology is measurable with respect to the branching $\mathrm{SLE}_6$ tree, but that is more cumbersome than using the continuum loop ensemble, especially that \cite{\CamiaNewmanFull} has already established the uniqueness and several properties of the latter, hence we have not written up that proof.)

Let $D\subset \hat\C=\C\cup\{\infty\}$ be open. A {\bf quad} in the domain $D$ can be considered as a homeomorphism $\Quad$ from $[0,1]^2$ into $D$. The space of all quads in $D$, denoted by $\QUAD_D$, can be equipped with the following metric: $d_\QUAD(\Quad_1,\Quad_2):=\inf_{\phi}\sup_{z\in \p [0,1]^2} |\Quad_1(z)-\Quad_2(\phi(z))|$, where the infimum is over all homeomorphisms $\phi: [0,1]^2 \longrightarrow [0,1]^2$ which preserve the 4 corners of the square. A {\bf crossing} of a quad $\Quad$ is a connected closed subset of $[\Quad]:=\Quad([0,1]^2)$ that intersects both $\p_1\Quad=\Quad(\{0\}\times[0,1])$ and $\p_3\Quad=\Quad(\{1\}\times[0,1])$.

\begin{figure}[htbp]
\SetLabels
(.25*.02)$\Quad_1$\\
(1.02*.74)$\Quad_2$\\
\endSetLabels
\centerline{
\AffixLabels{
\includegraphics[height=1.5 in]{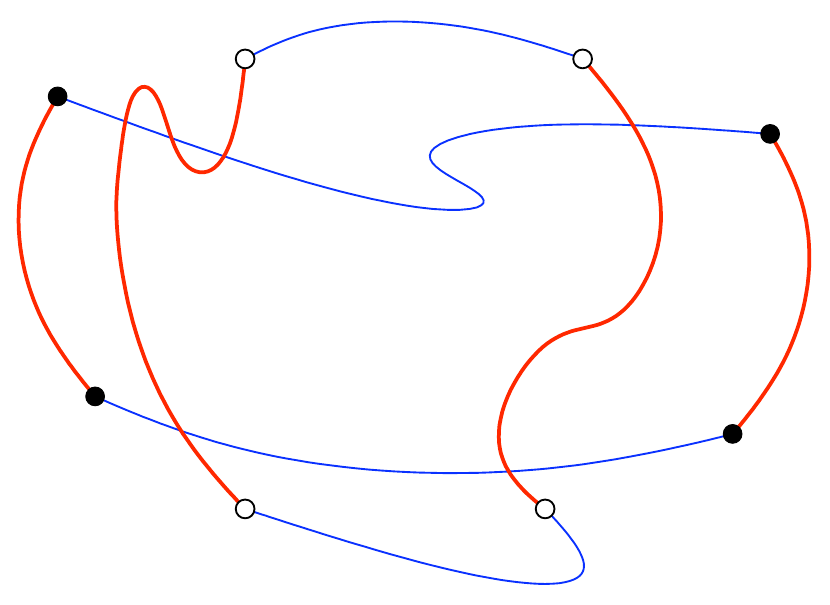}
}}
\caption{Two quads, $\Quad_1 \leq \Quad_2$.}
\label{f.poset}
\end{figure}

From the point of view of crossings, there is a natural partial order on $\QUAD_D$: we write $\Quad_1 \leq \Quad_2$ if any crossing of $\Quad_2$ contains a crossing of $\Quad_1$.  See Figure~\ref{f.poset}. Furthermore, we write $\Quad_1 < \Quad_2$ if there are open neighborhoods $\mathcal{N}_i$ of $\Quad_i$ (in the uniform metric) such that $ N_1\leq N_2$ holds for any $N_i\in \mathcal{N}_i$.  A subset $S\subset \QUAD_D$ is called {\bf hereditary} if whenever $\Quad\in S$ and $\Quad'\in\QUAD_D$ satisfies $\Quad' < \Quad$, we also have $\Quad'\in S$. The collection of all closed hereditary subsets of $\QUAD_D$ will be denoted by $\HH_D$. Any discrete percolation configuration $\omega_\eta$ of mesh $\eta>0$, considered as a union of the topologically closed percolation-wise open hexagons in the plane, naturally defines an element $S(\omega_\eta)$ of $\HH_D$: the set of all quads for which $\omega_\eta$ contains a crossing. Thus, in particular, critical percolation induces a probability measure on $\HH_D$, which will be denoted by $\P_\eta$. Note that since each quad is closed in $D$ and the elements of $\HH_D$ are closed hereditary sets, crossings of a quad $\Quad$ are allowed to use the boundary $\p\Quad$. 

Hereditary subsets can be thought of as Dedekind cuts in the setting of partially ordered sets (instead of totally ordered sets, as usual). It can be therefore hoped that by introducing a natural topology, $\HH_D$ can be made into a compact metric space. Indeed, let us consider the following subsets of $\HH_D$. For any $\Quad\in \QUAD_D$, let $\boxminus_\Quad:=\{S\in\HH_D:\Quad\in S\}$, and for any open $U\subset \QUAD_D$, let $\boxup_U:=\{S\in\HH_D: S\cap U=\emptyset\}$. Observe that it is natural to declare the sets $\boxup_U^c$ to be open, and (since the quads themselves are closed in $D$) the sets $\boxminus_\Quad$ to be closed. 
Therefore, we define $\T_D$ to be the minimal topology that contains every $\boxminus_\Quad^c$ and $\boxup_U^c$ as open sets. It is proved in \cite[Theorem 3.10]{\SchrammSmirnovNoise} that for any nonempty open $D$, the topological space $(\HH_D,\T_D)$ is compact, Hausdorff, and metrizable.
Furthermore, for any dense $\QUAD_0 \subset \QUAD_D$, the events $\{\boxminus_\Quad : \Quad\in\QUAD_0\}$ generate the Borel $\sigma$-field of $\HH_D$. Now, since Borel probability measures on a compact metric space are always tight, we have subsequential scaling limits of $\P_\eta$ on $\HH_\C$, as $\eta_k\to 0$, denoted by $\P_0=\P_0(\{\eta_k\})$. It is proved in \cite[Corollary 5.2]{\SchrammSmirnovNoise}  that, for any subsequential limit $\P_0$, and any finite collection of events $\boxminus_\Quad$ and $\boxminus_\Quad^c$, the joint probabilities converge along the subsequence; the main reason is that the boundary of the joint event can be shown to have zero measure under $\P_0$, so the weak convergence of the measures $\P_{\eta_k}$ to $\P_0$ is good enough.  

Of course, as explained carefully in \cite{\SchrammSmirnovNoise}, the choice of the space $\HH_D$ (or any other setup for the scaling limit) already poses restrictions on what events one can work with.
Note, for instance, that $\A:=\{ \exists$ neighborhood $U$ of the origin $0\in\C$ s.t.~all quads $\Quad \subset U$ are crossed$\}$ is clearly in the Borel $\sigma$-field of $(\HH_D,\T_D)$, and it is easy to see that $\P_0[\A]=0$ under any limit $\P_0$, but if the sequence of $\eta$-lattices is such that $0$ is always the center of a hexagonal tile, then $\P_\eta[\A]=1/2$. (The above machinery breaks down because the boundary of the event $\A$ has positive, or in fact, full measure under $\P_0$.) Similarly, one can talk about the Borel-measurable event $\mathcal{M}_\eta$ in $(\HH_D,\T_D)$ that counts whether there are more open than closed $\eta$-hexagons in the percolation configuration $\omega_\eta$ in $D$, but there is no good candidate for such an event that would be non-trivial under $\P_0$. (And there is good reason for this failure: it is proved in \cite[Section 4]{\BKS} that crossing events are asymptotically independent from majority events.)
Therefore, it is important to know that for a lot of natural events such problems do not occur. This will be the goal of Subsection~\ref{ss.limitarms} below.
\medskip

As mentioned above, $(\HH_D,\T_D)$ is a compact metrizable topological space. Unfortunately, we have not found a natural explicit metric that would yield this topology. However, in Subsection~\ref{ss.limitarms} below, it will be important to understand well what it means for two configurations to be close to each other. To this end, we will now define an explicit countable subbase of $\T_D$.

\begin{definition}\label{d.kquads}(A dyadic family of quads)
For any $k\geq 1$, let $(Q_{n}^{k})_{1\le n \le N_{k}}$ be the family of all quads which are {\it polygonal quads} in $D\cap 2^{-k} \Z^{2}$, i.e., their boundary $\p Q_n^k$ consists of edges of $D\cap 2^{-k} \Z^{2}$ and the four marked vertices are vertices of $D\cap 2^{-k} \Z^{2}$. (For fixed $k$, there are finitely many such quads since the domain $D$ is assumed to be bounded). We will denote by $\QUAD^k=\QUAD^k_D$ this family of quads. Notice that $\QUAD^k\subset \QUAD^{k+1}$. 
\end{definition}

Clearly, the family $\QUAD_\N:= \bigcup_k \QUAD^k$ is dense in the space of quads $(\QUAD_D,d_\QUAD)$. Thus, as mentioned above, the events $\{\boxminus_\Quad : \Quad\in \QUAD_\N\}$ generate the Borel $\sigma$-field of $\HH_D$. 


In order to use the open sets of the form $\boxup_U^c$, where $U$ is an open set of $\QUAD_D$, we will associate to each $Q\in \QUAD^k,\, (k\geq 1)$, the open set $\hat Q_k:= B_{d_\QUAD}(Q,2^{-k-10})$, and with a slight abuse of notation, we will  write $\boxup_{\hat Q}^c$ for the open set $\boxup_{\hat Q_k}^c$. 

\begin{figure}[!htp]
\begin{center}
\includegraphics[width=0.6\textwidth]{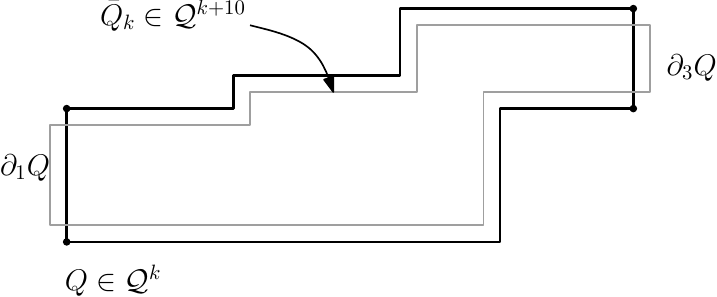}
\end{center}
\caption{A dyadic quad $Q\in \QUAD^k$ and the associated harder-to-cross quad $\bar Q_k$.}\label{f.barQ}
\end{figure}

Also, to each quad $Q\in \QUAD^k$, we associate the quad $\bar Q_k\in \QUAD^{k+10}$ which among all quads $Q'\in\QUAD^{k+10}$ satisfying $Q'>Q$ is the smallest one. Even though $>$ is not a total order, it is not hard to check that $\bar Q_k$ is uniquely defined. See Figure \ref{f.barQ} for an illustration. Note furthermore that $\bar Q_k$ satisfies $d_\QUAD(\bar Q_k, Q)\in [2^{-k-10},2^{-k-5}]$. Since, by definition $\bar Q_k > Q$, one has that $\boxminus_{Q}^c \subset \boxminus_{\bar Q_k}^c$. With a slight abuse of notation, we will write $\boxminus_{\bar Q}^c$ for the open set $\boxminus_{\bar Q_k}^c$.

\begin{definition}[A family of neighborhoods]
For each $k\geq 1$ and each point $\omega \in \HH$, let $\O_k(\omega)$ be the following open set:
\begin{equation}
\O_k(\omega):= \big( \bigcap_{Q\in \QUAD^k\,\text{s.t.} \, Q\notin \omega }\, \boxminus_{\bar Q}^c \big) \bigcap 
\big( \bigcap_{Q \in \QUAD^k \,\text{s.t.} \, Q\in \omega } \boxup_{\hat Q}^c\big)\,.
\end{equation}
Let also $\O_0(\omega)$ be the space $\HH$ for any $\omega\in \HH$.
\end{definition}

In words, these are the configurations $\omega'$ that are close to $\omega$ in the sense that for each fine dyadic quad that is crossed in $\omega$, there is a close-by dyadic quad crossed in $\omega'$, and for each fine dyadic quad that is not crossed in $\omega$, there is just a slightly harder-to-cross quad that is not crossed in $\omega'$.

\begin{remark}
Since for any $Q\in \QUAD^k$, $\boxminus_Q^c$ is already an open set, one might wonder why we have chosen here to relax $\boxminus_Q^c$ into $\boxminus_{\bar Q}^c$. This choice makes the definition and hence the forthcoming proofs more symmetric. 
\end{remark}

\begin{remark}
Let us point out that for any $\omega\in \HH$ and any $k\geq 0$, we have that $\omega\in \O_{k+1}(\omega) \subset \O_{k}(\omega)$. Hence the {\bf finite} coverings $\Cov^k:= \{  \O_k(\omega) :  \omega\in \HH \}$ are finer and finer as $k\to \infty$. 
\end{remark}




We will use this family of neighbourhoods through the following lemma.

\begin{lemma}\label{l.CovrtoCovk}
Let $d_\HH$ be any metric that generates the topology of $(\HH_D,\T_D)$. Then there exists a function 
$\r: \N_+ \longrightarrow \R_+$ such that, for any $k\geq 1$ and any $\omega,\omega'\in \HH$, if $d_\HH(\omega,\omega')\le \r(k)$, then $\omega'\in \O_k(\omega)$ and $\omega\in \O_k(\omega')$. 

In words, if $\omega$ and $\omega'$ are sufficiently close ($\r(k)$-close) to each other, then (up to small perturbations of $2^{-k-10}$) they have the same crossings for all quads in $\QUAD^k$.
\end{lemma}

\proof
Let us fix some integer $k\geq 0$.
Let $k'$ be a slightly larger integer, say $k+20$.  There are finitely many quads in $\QUAD^{k'}$, hence there are only finitely many possible open sets of the form $\O_{k'}(\omega)$, and the union of these covers $\HH$ (because $\omega\in \O_{k'}(\omega)$ for any $\omega\in \HH$). It follows that, to any point $\omega\in \HH$, one can associate a radius $r_\omega>0$ so that the ball $B_{d_\HH}(\omega,2 \,r_\omega)$  is included in at least one of the open sets $\O_{k'}(\bar\omega)$. 
Consider now the covering $\{ B_{d_\HH}(\omega, r_\omega) :  \omega\in \HH \}$, from which one can extract a finite covering 
\[
\{ B_{d_\HH}(\omega_i, r_i) : i=1,\ldots, N_{k'} \}\,.
\]
Let us define $\r(k):= \min_{1\le i \le N_{k'}} \{ r_i \}$ and check that it satisfies the desired properties.
Let $\omega,\omega'$ be any points in $\HH$ such that $d_\HH(\omega, \omega') \le \r(k)$. By our choice of $\r(k)$, one can find at least one ball $B_{d_\HH}(\omega_i, r_i)$ in the above covering such that both $\omega$ and $\omega'$ lie in the ball $B_{d_\HH}(\omega_i, 2\, r_i)$. In particular this means that one can find some $\bar \omega \in \HH$ such that both $\omega$ and $\omega'$ lie in $\O_{k'}(\bar \omega)$. Let us now prove that $\omega' \in \O_k(\omega)$; the other condition is proved similarly. 
Consider any quad $Q\in \QUAD^k$. We will distinguish the following cases:

\begin{itemize}
\item[(a)] Suppose $Q \in \bar\omega$ and $Q\in \omega$. Since $\omega'\in \O_{k'}(\bar \omega)$, we have that $\omega'\in \boxup_{\hat Q_{k'}}^c \subset \boxup_{\hat Q_k}^c$.

\item[(b)] Suppose $Q \in \bar\omega$ and $Q\notin \omega$. We need to show that $\bar Q_k \notin \omega'$. For this, note that one can find a quad $R$ in $\QUAD^{k'}$ which is such that $\bar Q_k > \bar R_{k'}$ and $\hat R_{k'}>Q$ (in the sense that all the quads in the open set $\hat R_{k'}$ are larger than $Q$). If $R$ happened to be in $\bar \omega$, then since $\omega\in \O_{k'}(\bar \omega)$, $Q$ would necessarily belong to $\omega$. Hence $R\notin \bar \omega$ and thus $\bar R_{k'}\notin \omega'$ which implies $\bar Q_k \notin \omega'$. 

\item[(c)] Suppose $Q\notin \bar \omega$ and $Q\in \omega$. We need to show that $\omega' \in \boxup_{\hat Q_k}^c$. Similarly to case (b), one can find a quad $R\in \QUAD^{k'}$ such that $Q>\bar R_{k'}$ and $\hat R_{k'} \subset \hat Q_k$. If $R$ was not in $\bar \omega$, then $Q$ would not be in $\omega$ either. Hence $R\in \bar \omega$ and thus $\omega'\in \boxup_{\hat R_{k'}}^c \subset  \boxup_{\hat Q_k}^c$.

\item[(d)] Finally, suppose $Q \notin \bar \omega$ and $Q\notin \omega$. Note that  $\bar Q_k > \bar Q_{k'}$. Since $\omega' \in \O_{k'}(\bar \omega)$, $\bar Q_{k'}$ is not in $\omega'$ and thus $\omega' \in  \boxminus_{\bar Q_k}^c$, which ends the proof of the lemma. \qed
\end{itemize}

\subsection{Uniqueness of the quad-crossing percolation limit}\label{ss.black}

In \cite{\SchrammSmirnovNoise}, only subsequential limits $\lim_{k\to\infty} \P_{\eta_k}=\P_0(\{\eta_k\})$ are considered.
As explained in the introduction, it is crucial for us to know that there is a unique and conformally invariant scaling limit in this quad-crossing space.
It turns out that this uniqueness property can be recovered from another notion of the scaling limit of percolation: 
the {\bf continuum non-simple loop ensemble} of \cite{\CamiaNewmanFull}. The point of view, or rather the topology considered in \cite{\CamiaNewmanFull} 
is quite different from the approach in \cite{\SchrammSmirnovNoise}. Instead of keeping track of the crossing events for quads, one keeps track 
of all the interface loops (that lie between some open cluster and some closed one). The loops are oriented depending 
whether their inside is open or not. To a discrete percolation configuration $\omega_\eta$ on $\eta \Tg$,
one thus associates a collection of oriented loops $\{\ell_i(\omega_\eta)\}_{i\in \N}$.

Let us briefly describe what is the topology of convergence in \cite{\CamiaNewmanFull}. First of all, the loops correspond here to continuous maps from the circle to the plane $\C$.
They are identified up to reparametrization by homeomorphisms of the circle with positive winding.
Since we want to consider percolation configurations in the whole plane, it is convenient
to compactify $\C$ as usual into $\C\cup \{\infty\} \simeq \mathbb{S}^2$. Let 
$d_{\mathbb{S}^2}$ be the induced metric on  $\C\cup \{\infty\}$.

We now equip the space $L$ of such loops with the following (pseudo-)distance
 \[
d_{L}( \ell ,\ell') := \inf_{\phi} \sup_{t \in \R / \Z}  d_{\mathbb{S}^2}\bigl(\ell(t) , \ell'(\phi(t))\bigr)\,,
\] 
where the infimum is taken over all homeomorphisms of the circle which have positive winding (recall that 
loops are oriented).

Now, let $\mathcal{L}$ be the space of countable collections of loops in $L$.
As mentioned above, any discrete percolation configuration $\omega_\eta$ gives rise to an element of $\mathcal{L}$.
We consider the Hausdorff topology on $\mathcal{L}$ induced by $d_L$. Namely if $c,c' \in \mathcal{L}$, let 
\[
d_{\mathcal{L}}(c,c') := \inf \{ \eps>0 : \forall \, \ell \in c \, \exists \, \ell'\in c' \text{ and } 
\forall \, \ell' \in c' \, \exists \, \ell\in c \text{ such that } 
d_{L}(\ell,\ell') \le \eps\}\,.
\]

In \cite{\CamiaNewmanFull}, the following theorem is shown:

\begin{theorem}[\cite{\CamiaNewmanFull}]
Let $\omega_\eta$ be a critical percolation configuration on the rescaled triangular grid $ \eta \Tg$ (on the whole plane).
The random configuration $\omega_\eta$ considered as a random element $\iota_\eta$ in  $\mathcal{L}$ converges in law, under the topology induced by $d_{\mathcal{L}}$, to a {\bf continuum percolation} $\iota \in \mathcal{L}$. 
\end{theorem}

Let $\P_0^{\mathcal{L}}$ denote the law of the above continuum percolation $\iota \in \mathcal{L}$. Several properties of this scaling limit $\iota$ are established; most importantly:

\begin{proposition}[{\bf ``finite chaining''}, \cite{\CamiaNewmanFull}, third item in Theorem 2]\label{pr.finitechaining}
Almost surely (under $\P_0^\mathcal{L}$), for any pair of loops $\ell,\ell'\in\iota$, there is a finite set a loops 
$\ell_0=\ell,\ell_1,\ldots, \ell_k=\ell'$ in $\iota$, such that for all $j\in [1,k]$, $\ell_{i-1}\cap \ell_i \neq \emptyset$.
\end{proposition}

We will not enter in more details here, but, using the fact that at the scaling limit there are no 6-arms events, it is easy to see that almost surely if $\ell,\ell'\in \iota$ touch each other and have the same orientation, then one loop cannot lie ``inside'' the other one. Conversely, if two loops of different orientations touch each other, then one has to be ``inside'' the other one. The meaning of ``inside'' needs some care here, since on a typical loop there are inside ``fjords'' and exterior ones. See \cite{\CamiaNewmanFull} for more details. This condition implies a parity constraint on the above finite chaining property.
\medskip

\ni{\bf Uniqueness of the scaling limit under the quad-crossing topology.}
To establish the uniqueness, it is enough to characterize uniquely the subsequential scaling limits $\P_0$. 
To achieve this, one has to prove that for any finite set of quads $\Quad_1,\Quad_2, \ldots,\Quad_p$, the {\bf joint law} of the crossing events 
$\boxminus_{\Quad_1}(\omega_\eta), \ldots, \boxminus_{\Quad_p}(\omega_\eta)$ converges in law to a unique random $p$-vector. We can restrict ourselves to 
piecewise smooth quads, since they form a dense family in $\QUAD_\C$. Notice that the 
$p=1$ case corresponds to Cardy's formula (and thus is not known on $\Z^2$). If we had at our disposal a joint Cardy's formula, we would be done.
For example, Watts' formula, proved in \cite{\DubedatWatts,\SchrammWatts}, is a special case of $p=2$. Such a general Cardy's fomula is not known, but fortunately we do not need 
an explicit formula here --- we just need to show that the joint laws converge. Let us fix our finite set of quads $\Quad_1,\Quad_2, \ldots,\Quad_p$, and consider 
percolation configurations $\omega_\eta$ with mesh $\eta\to 0$. Viewed as sets of loops, they converge to a continuum set of loops $\iota\sim \P_0^{\mathcal{L}}$. For any $j\in[1,p]$, as is well-known, the crossing event $\boxminus_{\Quad_j}(\omega_\eta)$ can be recovered from the discrete exploration path 
starting at one of the 4 corners of $\Quad_j$. It is known that such an exploration path converges towards an $\mathrm{SLE}_6$ curve in $\Quad_j$. Now the 
loops in $\iota$ are made of $\mathrm{SLE}_6$ curves, hence at first sight one seems to be in good shape. The trouble is that these $\mathrm{SLE}_6$ curves 
sit in the plane and do not see $\p \Quad_j$. Therefore, the continuum exploration path which would give us whether $\boxminus_{\Quad_j}$ is satisfied or not
is made of many pieces of $\mathrm{SLE}_6$ loops coming from $\iota$. Thus it is not a priori clear that $\boxminus_{\Quad_j}$ is measurable with respect to $\iota$.
The answer to this question is provided in \cite{\CamiaNewmanFull}; the proof (using the finite chaining lemma and the parity condition mentioned above) can be found right after Theorem 7 in \cite{\CamiaNewmanFull}:

\begin{proposition}
Let ($\Quad,a,b,c,d$) be a piece-wise smooth quad in the plane. Let $\omega_\eta$ be a sequence of percolation configurations which converge 
pointwise (for $d_\mathcal{L}$) towards $\iota \sim \P_0^\mathcal{L}$. (This can be done using Skorokhod representation theorem). The exploration paths
$\gamma_\eta$ from $a$ to $b$ converge pointwise (for the Hausdorff topology) to an $\mathrm{SLE}_6$ curve $\gamma$ from $a$ to $b$ which 
can be a.s.~recovered from $\iota$ in the sense that $\gamma=\gamma(\iota)$ can be constructed as a measurable function of $\iota$.
\end{proposition}

Hence, for each quad $\Quad_j$, simultaneously, one can extract from the continuum percolation $\iota\sim \P_0^\mathcal{L}$ an exploration path $\gamma_j$ which 
almost surely tells us about the crossing event in $\Quad_j$. (Only almost surely, since it could happen, for example, that $\gamma_j$ intersects one side at just a single point, in which case $\Quad_j$ could be crossed in some of the $\omega_\eta$ configurations and not crossed in others, so we would not get convergence of the crossing events, or could get the wrong limit. But this happens only with probability zero, due to the half-plane 3-arm event having exponent larger than 1.) This shows the uniqueness of the quad-crossing limit.
\medskip

\ni{\bf Conformal invariance of the scaling limit under the quad-crossing topology.}
The previous paragraph can be rephrased as follows: the continuum limit $\omega$ (under the quad topology \cite{\SchrammSmirnovNoise}) is a measurable 
function of the continuum limit $\iota$ (under the above topology \cite{\CamiaNewmanFull}). The conformal invariance of $\iota$ is provided by Theorem 7 in 
\cite{\CamiaNewmanFull}. In this section, we introduced the scaling limit $\iota$ only in the plane $\C$, but it can be defined for simply connected domains as well.
See \cite{\CamiaNewmanFull} for more details. If $\iota$ is a continuum percolation in some domain $\Omega$ and $f : \Omega \longrightarrow \tilde \Omega$ is a
conformal map, then $\tilde \iota:= f(\iota)$ (where each loop $\ell$ in $\iota$ is mapped conformally to a loop $f(\ell)$
in $\tilde \iota$) has the law of a continuum percolation in $\tilde \Omega$. This conformal invariance together with the above discussion imply that 
if $\omega$ is a scaling limit for $(\HH_\Omega,\T_\Omega)$, then $\tilde \omega := f(\omega)$ (defined by $\boxminus_{f(\Quad)}(\tilde \omega):= \boxminus_\Quad(\omega)$ for a dense family of quads $\Quad \subset \Omega$) has the law of a continuum percolation for $(\HH_{\tilde \Omega},\T_{\tilde \Omega})$. 
This type of invariance will be crucial in Section~\ref{s.covariance}.

\subsection{The limit of arm events}\label{ss.limitarms}

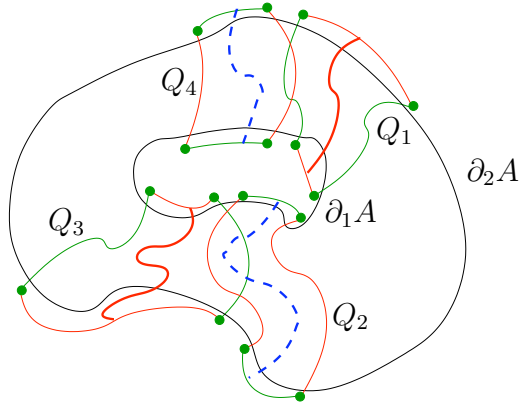
\begin{figure}[htbp]\label{f.def4arms}
\SetLabels
(.73*.45)$\p_1A$\\
(1.02*.55)$\p_2 A$\\
(.15*.43)$Q_3$\\
(.73*.2)$Q_2$\\
(.38*.77)$Q_4$\\
(.82*.65)$Q_1$\\
\endSetLabels
\centerline{
\AffixLabels{
\includegraphics[width=0.4\textwidth]{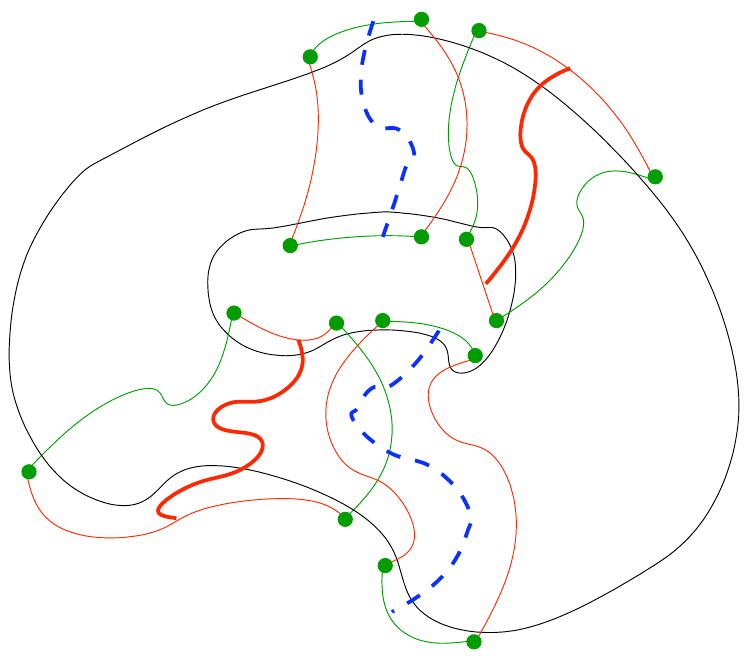}
}}
\caption{Defining the alternating 4-arm event using quads crossed or not crossed.}
\end{figure}

We need to define the usual arm events in a way that makes them measurable in the Borel $\sigma$-field of the quad-crossing topology $\T_D$ and as amenable as possible for proving the convergence of their probabilities as $\eta\to 0$. Namely, for any topological annulus $A\subset D$ with piecewise smooth inner and outer boundary pieces $\p_1 A$ and $\p_2 A$, we define the {\bf alternating 4-arm event} in $A$ as $\A_4=\bigcup_{\delta>0} \A_4^\delta$, where $\A_4^\delta$ is the existence of quads $\Quad_i \subset D$, $i=1,2,3,4$, with the following properties (see Figure~\ref{f.def4arms}):
\begin{itemize}
\item[(i)] $Q_1$ and $Q_3$ are disjoint and are at distance at least $\delta$ from each other; the same for $Q_2$ and $Q_4$; 
\item[(ii)] for $i\in \{1,3\}$, the sides $Q_i(\{0\}\times [0,1])$  lie inside $\p_1 A$ and the sides $Q_i(\{1\}\times [0,1])$ lie outside  $\p_2 A$; for $i\in \{2,4\}$, the sides $Q_i([0,1] \times \{0\})$ lie inside $\p_1 A$ and the sides $Q_i([0,1]\times \{1\})$ lie outside  $\p_2 A$; all these sides are at distance at least $\delta$ from the annulus $A$ and from the other $Q_j$'s;
\item[(iii)] the four quads are ordered cyclically around $A$ according to their indices;
\item[(iv)]  For $i\in \{1,3\}$, we have $\omega\in \boxminus_{Q_i}$, while for $i\in \{2,4\}$, we have $\omega\in  \boxminus_{Q_i}^c$. In plain words, the quads $Q_1,Q_3$ are crossed, while the quads $Q_2,Q_4$ are dual crossed between the boundary pieces of $A$ (with a margin $\delta$ of safety).
\end{itemize}

The definitions of general {\bf (mono-  or polychromatic) $k$-arm events} in $A$ are of course analogous: for arms of the same color we require the corresponding quads to be completely disjoint, and we still require all the boundary pieces lying outside the annulus $A$ to be disjoint.
 
\bl\label{l.meas4arm}
Let $A\subset D$ be a piecewise smooth topological annulus (with finitely many non-smooth boundary points). Then the 1-arm, any polychromatic 2- and 3-arm, the alternating 4-arm and any polychromatic 6-arm event in $A$, denoted by $\A_1$, $\A_2$, $\A_3$, $\A_4$ and $\A_6$, respectively, are measurable w.r.t.~the scaling limit of critical percolation in $D$, and $\lim_{\eta \to 0}\P_\eta[\A_i]=\P_0[\A_i]$. Moreover, in any coupling of the measures $\{\P_{\eta}\}$ and $\P_0$ on $(\HH_D,\T_D)$ in which $\omega_{\eta}\to\omega$ a.s.~as $\eta\to 0$,  we have 
\[
\Pb{\{ \omega_\eta\in \A_i \} \Delta \{\omega \in \A_i\}} \to 0\qquad (\text{as }\eta\to 0)\,.
\]
\el

\proof
Using a countable dense subset of $\QUAD_D$, it is clear that $\A_i$ is in the Borel $\sigma$-field of $(\HH_D,\T_D)$. In fact, it is also easy to see that it is an open set. It follows that, in any coupling $\omega_\eta\to\omega$, we have $\liminf_{\eta\to 0}\1_{\A_i}(\omega_\eta) \geq \1_{\A_i}(\omega)$ and hence $\Pb{\{\omega\in\A_i\}\setminus\{\omega_\eta\in\A_i\}}\to 0$.
The nontrivial direction is that 
\begin{equation}\label{e.setminus}
\Pb{\{\omega_\eta\in\A_i\}\setminus\{\omega \in\A_i\}}\to 0\,.
\end{equation}
The main step will be to prove that for each $i=1,2,3,4,6$ there exists a function $\eps_i:[0,1]\longrightarrow [0,1]$ with $\lim_{\delta\to 0} \eps_i(\delta)=0$ such that  
\begin{equation}\label{e.quant}
\P_\eta[\A_i^\delta \,|\, \A_i] > 1-\eps_i(\delta)\qquad \text{for all }\eta>0\,,
\end{equation}
so that we can detect the $i$ arms with quads satisfying the $\delta$-disjointness criteria with a high probability that does not depend on the mesh $\eta$.
\medskip

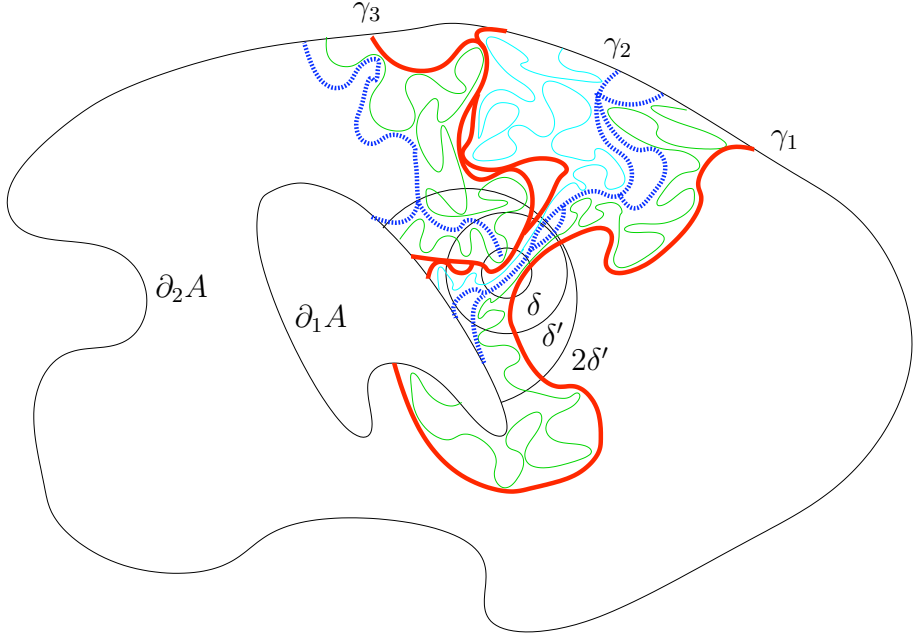
\begin{figure}[htbp]
\SetLabels
(.35*.5)$\p_1A$\\
(.2*.55)$\p_2A$\\
(.58*.52)$\delta$\\
(.6*.47)$\delta'$\\
(.64*.43)$2\delta'$\\
(.85*.80)$\gamma_1$\\
(.67*.95)$\gamma_2$\\
(.40*1.01)$\gamma_3$\\
\endSetLabels
\centerline{
\AffixLabels{
\includegraphics[height=3.2 in]{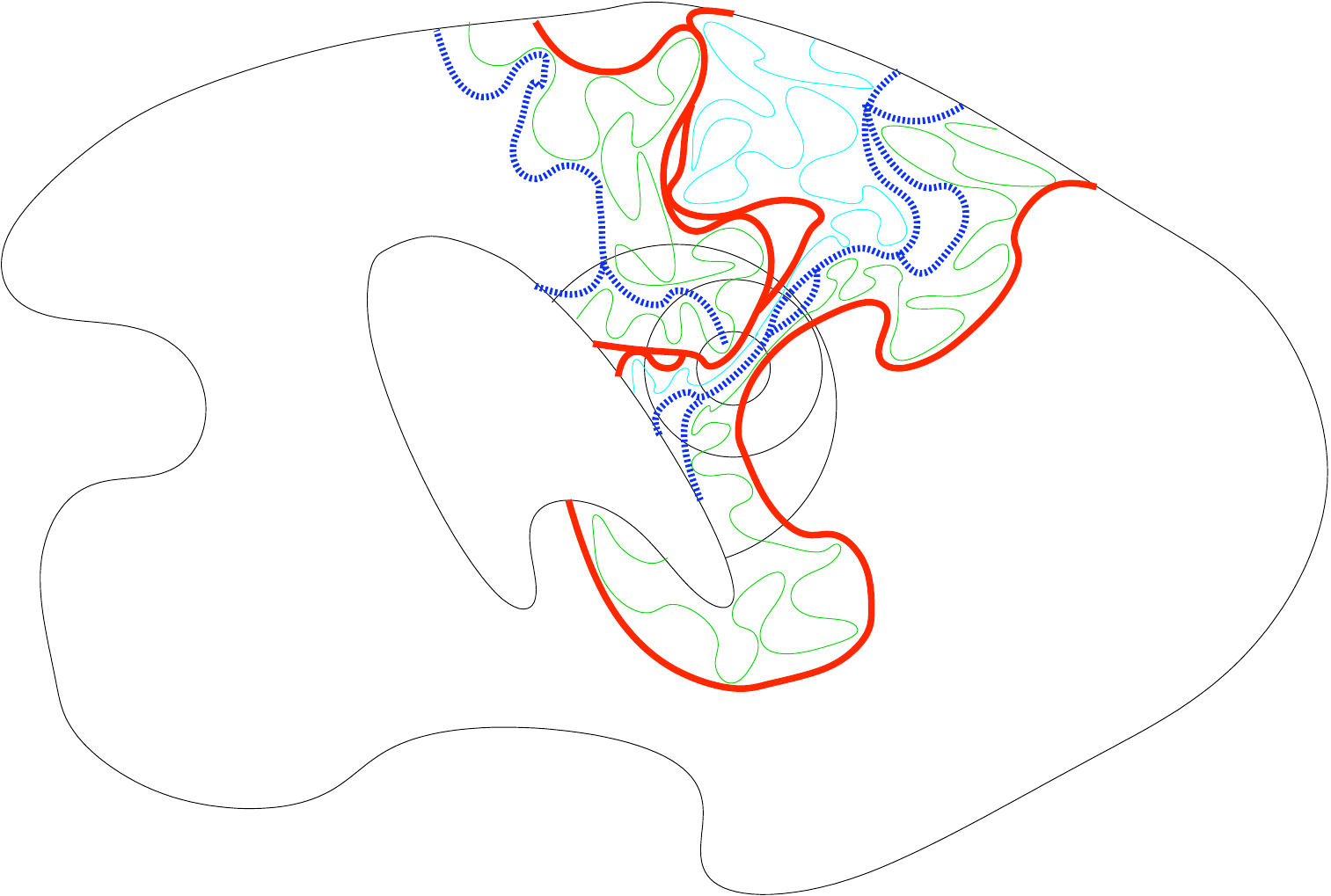}
}}
\caption{Convergence for $\A_4$. Closeness implies here a 7-arm whole plane event from radius $\delta$ to $\delta'$, and a 4-arm half-plane event from $2\delta'$ to a unit order.}
\label{f.474}
\end{figure}

Let us start with the $i=4$ case. If we have $\A_4$ in $\omega_\eta$, then choose 4 alternating arms arbitrarily,  consider their connected components (in their own color), and take the counterclockwisemost simple path in each component that connects $\p_1 A$ to $\p_2 A$, denoted by $\gamma_i$, $i=1,2,3,4$, in counterclockwise order. Clearly, if some $\gamma_i$ and $\gamma_{i+2}$ 
come $\delta$-close to each other somewhere inside $A$, at distance $\delta'$ from the closer boundary piece $\p_ i A$, then we have a polychromatic 7-arm event from that $\delta$-ball to distance $\delta'$,
and, assuming that the boundary is smooth here in a neighbourhood of macroscopic size, we have a half-plane 4-arm event from the $2\delta'$-neighborhood  of $\p_i A$ to a distance of unit order. See Figure~\ref{f.474}.
Around points of non-smoothness at the boundary, possibly with an inner angle as large as $2\pi$, we have to consider the 4-arm event in the ``half-plane squared'', giving the approximate square root of the half-plane probability on $\eta\Tg$ (by conformal invariance).
Using that the boundary is composed of finitely many smooth pieces, we get that the above $\delta$-closeness event happens with probability at most $O\big(\delta^{-2}\alpha_7(\delta,1) + \delta^{-1}\alpha_4^+(\delta,1)+\alpha_4^+(\delta,1)^{1/2+o(1)}\big)$ on $\eta\Tg$. 

Furthermore, if the endpoints of any two of our four arms come $\delta$-close on a smooth part of $\p_i A$, then we have a half-plane 3-arm event from a $\delta$-ball on that $\p_i A$ to a distance of unit order, or, at a non-smooth boundary point, a 3-arm event in the half-plane squared, in the worst case. So, this $\delta$-closeness event happens with probability at most $O\big(\delta^{-1}\alpha_3^+(\delta,1)+\alpha_3^+(\delta,1)^{1/2+o(1)}\big)$.

Finally, since we want our four quads to overhang the boundary pieces $\p_1A$ and $\p_2A$,  we need to continue our four arms to an additional distance $\delta$ at each end. This can be done unless we are stopped by a  crossing of the opposite colour, which, together with our arm, gives a 3-arm event at the boundary, so the probability is bounded again by $O\big(\delta^{-1}\alpha_3^+(\delta,1)+\alpha_3^+(\delta,1)^{1/2+o(1)}\big)$.

In summary, using the 7-arm probability in the plane (exponent $> 2$) and the 3-arm probabilities in the half-plane and the whole plane (exponent $=2$ and $>0$, respectively), we get that, with probability $1-O(\delta^\gamma)$, $\gamma>0$, none of the above $\delta$-closeness events happen, hence we can clearly find the four quads required, and we are done.

(Note that if some $\gamma_i$ and $\gamma_{i+1}$ come close to each other, then we get only a polychromatic 5-arm event, which still happens sometimes in the annulus. That is why we allow neighbouring quads to overlap in our definition of $\A_i^\delta$.)
\medskip

Assume now that $\A_3$ holds in $\omega_\eta$; say, there are two disjoint closed arms and one open. Viewed from the given  open arm $o$, take the clockwisemost and counterclockwisemost simple closed arms connecting $\p_1 A$ to $\p_2 A$, denoted by $c_1$ and $c_2$. If we take any $x_i$ on each $c_i$, then there is a open path $o'=o'(x_1,x_2)$ connecting them. If $c_1$ and $c_2$ come $\delta$-close to each other without sandwiching $o$, with some points $x_i$ on $c_i$, then we get a 6-arm event around the pair $x_1,x_2$: four closed arms given by $c_1$ and $c_2$, and two open arms given by $o'$. If $c_1$ and $c_2$ come $\delta$-close to each other while sandwiching $o$, then we have 6 arms given by $c_1$, $c_2$ and $o$. Of course, if the points of closeness are close to either boundary piece $\p_i A$, then we have to consider half-plane 3-arm events, as well. Altogether, we get that,  with good probability, none of these $\delta$-closeness events happen, moreover, just as for $\A_4$, we can extend the crossings to additional distances $\delta$ at the ends, and we are done.
\medskip

The case of $\A_6$ clearly follows from the combination of the above arguments for $\A_3$ and $\A_4$.
\medskip

\begin{figure}[htbp]
\SetLabels
(.65*.38)$\p_1A$\\
(.92*.88)$\p_2 A$\\
(.43*.52)$\gamma_1$\\
(.53*.50)$\gamma_2$\\
(.38*.93)$\gamma_1^c$\\
(.7*.88)$\gamma_2^c$\\
\endSetLabels
\centerline{
\AffixLabels{
\includegraphics[height=3.2 in]{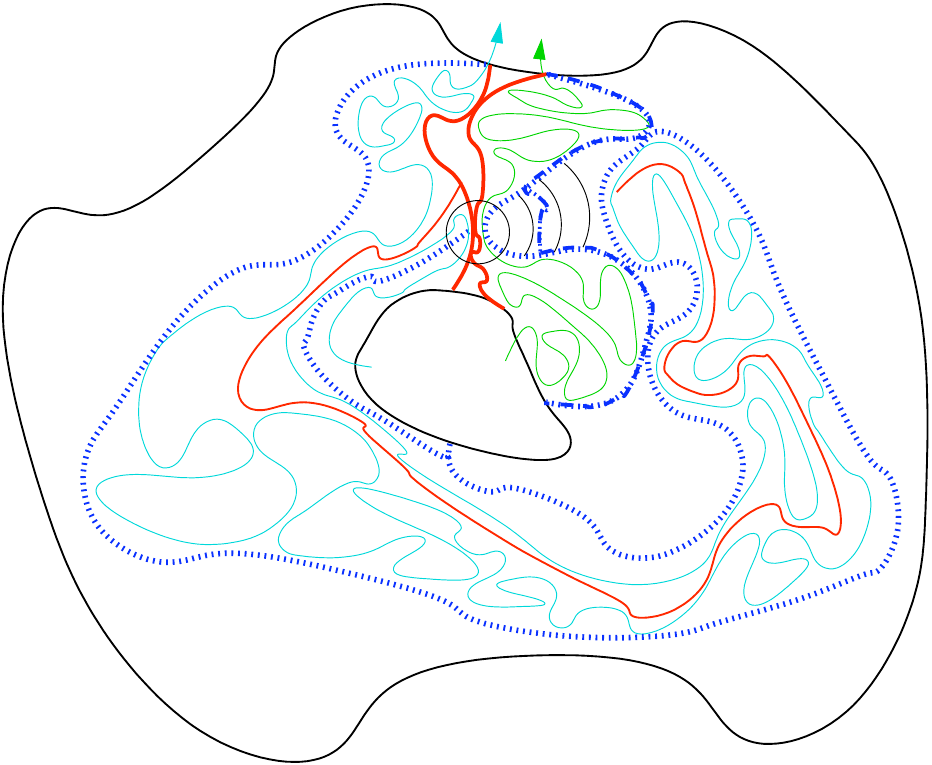}
}}
\caption{Convergence for $\A_2$. How to avoid $\delta$-closeness despite a 5-arm event around a $\delta$-ball produced by $\gamma_1^o$ and $\gamma_2^c$.}
\label{f.2arm}
\end{figure}
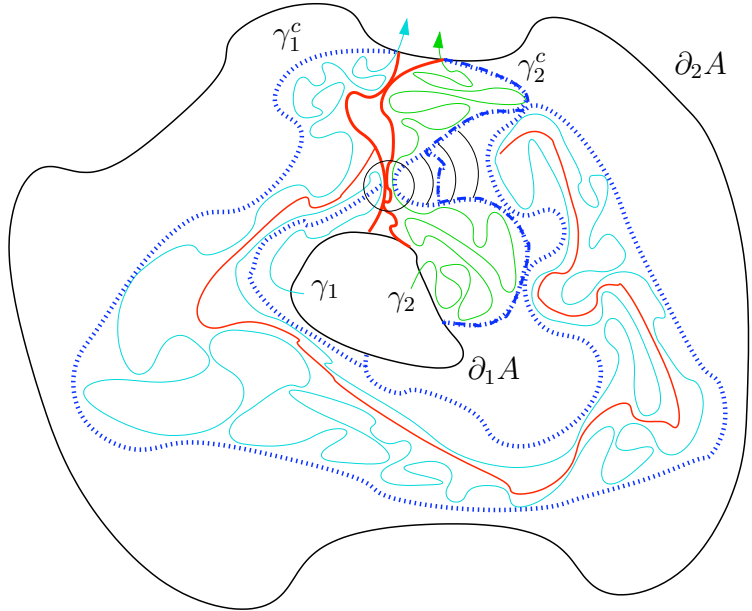

Now assume that $\A_2$ holds; then there are at least two interfaces from $\p_1 A$ to $\p_2 A$. Condition on the set of all the interfaces: its cardinality is a positive even number. If there are at least 4 interfaces, then also $\A_4$ holds, so, by our above result, we can detect even four arms, thus did not lose any probability in the scaling limit. Therefore, we may assume that there are exactly two interfaces, denoted by $\gamma_1$ and $\gamma_2$. Let the right (clockwise) side of $\gamma_1$ and the left (counterclockwise) side of $\gamma_2$ be the open sides. Denote the outer boundaries of the two $\gamma_i$'s, two open and two closed simple paths, by $\gamma_1^c,\gamma_1^o,\gamma_2^o,\gamma_2^c$, in the clockwise order. Of course, $\gamma_1^o$ and $\gamma_2^o$ do not have to be disjoint from each other at all, and the same for the $\gamma_i^c$'s. Consider now $\gamma_1^o$ and $\gamma_2^c$. It is possible that these two paths come very close to each other by sandwiching $\gamma_2$: away from the boundaries $\p_i A$ this requires a 5-arm event, which has exponent 2. More precisely, the expected number of 5-arm events from a grid of $\delta$-boxes to unit distance is finite and positive, uniformly in the mesh $\eta$, see \cite[First exercise sheet]{\WWperc}; therefore, the number of these 5-arm $\delta$-closeness events (for any fixed $\delta$) is tight in $\eta$. Denote these $\delta$-balls by $B_1,\dots, B_t$, $t=t(\delta,\eta)$. We want to find a closed arm that goes around these balls. For this, consider $\gamma_1^c$: it is unlikely that it comes very close to any of these balls (in particular, that it intersects any of them), since that would produce already a 6-arm event by bringing in an extra open arm (coming from the open side of $\gamma_1$). So, around any of the $B_j$, there is a region of ``macroscopic radius'' that is bordered by the two closed arms given by $\gamma_2^c$ and by $\gamma_1^c$. In this region between $\gamma_1^c$ and $\gamma_2^c$ that does not contain the interfaces $\gamma_i$, all the information that we have about the percolation configuration is that there is no open arm from $\p_1 A$ to $\p_2 A$. This is a monotone increasing event for ``closed'', hence, by the FKG inequality, we can apply RSW in consecutive disjoint annuli around $B_j$ to produce a closed shortcut between the two closed arms given by $\gamma_1^c$. Since the number $t(\delta,\eta)$ of the 5-arm balls that we need to go around is tight in $\eta$, with large probability we can construct a closed arm that does not come $\delta$-close to $\gamma_1^o$, anywhere far from the boundaries $\p_i A$. As usual, boundary closeness can easily be excluded using the 3-arm half-plane event.
\medskip

\begin{figure}[htbp]
\SetLabels
(.23*.27)$\p_1A$\\
(.43*.75)$\p_2 A$\\
(.84*.73)$\delta$\\
(.81*.52)$F$\\
\endSetLabels
\centerline{
\AffixLabels{
\includegraphics[height=2.5 in]{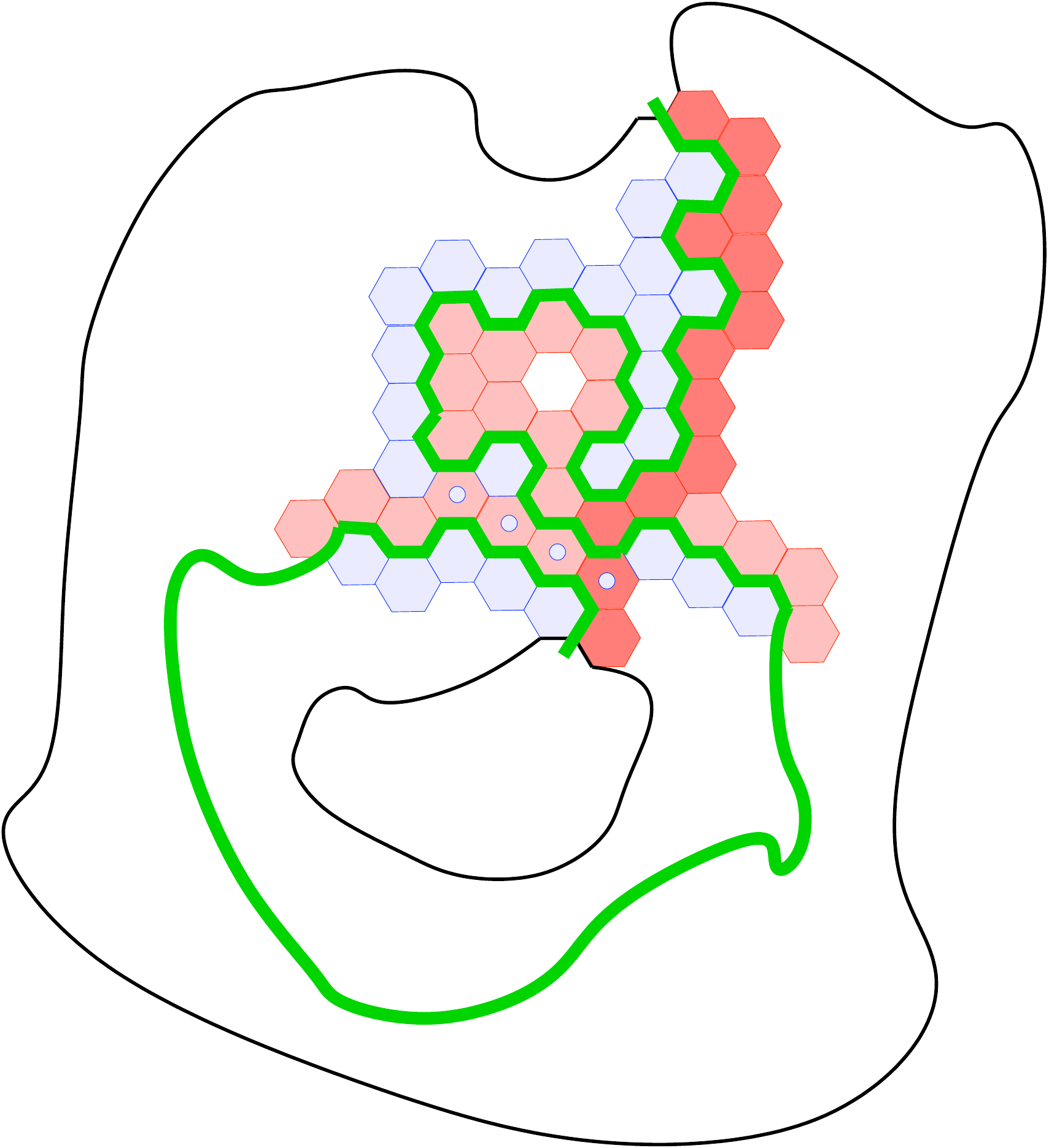}
\hskip 0.5 cm
\includegraphics[height=2.5 in]{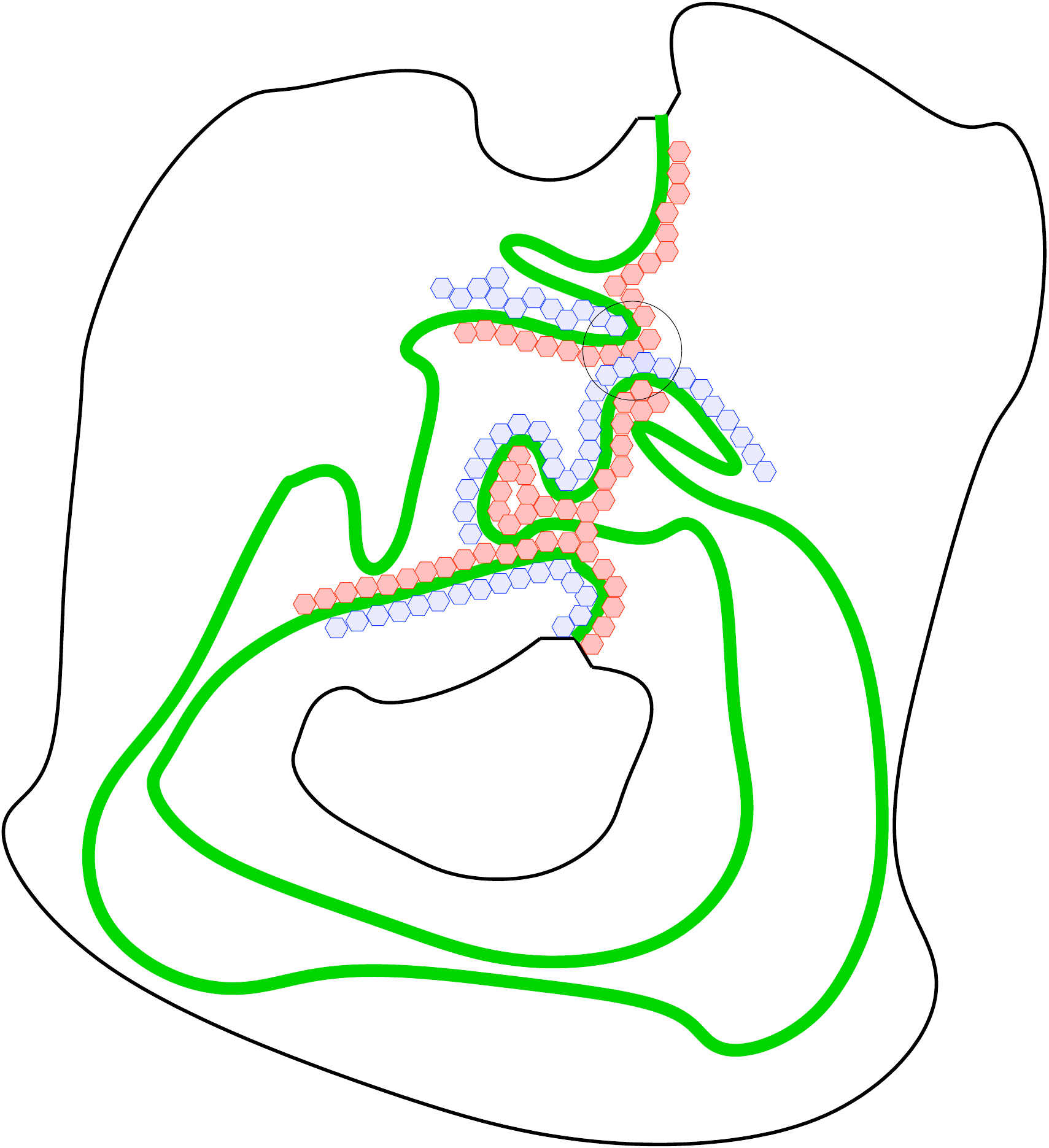}
}}
\caption{Convergence for $\A_1$. The left picture shows how to find an open (red) simple arm using the radial exploration process. The right picture shows a typical way how a $\delta$-almost self-touch of that open simple arm would imply a 6-arm event from radius $\delta$.}
\label{f.1arm}
\end{figure}

For $\A_1$, there is a new difficulty compared to the previous cases: if there is no closed crossing in $A$, i.e., there is an open circuit, then things like the ``clockwisemost simple path in the component of the open crossing'' are not well-defined. So, let us take a {\bf radial exploration interface}: start from the inner boundary $\p_1 A$, with open hexagons on the right, closed hexagons on the left of the interface. If the interface makes a clockwise loop around $\p_1 A$, then we have discovered a closed circuit, so there is no open crossing. Whenever the interface makes a counterclockwise loop around $\p_1 A$, let us pretend that the previously visited open hexagon which the interface has just bumped into was closed, and continue the exploration. After the first such forced turn, the interface can bump into itself not only because of counterclockwise loops around $\p_1 A$, but also because of the recolored hexagons; nevertheless, the rule for recoloring hexagons remains the same. This interface reaches $\p_2A$ if and only if there is an open crossing. Given this event, erase chronologically the  counterclockwise loops around $\p_1 A$ it has completed. The resulting interface has an open (non-simple) crossing path on its right, without loops around $\p_1 A$. Taking the right (clockwise) boundary of it, we get a simple open path. See the left picture on Figure~\ref{f.1arm}. We now claim that the probability that this path comes to $\delta$-close to itself (after a macroscopic detour) goes to 0 as $\delta\to 0$, uniformly in $\eta$, hence it does not happen in the scaling limit. In particular, our open simple path does not make loops around $\p_1 A$ in the scaling limit, and hence we can detect it with the event of crossing a quad.

We still need to prove the claim. If our simple open path was the outer boundary of a chordal exploration path, then coming $\delta$-close would produce a 6-arm event from radius $\delta$ to order 1, and we would be done. Now the exploration path can also go around $\p_1 A$, and come $\delta$-close to itself from the other side. This would again cause a 6-arm event. See the right picture of Figure~\ref{f.1arm}. However, because of the forced turns (denoted by $F$ on the same picture), some of the closed arms among those 6 arms might be faulty: open hexagons appear at the forced turns, usually building large pieces of mixed color arms that we pretended to be closed. (For instance, on the left picture of Figure~\ref{f.1arm}, we have four faulty hexagons next to each other.) What saves us is that the first time a closed arm emanating from a ball of $\delta$-closeness reaches a faulty hexagon (on either side of a faulty piece of an arm), there are at least two {\it additional} genuine arms emanating from that neighbourhood (since we have to produce the faulty hexagon, by the radial exploration path touching itself). Combining with the fact that the number of macroscopically different faulty places is tight as $\eta\to 0$ (because they imply 5-arm events), one gets that these faulty places do not help noticeably our open path to come $\delta$-close to itself. This finishes the proof of~(\ref{e.quant}).
\medskip

The point of (\ref{e.quant}) is that if $\A_i^\delta(\omega)$ holds for some $\omega\in\HH$ and $\delta>0$, and if we choose $k$ such that $2^{-k}<\delta/2^{10}$, then any $\omega'\in\O_k(\omega)$ will clearly satisfy $\A_i$, as well. Thus,~(\ref{e.quant}) implies that, for all $\eta>0$,
\begin{equation}\label{e.quantt}
\P_\eta\big[ \A_i(\omega') \text{ for all }\omega'\in\O_k(\omega) \,\big|\, \A_i(\omega) \big] > 1-\eps_i(\delta)\,.
\end{equation}
Now, take a coupling in which $\omega_\eta\to\omega$ almost surely. By Lemma~\ref{l.CovrtoCovk}, for $\eta=\eta(k,\delta)$ small enough, we have $\Pb{\omega \in \O_k(\omega_\eta)}>1-\eps_i(\delta)$. Therefore,
\begin{align*}
\Pb{\{\omega_\eta\in\A_i\}\setminus\{\omega \in\A_i\}} &\leq \Pb{\omega_\eta\in\A_i,\ \omega\in\O_k(\omega_\eta),\ \omega \not\in\A_i } + \eps_i(\delta)\\
&\leq 2\, \eps_i(\delta)\,,\qquad \text{using }(\ref{e.quantt})\,.
\end{align*}
By taking $\delta\to 0$, we obtain~(\ref{e.setminus}), finishing the proof of the lemma.
\QED


For future reference, let us reiterate that the proof of Lemma~\ref{l.meas4arm} implies the following quantitative bound on how close the event $\A^\delta_4$ is to the event  $\A_4$. 

\begin{corollary}
On $\eta\Tg$, for any piecewise smooth annulus $A$ and any $\gamma<1$, there is a constant $c=c_{A,\gamma}>0$ such that, 
for any $\delta>0$:
\[
\Pb{\A_4^\delta \md \A_4} \geq 1- c\, \delta^\gamma\,.
\]
It is easy to see that the same result holds on $\Z^2$, as well, only with some undetermined exponent $\gamma>0$ instead of it being arbitrarily close to 1. 
\end{corollary}

\begin{remark}\label{r.quadarm}
The above proofs clearly apply also to the case when the arms have to end at prescribed arcs on one or both of the boundaries $\p_i A$.
\end{remark}

\begin{remark}\label{r.joint4arm}
As mentioned above, \cite[Corollary 5.2]{\SchrammSmirnovNoise} proved the $k\to\infty$ convergence of the $\P_{\eta_k}$-probability of the intersection of any finite set of quad crossing and non-crossing events, where $\{\eta_k\}$ is a subsequence for which the limit $\P_0(\{\eta\})$ in $(\HH_D,\T_D)$ exists. Therefore, the proof of the above lemma, together with the uniqueness of the quad-crossing limit, imply the convergence of the joint $\P_{\eta}$-distribution of any finite set of macroscopic arm events, and the convergence in probability of the joint indicator in any coupling where $\omega_\eta\to\omega$ as $\eta\to 0$; the main point is that, in the scaling limit, the 6-arm and half-plane 3-arm events do not occur anywhere in the finite set of annuli, hence the convergence of quad crossings and non-crossings really measures what we want. 
\end{remark}
\medskip

Expressing arm events using the quad-crossing topology $\T$ is a natural first step towards proving that other notions of the percolation scaling limit are measurable w.r.t.~$\T$. For instance, we immediately get the following partial result towards extracting an $\mathrm{SLE}_6$ from the quad-crossing scaling limit. 

\begin{corollary} \label{cor.SLE6quad}
In a smooth domain $\Omega$ with boundary points $a,b\in\Omega$, the joint law $(\omega_\eta,\gamma_\eta)$ of the discrete quad-crossing configuration $\omega_\eta$ and the exploration path trace $\gamma_\eta$ from $a_\eta$ to $b_\eta$ (where $a_\eta\to a$ and $b_\eta\to b$) converges to a coupling $(\omega,\gamma)$, where $\omega$ is a percolation  quad-crossing scaling limit and $\gamma$ is the {\bf trace of a chordal $\mathrm{SLE}_6$} from $a$ to $b$ that is a function of $\omega$. The topology for the convergence of $\gamma_\eta$ to $\gamma$ is the Hausdorff distance of closed sets in $\Omega$.  
\end{corollary}

\proof 
Note that for any $\eps>0$, a closed $\eps$-box $Q$ inside $\Omega$ is touched by the exploration path $\gamma_\eta$ if{f} there are two arms from $\p Q$ to the boundary arcs $ab$ and $ba$, in the appropriate colors. Consider any finite covering of $\Omega$ by $\eps$-boxes. By Remarks~\ref{r.quadarm} and~\ref{r.joint4arm} applied to the 2-arm events from these small $\eps$-boxes to the smooth boundary arcs $ab$ and $ba$, and by the convergence in law of the exploration path $\gamma_\eta$ to a chordal $\mathrm{SLE}_6$ curve $\gamma$ (see Subsection~\ref{ss.interface} below), the set of $\eps$-boxes touched by $\gamma_\eta$ converges a.s., in any coupling with $\gamma_\eta\to\gamma$, to the set of $\eps$-boxes touched by $\gamma$, and this set is measurable w.r.t.~the quad-crossing topology $\T$. Taking a decreasing intersection of these compact sets as $\eps\to 0$, we get the measurability of the trace of the $\mathrm{SLE}_6$ path $\gamma$.
\QED

However, this convergence of the closed sets $\gamma_\eta$ to $\gamma$ is weaker than what one would really want for curves that naturally have their points ordered. ``Cutting'' the closed set $\gamma$ into a continuous increasing family in such a way that it gives the right curve would require work. So, the following question remains. See Subsection~\ref{ss.interface} for some related discussion.

\begin{question} \label{q.SLEmeas}
Is it the case that the chordal $\mathrm{SLE}_6$ {\bf curve} is measurable with respect to the quad-crossing continuum percolation $\omega$ inside $(\Omega, a, b)$?
\end{question}

Assuming a positive answer to Question~\ref{q.SLEmeas}, one could also try to do the measurability of the Camia-Newman description of the full scaling limit using iterated chordal $\mathrm{SLE}_6$ paths \cite{\CamiaNewmanFull}, but there is an additional issue here: the iterative process requires running the $\mathrm{SLE}_6$ paths in domains with fractal boundary, while we required smooth boundaries above.

\section{Coupling argument}\label{s.coupling}

Let us consider some annulus $A=A(u,v)$ of inner radius $u$ and outer radius $v>u$. Let $\omega_{\eta}$ be a percolation
configuration on the triangular grid of mesh $\eta$ inside $A$. Call $\Gamma$ the set of percolation interfaces which cross $A$
(it might be that there are no such interfaces, in which case $\Gamma=\emptyset$). We will need to measure how well separated
the interfaces are on the inside boundary $\p_1 A$. For that purpose we define a measure of the {\bf (interior) quality} $\Qual(\omega_{\eta})
=\Qual(\Gamma)$ to be the least distance between the endpoints of $\Gamma$ on $\p_1 A$ normalized by $u$. More precisely, if there are $p\geq 2$
interfaces crossing $A$ and if $x_1,\ldots,x_p$ denote the endpoints of these interfaces on $\p_1 A$, then we define
$$
\Qual(\omega_{\eta}) = \Qual(\Gamma)= \frac{1}{u} \inf_{k\neq l} |x_k - x_l |,
$$
where $|\cdot|$ denotes Euclidean distance. If $\Gamma=\emptyset$, we define $\Qual(\Gamma)$ to be zero.

We define similarly the {\bf exterior quality} $\Qual^+(\omega_{\eta}) = \Qual^+(\Gamma)$ to be the least distance between
the endpoints of $\Gamma$ on $\p_2 A$ normalized by $v$ (and set to be zero if there are no such interfaces).
  For any $\alpha>0$, let $\Sep^{\alpha}(u,v)$ be the event that
 the quality $\Qual(\Gamma)$ is bigger than $\alpha$; also, $\Sep_+^{\alpha}(u,v)$ will denote the event that $\Qual^+(\Gamma) > \alpha$.

  For a square $B_R$ of radius $R$, we define a notion of {\bf faces} around that square.
  Let $x_1,\ldots,x_4$ be four distinct points on $\p B_R$ chosen in a counterclockwise order.
  We will adopt here cyclic notation, i.e., for any $j\in \Z$, we have $x_j=x_i$ if $j\equiv i  [4]$.
  For any $i\in\Z$, let $\theta_i$ be a simple path of hexagons joining $x_i$ to $x_{i+1}$, i.e., a sequence of hexagons $h_1,\ldots,h_n$ such that $h_i$ and $h_j$ are neighbors if and only if $|i-j|=1$. We assume furthermore
  that there are no hexagons in $\theta_i$ which are entirely contained in the square of radius $R$ (they might still
  intersect $\p B_R$) and that all the hexagons in $\theta_i$ are white (=open) if $i$ is odd and black (=closed) if $i$ is even (the endpoints $x_i$ could be of either colour). If a set of paths
  $\Theta=\{\theta_1,\ldots,\theta_4 \}$ satisfies the above conditions, $\Theta$ will be called a configuration of {\bf faces}
 with endpoints $x_1,\ldots,x_4$. We define similarly the quality of a configuration of faces $\Qual(\Theta)$ to be the least distance between the endpoints of the faces, normalized by $R$.

  Let $\mathcal{G}(u,v)$ be the event that there exist 4 alternating arms from radius $u$ to $v$ and there is no extra disjoint arm crossing $A$. In particular, on the event $\mathcal{G}(u,v)$, the set of interfaces $\Gamma$ consists of exactly 4 interfaces, and furthermore, any
 two consecutive interfaces have to share at least one hexagon
 (if not, there would be at least 5 arms from $u$ to $v$).
  It is easy to see that on the event $\mathcal{G}(u,v)$, the 4 interfaces of $\Gamma$ induce a natural
  configuration of faces $\Theta=\Theta(\Gamma)$ at radius $u$ with same endpoints $x_1,\ldots,x_4$ as $\Gamma$.
  More precisely,
  if $\mathcal{H}$ is the set of all the hexagons neighboring the 4 interfaces in $\Gamma$,
  then by the definition of $\mathcal{G}(u,v)$, the connected component of $\C \setminus \mathcal{H}$ which contains
  the center of the annulus $A$ is a bounded domain; the set of hexagons which lie on the boundary of this domain
  form the 4 faces of $\Theta$. 
  In what follows, when we condition on the event $\mathcal{G}(u,v)$, 
  we will often condition also on the configuration of faces $\Theta$.  
  (Notice that, by definition, on the event $\mathcal{G}(u,v)$
  we have $\Qual(\Gamma) = \Qual(\Theta)$).

 If we are given a square of radius $R$ and a configuration of faces $\Theta=\{\theta_1,\ldots,\theta_4\}$ around
 that square, $\mathcal{D}=\mathcal{D}_{\Theta}$ will denote the bounded component of $\C\setminus \Theta$
 (this is a finite set of $\eta$-hexagons). Let $U=U_{\Theta}$ be the random variable which
 is set to be one if there is an open crossing from $\theta_1$ to $\theta_3$ inside $\mathcal{D}_{\Theta}$ and 0 otherwise.
 For any radius $r<R$, let $\mathcal{A}_{\Theta}(r,R)$ be the event that there are open arms from the box of radius $r$ to the open
 faces $\theta_1$ and $\theta_3$ and there are closed arms from $r$ to the closed faces $\theta_2,\theta_4$ (in other words, if one starts
 4 interfaces in $\mathcal{D}_{\Theta}$ at the endpoints of $\Theta$ at radius $R$, then these interfaces go all the way to radius $r$).

 Let $Q$ be the square of radius $1$ around the origin. As a particular case of the above, for any $\eta<r<1$, we define $\mathcal{A}_0(r,1)$
 to be the event that there are open arms from the square of radius $r$ to the left and right edges of $Q$ and there are closed arms from
 $r$ to the top and bottom edges of $Q$. Similarly, $U_0$ will be the indicator function of a left-right crossing in $Q$.

\begin{figure}[htbp]
\centerline{
\includegraphics[height= 2.6 in]{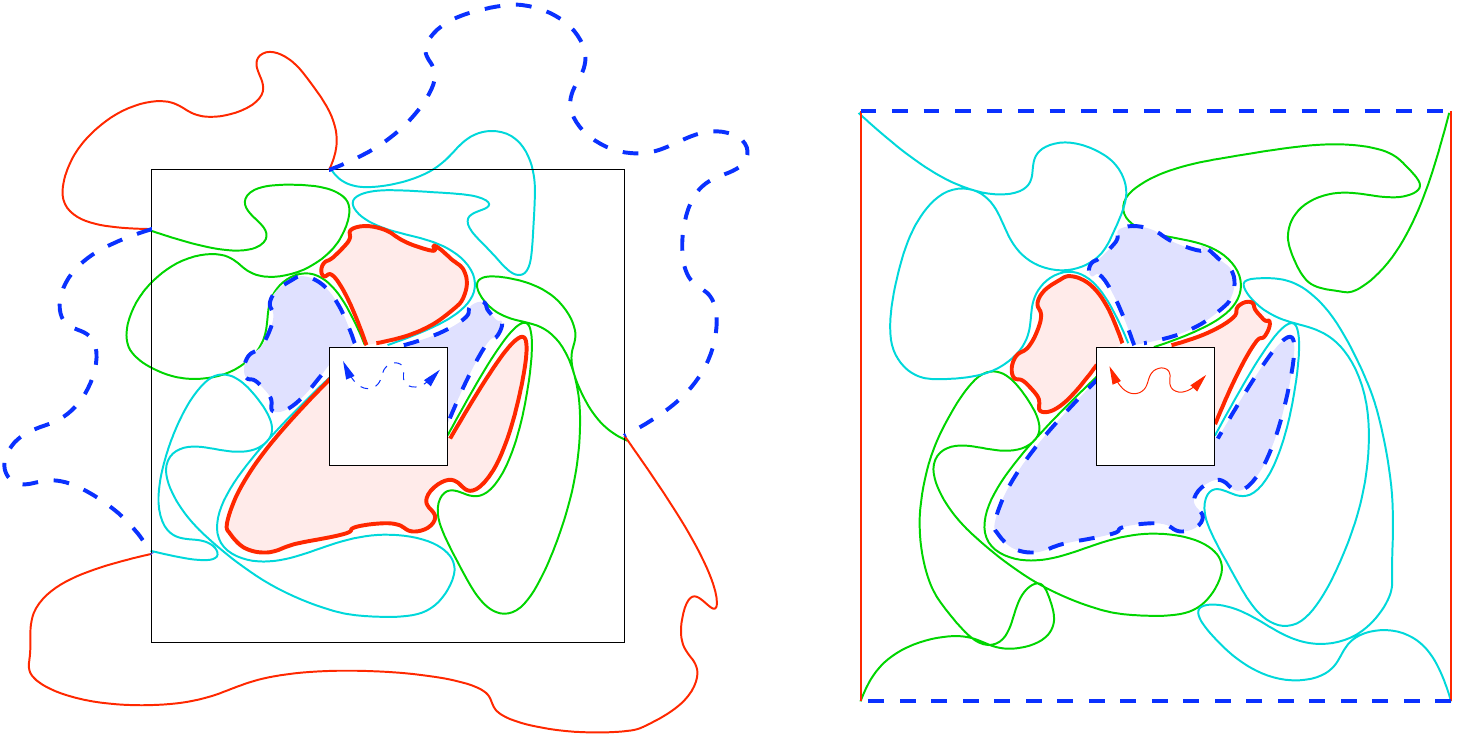}
}
\caption{Successful coupling of $\nu_\Theta$ and $\nu_0$ when $U_\Theta=\tau=0$: the two four-tuples of interfaces induce the same configuration of faces at radius $r$, only with reversed colors.}
\label{f.coupling}
\end{figure}

 \begin{proposition}[Coupling property]\label{pr.couplingfaces}
 Let $\tau \in\{0,1\}$, and let $A=A(r,R)$ be the annulus centered around 0 with radii  $100 \eta<  10 r <  R \leq 1$.
 Assume we are given any configuration of faces $\Theta = \{\theta_1,\ldots,\theta_4\}$ around the square of radius $R$.
 Let $\nu_{\Theta}$ be the law of the percolation configuration inside $\mathcal{D}_{\Theta}$ conditioned on the events
 $\mathcal{A}_{\Theta}(r,R)$ and $\{U_{\Theta} = \tau \}$.

 Let $\nu_0$ be the law of the percolation configuration in the square $Q$ under the conditioning
 that $\mathcal{A}_0(r,1)$ and $\{U_0 = 1\}$ hold.

There exists an absolute exponent $k>0$ such that the following holds.
 If $\tau=1$, there is a {\bf coupling} of the conditional laws $\nu_{\Theta}$ and $\nu_0$ so that with (conditional)
 probability at least $1 - (r/R)^{k}$,
 the event $\mathcal{G}(r,R)$ is satisfied for both configurations, and the induced faces at radius $r$, $\Theta(r)$ and $\Theta_0(r)$, are identical. If $\tau=0$, there is also such a coupling except that the faces $\Theta(r)$ and $\Theta_0(r)$ are (with probability at least $1- (r/R)^{k}$)
 identical but with reversed color (see Figure~\ref{f.coupling}).
 \end{proposition}

\begin{remark}
The proposition also holds if the annulus $A$ is not centered around the origin, but there are issues here coming
from the discrete lattice: indeed, the set of $\eta$-hexagons intersecting a square $S$ is not invariant under translations of the square $S$.
We will deal with this issue in the next section. 
\end{remark}

The proof of this coupling property will rely on the very useful {\it separation of interfaces phenomenon}, which has several versions in the literature, for several different planar models (percolation, random walks, Gaussian free field, loop-erased random walk). The following percolation lemma is essentially from \cite{\Sapozhnikov}, but, for completeness, the proof of this exact version is included in Appendix~\ref{s.appendix} below, together with a discussion of the general phenomenon.

\begin{lemma}[Strong Separation Lemma]\label{l.separation}
Let $A=A(R_1,R_2)$ be the annulus centered around 0 with radii $100 \eta < 10r < 2 R_1 < R_2$. Assume we are given some 
configuration of faces $\Theta=\{\theta_1,\ldots,\theta_4\}$ around the square of radius $R_2$. Note that the initial quality $\Qual(\Theta)$ might be arbitrarily small
(it could be of order $\eta$, for example). 

Let $\nu_{\Theta}$ be the conditional law of the percolation configuration inside $\mathcal{D}_\Theta$ 
conditioned on the event $\mathcal{A}_\Theta(r, R_2)$, and let $\Gamma = \Gamma(R_1)$ be the four-tuple of interfaces in $\mathcal{D}_\Theta$ which start at the four endpoints of $\Theta$ until they reach radius $R_1$.

There is some absolute constant $c>0$ such that 
\[
\PB{ \Qual(\Gamma) >\frac 1 4 \bigm | \mathcal{A}_\Theta (r,R_2)} = \nu_\Theta \Bigl[ \Qual(\Gamma) > \frac 1 4 \Bigr] > c \,. 
\]
In other words, conditioned on the event $\mathcal{A}_\Theta(r,R_2)$, interfaces which might be initially very close to each other tend to quickly ``separate''. 

If we further condition on $\{ \tau=1\}$ or $\{\tau=0\}$, the result is still valid (with possibly a smaller absolute constant $c$).

\end{lemma}

\proofof{Proposition \ref{pr.couplingfaces}}
We first prove the proposition in the case where $\tau=1$; the case $\tau=0$ will need an additional color switching argument.
Let $N=\floor{\log_8(\frac {R}{r})}$. For $0\leq i \leq N$, let $r_i = 8^{N-i}r$.


We will denote by $\omega_{\Theta}$ and $\omega_0$ the percolation configurations in $\mathcal{D}_{\Theta}$ and the square $Q$, sampled according to $\nu_{\Theta}$ and $\nu_0$, respectively.

Let $0\leq i < N$; assume that we have sampled under $\nu_{\Theta}$ the set of interfaces $\Gamma(r_i)$ which start from the four endpoints of $\Theta$ until they reach radius $r_i$ and that we sampled under $\nu_0$
the set of interfaces $\Gamma_0(r_i)$ defined similarly, and these two marginals ( $\Gamma(r_i)$ and  $\Gamma_0(r_i)$) are already coupled in some manner. These interfaces might be badly separated, which would prevent us from using RSW technology at the present scale $r_i$. 
Therefore, we rely on the Separation Lemma~\ref{l.separation}, which claims that if these interfaces are continued (under $\nu_\Theta$)
until radius $\frac {r_i} 2$, then with positive probability, they are well-separated around radius $\frac {r_i} 2$. More precisely,
if $\Gamma^i:=\Gamma(\frac {r_i} 2)$ 
denotes the continuation (under $\nu_\Theta$) of the interfaces $\Gamma(r_i)$ until radius $\frac {r_i} 2 = 4 r_{i+1}$ ($\Gamma^i_0$
is defined similarly), then with positive (conditional) probability $c>0$, one has
\[
\nu_{\Theta,\Gamma(r_i)} \bigl[ \Qual(\Gamma^i) > \frac 1 4 \bigr] > c\,.
\]
The same statement holds for $\Gamma_0(r_i), \Gamma_0^i$. From $\Gamma(r_i)$ and $\Gamma_0(r_i)$ already given, we are continuing the interfaces to get $\Gamma^i$ and $\Gamma^i_0$ independently from each other.


It is straightforward to check that the
initial configurations of faces, $\Theta$ and $\p Q$, plus the sets of interfaces $\Gamma^i$ and $\Gamma^i_0$ induce two
configurations of faces, $\Theta^i$ and $\Theta^i_0$, around radius $\frac {r_i} 2$. Now the law of $\nu_{\Theta}$ in $\mathcal{D}_{\Theta^i}$
conditioned on the explored interfaces $\Gamma^i$ is exactly the same as the law of the percolation in $\mathcal{D}_{\Theta^i}$ conditioned on
$\mathcal{A}_{\Theta^i}(r,r_{i})$ and $\{U_{\Theta^i}=\tau =1\}$ (a law which depends only on the faces in $\Theta^i$).
Let $\nu^i$ (and similarly $\nu^i_0$) denote this law. 

If the configurations of faces $\Theta^i$ and $\Theta^i_0$ happen to be well separated (i.e., $\Qual(\Theta^i)=\Qual(\Gamma^i)>1/4$, and the same for $\Theta^i_0$),
then we will see below that in some sense, the measures $\nu^i$ and $\nu^i_0$ are ``absolutely continuous'' with each other
within the annulus $A(r_{i+1}, 2r_{i+1})$. This property 
will allow us to couple these two measures at that scale (with positive probability).

Let $\mathcal{W}_i:=\{ \Qual(\Gamma^i) > 1/4\} \cap \{ \Qual(\Gamma^i_0) >1/4\}$ be the event that the interfaces are well-separated at radius $r_i/2$.
We have just seen above that conditioned on the interfaces up to radius $r_i$, $\Gamma(r_i)$ and $\Gamma_0(r_i)$, the event $\mathcal{W}_i$ is satisfied 
with conditional positive probability $c^2$ (recall that the conditioning here is under $\mathcal{A}_\Theta(r,R)$ $\tau=1$, $\Gamma(r_i)$ and $\Gamma_0(r_i)$). This means that if we discover the set of interfaces inductively along scales all the way from $R$ to $r$, there will be a positive proportion 
of scales for which $\mathcal{W}_i$ will hold.

Let us now assume that the event $\mathcal{W}_i$ holds at the present scale $\frac {r_i} 2$. 
Somehow we would like to control the Radon-Nikodym 
derivative ``$\frac {d\nu^i} { d\nu^i_0}$'' within the annulus $A(r_{i+1}, 2r_{i+1})$. 

We will need a few definitions. 
First, let $\Upsilon^i$
denote the set of all interfaces crossing from $2r_{i+1}$ to $r_{i+1}$ (therefore, by the definition of $\nu_{\Theta}$, there are at least 4 such interfaces under $\nu_\Theta$).
The set of interfaces $\Upsilon^i_0$ is defined in the same way. 
Define also $\mathcal{R}_i$ to be the event that for a percolation configuration in the annulus $A(r_{i+1}, 2 r_{i+1})$, the
set of all its interfaces which cross the annulus satisfies the events  $\mathcal{G}(r_{i+1}, 2r_{i+1})$,
$\Sep^{1/4}(r_{i+1}, 2r_{i+1})$ and $\Sep^{1/4}_+(r_{i+1}, 2r_{i+1})$ (in particular, if $\mathcal{R}_i$ holds, then
there are exactly 4 interfaces crossing from $2r_{i+1}$ to $r_{i+1}$).

Now we sample the sets of interfaces $\Upsilon^{i}$ and $\Upsilon^{i}_0$
according to the laws $\nu^i$ and $\nu^i_0$. Note that conditioning on $\Upsilon^{i}$ (or $\Upsilon^{i}_0$)
determines the color of the hexagons neighboring the interfaces in $\Upsilon^{i}$ ($\Upsilon^{i}_0$).
Let then $S^{i}$ ($S^i_0$) be the union of all hexagons whose color is determined by $\Upsilon^{i}$ ($\Upsilon^{i}_0$), together with the colors.
Let $S$ be a possible value for $S^i$ such that $\mathcal{R}_i$ holds (i.e., the four interfaces determined by $S$ satisfy the event $\mathcal{R}_i$). 
Let $m_S$ be the number of hexagons that are in $S$. Clearly, without conditioning, $\Pb{S^i=S} = 2^{-m_S-1}$.
We claim that if $\mathcal{W}_i$ holds (as we assumed above), then there is a universal constant $c>0$ such that
\begin{equation}\label{e.coupling}
2^{-m_S}/ c \leq \Pb{S^i = S \bigm| \mathcal{A}_{\Theta^i}(r, r_{i}),\, U_{\Theta^i}=\tau,\,\Theta^i} = \nu^i [S^i=S] \leq  c\, 2^{-m_S}.
\end{equation}

One also has, with the same constant $c$,
\begin{equation}\label{e.coupling1}
2^{-m_S}/ c \leq \nu^i_0 [S^i_0=S] \leq  c\, 2^{-m_S}.
\end{equation}

Indeed, on the event $\{\Qual(\Gamma^i)> 1/4\} \supset \mathcal{W}_i$,
\begin{align*}
\nu^i [S^i=S ] &= \frac{\Pb{\mathcal{A}_{\Theta^i}(r, r_{i}),\, U_{\Theta^i} = \tau \bigm| S^i=S,\,\Theta^i }
\Pb{S^i = S \md \Theta^i} }{\Pb{\mathcal{A}_{\Theta^i}(r,r_{i}),\, U_{\Theta^i} =\tau \mid \Theta^i}}\\
&= 2^{-m_S-1} \frac{\Pb{\mathcal{A}_{\Theta^i}(r,r_{i}),\, U_{\Theta^i}=\tau \bigm| S^i=S ,\,\Theta^i}} {\Pb{\mathcal{A}_{\Theta^i}(r,r_{i}),\,
U_{\Theta^i}=\tau \mid \Theta^i}}\,.
\end{align*}
Since $S$ satisfies the event $\mathcal{R}_i$, conditioning on it implies that the arms which cross the annulus $A(r_{i+1}, 2r_{i+1})$ are well-separated. As in Figure~\ref{f.gluing} of Appendix~\ref{s.appendix}, this implies that we can glue interfaces (using RSW and FKG) on both ends of $A(r_{i+1}, 2r_{i+1})$ with a cost of only a constant factor, hence there is a positive constant $C >0$ such that
\begin{eqnarray}
\alpha_4(r,  r_{i}) /C & \leq & \Pb{\mathcal{A}_{\Theta^i}(r,r_{i}),\, U_{\Theta^i}=\tau \bigm| S^i=S,\,\Theta^i} \leq C  \alpha_4(r,r_{i})
\nonumber\\
\alpha_4(r,  r_{i}) /C &\leq &\Pb{\mathcal{A}_{\Theta^i}(r,r_{i}),\, U_{\Theta^i}=\tau \mid \Theta^i} \leq C \alpha_4(r,r_{i})\,.\nonumber
\end{eqnarray}
This implies our claim (\ref{e.coupling}). Now by summing the claim (\ref{e.coupling}) over the different $S$ for which $\mathcal{R}_i$
holds, we get that if $\mathcal{W}_i$ holds then:
\begin{equation}
\Pb{\mathcal{R}_i \bigm| \mathcal{A}_{\Theta^i}(r, r_{i}),\, U_{\Theta^i} =\tau ,\, \Theta^i} \geq \frac 1 c \, \Pb{\mathcal{R}_i}\,.\nonumber
\end{equation}
We now use the following easy fact, which is part of the folklore, though we are not aware of an explicit proof in the literature:

\bl\label{l.only4arms} In the annulus $A=A(r,2r)$, there is a uniformly positive probability to have exactly 4 disjoint alternating arms crossing $A$, with the endpoints of the resulting 4 interfaces on the two boundaries of $A$ having quality at least $1/4$. That is, $\Pb{\mathcal{R}_i}> c$ for an absolute constant $c>0$.
\el

\proof
Let us provide two proofs, the first using Reimer's inequality \cite{\Reimer}, the second being longer, but completely elementary, hence having a chance to be generalizable to other models, like the FK random cluster model.

Let $\A$ be the event that there exist four alternating arms in $A$, with the additional property that the four interfaces given by the clockwise boundaries of the connected monochromatic components of $A$ containing these arms have well-separated endpoints: their interior and exterior qualities are both at least $1/4$. It is clear from RSW techniques that $\Pb{\A}>c_1>0$. Let $\mathcal{B}= \{$there is an arm (open or closed) from $\p_1A$ to $\p_2A\}$. Then,
by Reimer, $\Ps{\A \square \mathcal{B}} = \P[$polychromatic 5 arms in $A$, four of them giving well-separated interfaces]$\leq \Ps{\A}\, \Ps{\mathcal{B}}$. But RSW implies $\Ps{\mathcal{B}^c} > c_2 >0$, so we get $\Ps{\A \setminus \A \square \mathcal{B}}=\P[$exactly 4 arms in $A$,  giving well-separated interfaces$] > c_2 \Ps{\A} > c_1c_2 >0$, and we are done.

\begin{figure}[htbp]
\SetLabels
(.23*.8)$N$\\
(0.23*0.14)$S$\\
(.07*.48)$W$\\
(.18*.36)$C_W$\\
(.24*.36)$C_C$\\
(.29*.36)$C_E$\\
(.39*.48)$E$\\
\endSetLabels
\AffixLabels{
\centerline{
\includegraphics[height= 2.5 in]{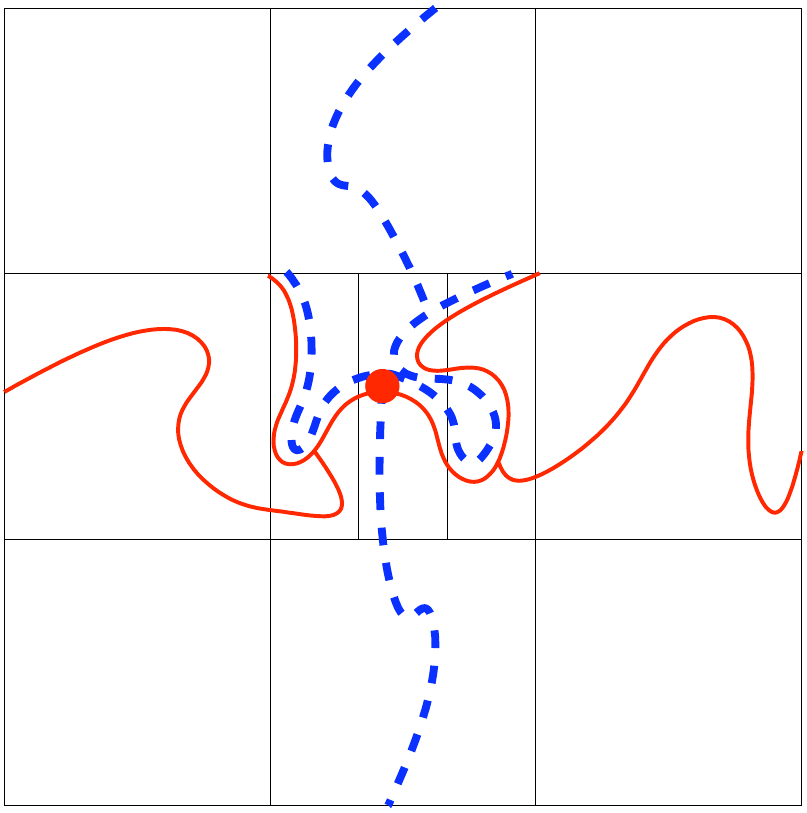}
\hskip 0.2 in
\includegraphics[height=2.5 in]{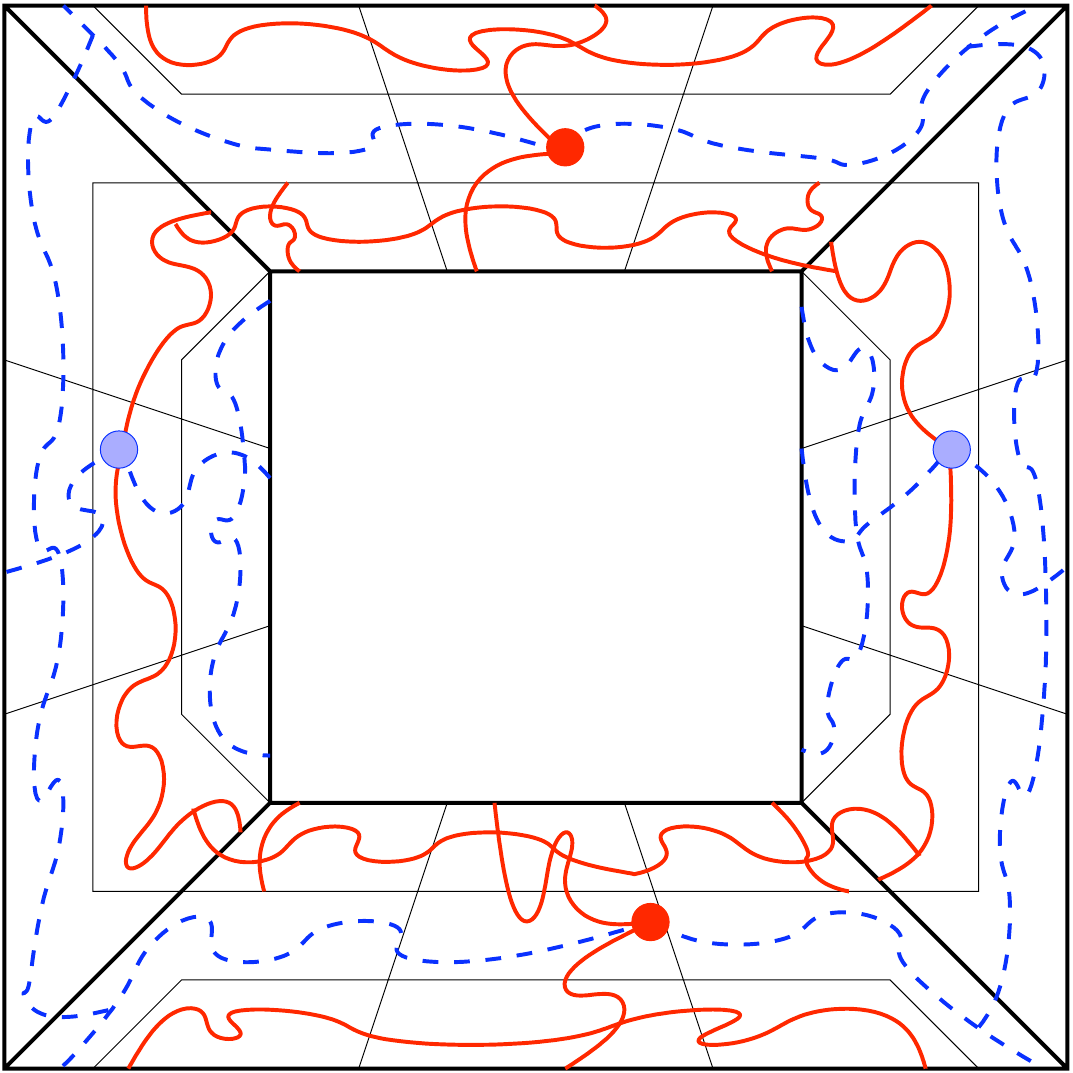}
}}
\caption{Constructing exactly four arms in an annulus, in a way that gives well-separated interfaces.}
\label{f.only4arms}
\end{figure}
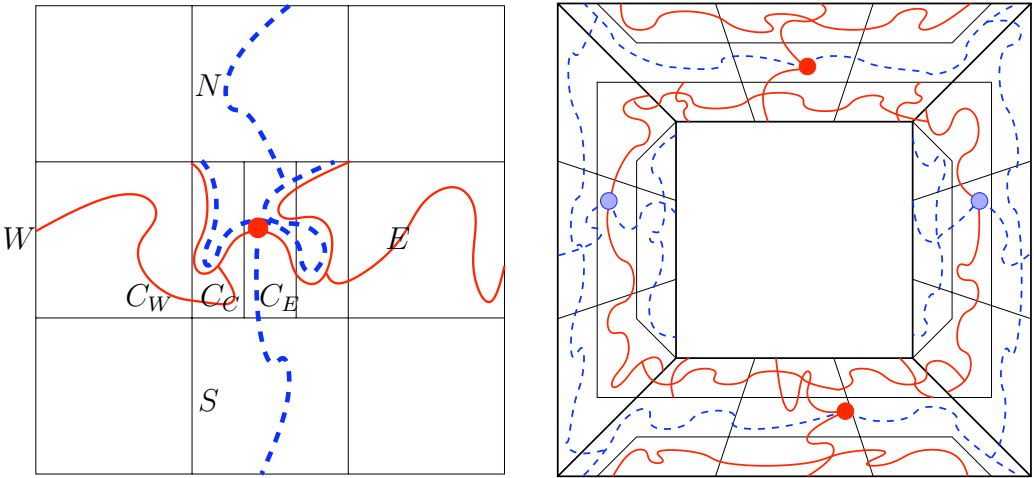

For the elementary proof, we first show that in the square there is a uniformly positive probability of the event that there is a left-right crossing but no two disjoint ones (which is the same as having at least one pivotal open bit), with the additional property that there are open arms from a pivotal bit to the left and right sides in the horizontal middle one-third of the square, and there are closed arms to the top and bottom sides in the vertical middle third. 
See the left hand picture of Figure~\ref{f.only4arms}. Let us denote the central ninth of the square by $C$, and its neighbouring ninths by $N,E,S,W$, in the obvious manner. Finally, let the left third of $C$ be $C_W$, the right third $C_E$, and the middle third $C_C$.

By RSW, there is a positive probability for a left-right open crossing in $C$; the uppermost such crossing can be found by running an exploration process from the upper left corner to the upper right corner of $C$, with open bits on the right (below the interface), closed bits on the left. So, the lower boundary of the set of explored bits is all open, the upper boundary is all closed, except the bits on the boundary of $C$, where they are unexplored (the interface only pretended that they satisfy the suitable boundary conditions). Furthermore, all bits beyond these boundaries are unexplored, hence unbiased. So, by RSW, with positive probability we get a closed path in $C_C\cup N$ connecting the upper side of the large square to the upper (closed) boundary of the interface, and another closed path in $C_C \cup S$ from the lower side of the large square to the lower (open) boundary of the interface. Where this closed arm hits the open side of the interface will be the pivotal we are looking for. Then, again by RSW, with positive probability there is an open path inside $W\cup C_W$ connecting the left side of the large square to the open side of the interface, and another open path in $C_E\cup E$ connecting the interface to the right side of the square. This finishes the construction for our first claim.

Now let us divide the annulus $A(r,2r)$ into four congruent sectors ($N$, $E$, $S$, $W$), and each sector into nine pieces, as shown on the right hand picture of Figure~\ref{f.only4arms}. All the sectors and the pieces are nice quads, so we can translate our previous claim to get the following: with positive probability, each of the odd sectors ($N,S$) contains a pivotal bit for an open crossing between the sides of the quad that are part of the annulus boundaries, each of the even sectors ($W,E$) contain a pivotal bit for a closed crossing between the sides that are part of the annulus boundaries, and all the open and closed arms emanating from these pivotals are contained in the middle pieces of their sectors. Then, using RSW and FKG, we can glue together these arms such that we get an open circuit in $A(r,2r)$ with just two closed bits on it: our closed pivotals in the sectors $W,E$; and can also get a closed circuit with just two open bits on it: our open pivotals in $N,S$. Finally, we can have additional connections in some of the pieces to ensure that the resulting four interfaces are well-separated; see again the picture. All of this happens with a uniformly positive probability, hence the proof is complete.
\QED

\proofcont{Proposition~\ref{pr.couplingfaces}} Summarizing what is above: we proved that at any scale $0\leq i < N =\floor{\log_8(\frac {R}{r})}$, if the event $\mathcal{W}_i$ holds
(which under $ \nu_\Theta \otimes \nu_0$ holds with positive probability), then there is some universal constant $c>0$ such that 

\begin{equation}\label{e.coupling2}
\Pb{\mathcal{R}_i \bigm| \mathcal{A}_{\Theta^i}(r, r_{i}),\, U_{\Theta^i}=\tau ,\, \Theta^i} \geq c\,,
\end{equation}
and
\begin{equation}\label{e.coupling3}
\Pb{\mathcal{R}_i \bigm| \mathcal{A}_{\Theta^i_0}(r, r_{i}),\, U_{\Theta^i_0}=1,\, \Theta^i_0} \geq c\,.
\end{equation}

\vskip 0.3 cm

We now have the tools which enable us to construct the coupling between $\nu_{\Theta}$ and $\nu_0$. We proceed by induction on the scale $r_i$, $i\geq 0$.
Assume we explored the sets of interfaces $\Gamma(r_i)$ and $\Gamma_0(r_i)$ (interfaces up to radius $r_i$) according to
some coupling of $\nu_{\Theta}$ and $\nu_0$. These interfaces might be badly separated at radius $r_i$. We sample their continuation (independently under $\nu_\Theta$ and $\nu_0$) until radius $r_i/2$ thus obtaining using the previous notations $\Gamma^i$ and $\Gamma^i_0$ (and the corresponding configurations of faces 
$\Theta^i$ and $\Theta^i_0$).
As we have seen, the event $\mathcal{W}_i$ that both configurations of faces are well-separated is satisfied with positive probability.

We now wish to sample the continuations $\Gamma(r_{i+1})$ and $\Gamma_0(r_{i+1})$ of these interfaces up to radius $r_{i+1}$ according to what have been explored so far, i.e., according to the laws $\nu^i$ and $\nu^i_0$.
If $\mathcal{W}_i$ does not hold, then the coupling ``fails'' at this scale and
we sample $\Gamma(r_{i+1})\setminus \Gamma^i$ and $\Gamma_0(r_{i+1})\setminus \Gamma^i_0$ independently, according to $\nu^i$ and $\nu^i_0$.
Else, if $\mathcal{W}_i$ holds, then by (\ref{e.coupling2}) and (\ref{e.coupling3}), with
positive probability for both $\nu^i$ and $\nu^i_0$, the sets of all interfaces crossing from $r_{i+1}$ to $2r_{i+1}$
satisfy $\mathcal{R}_i$. Furthermore, by (\ref{e.coupling}) and (\ref{e.coupling1}), any set of interfaces satisfying $\mathcal{R}_i$
has about the same probability under $\nu^i$ or $\nu^i_0$. This allows us to sample $\Gamma(r_{i+1})$
under $\nu^i$ and $\Gamma_0(r_{i+1})$ under $\nu^i_0$ and to couple these samples so that the conditional probability of
the event $\mathcal{S}_i:= \{ \Upsilon^{i} = \Upsilon^{i}_0 \in \mathcal{R}_i \}$ is greater than some absolute constant $c>0$.
Notice here that some care is needed since it is not $\Upsilon^{i}$ that we are sampling according to $\nu^i$, but rather $\Gamma(r_{i+1})\setminus \Gamma^i$;
indeed $\Upsilon^{i}$ could consist of more than 4 interfaces, but as we have seen, with positive probability there are only 4 of them,
and in that case we have $\Upsilon^{i} \subset \Gamma(r_{i+1})$. It is easy to check that, if $\mathcal{S}_i$ holds, then the interfaces $\Gamma(r_{i+1})$ and $\Gamma_0(r_{i+1})$ induce exactly the same configuration of faces $\Theta(r_{i+1}) = \Theta_0(r_{i+1})$ at radius $r_{i+1}$.

The induction stops when $\mathcal{S}_i$ has occurred  or when  $i$ reaches $N$, whichever happens first.
Let $i^*\leq N$ denote the index where the induction stopped. 
Call $\mathcal{S}$ the event that the induction
stopped at $i^*< N$, i.e., the event that the coupling succeeded. 
If  $i^*<N$, then the configurations of faces at radius $r_{i^*+1}$ are identical for both partially discovered configurations under
the coupled $\nu_{\Theta}$ and $\nu_0$. Let us call this configuration of faces $\Theta^*$. In both configurations, what remains to be
sampled in $\mathcal{D}_{\Theta^*}$ depends only on the faces $\Theta^*$, and since we assumed
$\tau=1$, they follow exactly the same law. Therefore, we can sample identically for $\nu_{\Theta}$ and $\nu_0$ and conditioned on $\Theta^*$
the 4 interfaces which start at the 4 endpoints of $\Theta^*$ until they reach radius $r$. These interfaces, plus the faces $\Theta^*$ if needed,  define the same faces around $r$ for the coupled $\nu_{\Theta}$ and $\nu_0$, as desired.

Let $\hat \nu$ denote the law of the coupling $(\nu_{\Theta}, \nu_0)$ we have just constructed.

It should be clear from what precedes that there is some universal constant $c>0$ such that, at each scale $0\leq i < N$, 
the coupling succeeds at scale $r_i$ with probability at least $c$. Indeed, knowing the interfaces $\Gamma(r_i), \Gamma_0(r_i)$,
up to radius $r_i$, one has positive probability to obtain well separated interfaces $\Gamma^i, \Gamma^i_0$ at radius $r_i/2$; then 
if these interfaces are indeed well separated, one couples the conditional measures $\nu^i$
and $\nu^i_0$ in such a way that the induced faces $\Theta(r_{i+1})$ and $\Theta_0(r_{i+1})$ are identical with positive probability. 
This implies that for any $0 \leq M \leq N$, the probability under the so-defined coupling $\hat \nu$
that $i^* \geq M$ is bounded above by $(1-C)^M$. Since $N=\lfloor\log_8(R/r)\rfloor$, there is indeed some absolute exponent $k>0$ such that 
\[
\hat \nu \bigl[  i^* < N\bigr] = \hat\nu \bigl[ \mathcal{S} \bigr] > 1 - (r/R)^k \,,
\]
which proves Proposition \ref{pr.couplingfaces} in the case $\tau=1$.
\vskip 0.3 cm

It remains to prove the case where $\tau=0$. The proof follows the exact same lines as in the case $\tau=1$, plus
the following color switching argument (we keep the same notations).
In the construction of the coupling $\hat\nu$, suppose we are at scale $i<N$ and that we already sampled $\Gamma^i$ and $\Gamma^i_0$.
If $\mathcal{W}_i$ holds, this allows us to sample $\Gamma(r_{i+1}) \setminus \Gamma^i$ under $\nu^i$
and $\Gamma_0(r_{i+1}) \setminus \Gamma^{i+1}_0$ under $\nu^i_0$ and to couple these samples so that the conditional probability
of the event $\tilde{\mathcal{S}}_i:=\{ \neg \Upsilon^{i} = \Upsilon^{i}_0 \in \mathcal{R}_i \}$ is greater than some absolute constant
$c>0$, where $\neg \Upsilon$ is the color-switched of $\Upsilon$. With the same proof, the coupling succeeds with probability at least $1-(r/R)^{k}$,
and in that case one ends up with two identical faces $\Theta(r)$ and $\Theta_0(r)$ with reversed colors.
\QED


\begin{remark}\label{r.Z2}
The proof adapts easily to the case of bond percolation on $\Z^2$. The only detail to handle is the color switching argument, since the duality on $\Z^2$ is not between open and closed hexagons in the same lattice, but open edges and closed dual edges (and vice versa). Since the dual of a percolation configuration in an $r\times (r+1)$ box lives in the $(r+1)\times r$ dual box, it is better to take annuli $A(r,R)=B(R)\setminus B(r)$, where $B(r)$ is an  $r\times (r+1)$ box. In Proposition~\ref{pr.couplingfaces}, we have the measure $\nu_\Theta$ with $\tau=0$, and want to couple it with $\nu_0$ that has $\tau=1$. When we do the color switching on $\nu_\Theta$, we get an annulus between $(R+1)\times R$ and $(r+1)\times r$ with the types of the faces switched and the conditioning changed to $\tau=1$. But now this is just a 90 degrees rotated version of a standard $\nu_\Theta$ measure with the $\tau=1$ condition, hence we can do our standard coupling, and the issue is solved.
\end{remark}

We will also need the following proposition, whose statement and proof are very similar to those of Proposition~\ref{pr.couplingfaces}.

\begin{proposition}\label{pr.coupling2}
Let $\Omega$ be some piecewise smooth simply connected domain with $0\in \Omega$. Let $d>0$ be the distance
from 0 to the boundary $\p \Omega$ and $d'= d\wedge 1$. For any $0<r<d'$, let $\mathcal{A}(r,\p \Omega)$ be the event that there
are four alternating arms from the square of radius $r$ to the boundary $\p \Omega$. Also, let $\mathcal{A}(r,1)$ be
the event that there are four arms from radius $r$ to radius 1.
For any $10 \eta < r < \frac {d'} {10}$, let $\nu$ be the law on percolation configurations
conditioned on the event $\mathcal{A}(r,\p \Omega)$, and let $\nu_0$ be the law
conditioned on the event $\mathcal{A}(r,1)$.

Then, there is a coupling of the conditional laws $\nu$ and $\nu_0$ such that with (conditional) probability
at least $1 - (r/d')^k$, the event $\mathcal{G}(r,d'/2)$ (as defined previously) is satisfied for both
configurations and the induced faces at radius $r$ are identical. Here $k>0$ is some absolute exponent.
\end{proposition}

\noindent{\bf Sketch of proof.} There is a small issue here compared to Proposition~\ref{pr.couplingfaces} which is that there is no initial configuration of faces to start with. In particular, if one explores the percolation configuration from $\p \Omega$ towards some radius $\gamma$, there might be more 
than four faces around $\gamma$. A way to deal with this minor issue is to choose an intermediate scale $\gamma$ (say, $:= \sqrt{r}$) and then argue (using quasi-multiplicativity) that with high conditional probability (under $\mathcal{A}(r,\p\Omega)$), there are only four arms from $\gamma$ to $\p \Omega$. In more detail, using $\mathcal{A}_4(\cdot,\cdot)$ for the event of having at least four arms of alternating colors, and $\mathcal{A}_5(\cdot,\cdot)$ for the event that there are at least five arms, with four of them of alternating colors, we have
\begin{align*}
\Pb{\A_4(\gamma,\p\Omega) \setminus \A_5(\gamma,\p\Omega) \md \A_4(r,\p\Omega)} &\geq 1 - \frac{\Pb{ \A_5(\gamma,\p\Omega),\ \A_4(r,\p\Omega)}}{\Pb{\A_4(r,\p\Omega)}}\\ 
&\geq  1 - \frac{\Pb{\A_4(r,\gamma)}\,\Pb{\A_5(\gamma,\p\Omega)}}{\Pb{\A_4(r,\p\Omega)}}\\ 
& \geq 1 - C \frac{\Pb{\A_5(\gamma,\p\Omega)}}{\Pb{\A_4(\gamma,\p\Omega)}}=1-o(1)\,,
\end{align*}
where the last step (comparing the 4- and 5-arm probabilities) is a typical application of the BK (van den Berg-Kesten) inequality, see \cite{\Grimm}.

On this large probability event, we have only four faces at radius $\gamma$, and can run the proof of Proposition~\ref{pr.couplingfaces} from $\gamma$ to $r$. \qed
\medskip

The same proofs work for the 2-arm case, just using two faces instead of four. The 1-arm case is different, since one arm does not produce faces that separate the inside from the outside. Instead, one may use open circuits there, and the proof is actually simpler, since, due to monotonicity, one can use just FKG, no separation of interfaces is needed. See Section~\ref{s.others} for more detail on these other cases.

\section[The {\it pivotal measure}]{A measure on the set of pivotal points which is measurable with respect to the scaling limit}\label{s.measurable}

\begin{minipage}{\textwidth}
\begin{center}
\includegraphics[width=0.5 \textwidth]{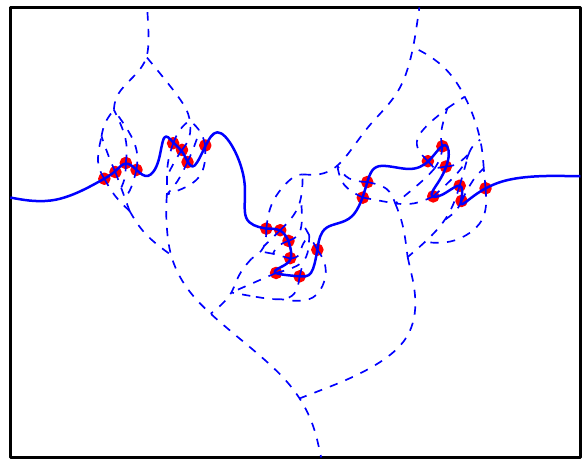}
\end{center}
\end{minipage}
\vskip 0.5 cm


\subsection{General setup and statement of convergence}\label{ss.A-imp}

In this section we define a natural scaling limit of the counting measure on pivotals, normalized so that the set of points in a nice bounded domain (say, $[0,1]^2$) that are 1-important (i.e., have the alternating 4-arm event to distance at least 1, as in Subsection~\ref{ss.results}), or which are pivotal for a nice quad comparable to that domain (say, the same $[0,1]^2$), have expected measure of order 1.

Of course, on the discrete level, for a given configuration $\omega_\eta$, any point is pivotal (or important) on some scale; for example, if the point happens to be surrounded by at least 5 hexagons of the same color, this is scale is only $\eta$. So, in order to keep a meaningful measure at the limit, we will only keep track of points in the mesh which are pivotals on a macroscopic scale.
 
This is consistent with our ultimate goal of proving that dynamical percolation, as well as near-critical percolation, correctly rescaled, has a scaling limit. Indeed, as we outlined in Subsection~\ref{ss.DPSL}, if one is interested in the trajectory $t\mapsto \omega^\eta_t$ up to a certain small precision, it is enough to know the initial configuration $\omega^\eta_0$, as well as what happened on the set of points which were initially at least $\rho$-important (the choice of $\rho$ depends on the amount of precision that one requires).
   
This step is made rigorous in \cite{\DPSL}. In the current paper, we focus
 on the scaling limit of the following (random) counting measures in the plane,
 parametrized by the scale ``cut-off'' $\rho>0$:
 \begin{equation}
 \bar\mu_{\eta} ^{\rho} = \bar\mu_{\eta}^\rho(\omega_\eta):= \sum_{x:\,\rho\textrm{-important}} \delta_x\, \frac {\eta^2} {\alpha_4^\eta(\eta,1)}\,.\label{e.rhomeasure}
 \end{equation}

Roughly speaking, $\rho$-important points are the same as points that are pivotal for some ``$\rho$-macroscopic'' quad with boundary ``not too close'' to the point. For convenience, we will actually consider measures which are defined in a somewhat less symmetric way, but will turn out to be less sensitive to local effects in our proof. This needs some definitions:

\begin{definition}
Let $A$ be some closed topological annulus of the plane; the bounded component of $A^c$ will be called its {\bf inside face}; $\p_1 A$ and $\p_2 A$ will be its inner and outer boundaries. We will say that $A$ is a {\bf proper annulus} if its inner and outer boundaries are piecewise smooth curves of finite length. 
\end{definition}

\begin{definition}
Let $\mathcal{H}$ be a family of proper annuli in the plane. It will be called
an {\bf enhanced tiling} if the collection of the inner boundaries form a locally finite tiling of the plane.
The {\bf diameter} of $\mathcal{H}$ will be
\begin{equation}
\diam\, \mathcal{H} := \sup_{A \in \mathcal{H}}\{ \diam\,{A} \}.
\nonumber
\end{equation}
\end{definition}

Let $A$ be some proper annulus of the plane with inside face $\Delta$; a point $x\in\Delta$ in the triangular grid of mesh $\eta$
is called {\bf $A$-important} if there are four alternating arms in $\omega_{\eta}$ from $x$ to the outer boundary $\p_2 A$ of $A$.
We define the measure
\begin{equation}
\mu_{\eta}^A := \mu_{\eta}^A (\omega_{\eta})= \sum_{x\in \Delta:\,A\textrm{-important}} \delta_x\, \frac {\eta^2}{ \alpha_4^\eta(\eta,1)}\,.
\label{e.Ameasure}
\end{equation}

Now, let $\mathcal{H}$ be some enhanced tiling, and define
\begin{equation}
\mu_{\eta}^{\mathcal{H}} := \mu_{\eta}^\mathcal{H} (\omega_{\eta})= \sum_{A \in \mathcal{H}}  \mu_{\eta}^A\,.
\label{e.Hmeasure}
\end{equation}
Note that in the measure $\mu_{\eta}^{\mathcal{H}}$ thus defined, we do not count the points which might lie on the
boundary $\p_1 A$ of some annulus $A\in \mathcal{H}$, but this has a negligible effect: indeed, for a fixed proper annulus $A$,
it is straightforward to check that when $\eta$ goes to 0, the probability that there is some $A$-important hexagon $x$
intersecting $\p_1 A$ is going to zero. 

We will now prove that for any fixed proper annulus $A$, the random measure $\mu_{\eta}^A(\omega_{\eta})$ has a scaling
limit $\mu^A$ when $\eta$ goes to zero, and moreover $\mu^A = \mu^A(\omega)$ is a measurable function of the continuum
percolation $\omega$, as defined in Section~\ref{ss.topology}. More precisely, we have the following theorem:

\begin{theorem}\label{th.measurable}
Let $A$ be a fixed proper annulus of the plane. When $\eta\to 0$, the random variable $(\omega_{\eta}, \mu_{\eta}^A)$ converges in law
to some $(\omega, \mu^A)$, where $\omega$ is the scaling limit of critical percolation, and the measure $\mu^A= \mu^A(\omega)$ is a measurable function of $\omega$.
\end{theorem}

A simple regularity property of the limit measure is the following:

\begin{proposition}\label{pr.Mom}
Let $A$ be a proper annulus with $\dist(\p_1 A,\p_2 A) \geq \eps$ and $\diam A<5\eps$ (the constant 5 is of course rather arbitrary here). Let $U$ be a bounded open subset of the inner face of $A$, with a boundary $\p U$ that has finite length. Then
$$
\Es{\mu^A(U)} \asymp \frac{\area(U)}{\alpha_4(\eps,1)} 
\asymp \eps^{-5/4}\,\area(U)\,.
$$
Assume now that $U$ is a ball or square of radius $r$, contained in the inner face of $A$. Then
$$
\Es{\mu^A(U)^2} 
\asymp \frac{1}{\alpha_4(\eps,1)} \frac{r^4}{\alpha_4(r,1)} 
\asymp  \eps^{-5/4}\,r^{11/4}.
$$
The constant factors implied by the notation $\asymp$ are independent of everything. 
\end{proposition}

The connection between the $\rho$-important measures $\bar \mu^\rho_\eta$ and the annulus-important and enhanced tiling measures $\mu^A_\eta, \mu^\mathcal{H}_\eta$ will be made in Subsection~\ref{ss.filtering}.

\subsection{Tightness and moment bounds}\label{ss.tightness}

Firstly, we show that the family of variables $\{(\omega_{\eta}, \mu_{\eta}^A)\}_{\eta>0}$ is tight. In Section~\ref{ss.topology}, we defined $\omega_\eta$ as a Borel measure on the compact separable metrizable space $(\HH_D,\T_D)$, hence $\{\omega_{\eta}\}_{\eta>0}$ is obviously tight.
It is a standard fact that if $\{X_{\eta}\}_{\eta>0},\, \{Y_{\eta}\}_{\eta>0}$ are tight families of variables, then, in any coupling, the coupled family of variables $\{(X_{\eta}, Y_{\eta})\}_{\eta>0}$ is tight, as well. So, it is enough to prove that the family of measures $\{\mu_{\eta}^A\}_{\eta>0}$ is tight. Since $\mu_{\eta}^A$
are finite measures supported on the inside face $\Delta$, proving tightness boils down to proving
$$
\limsup_{\eta\rightarrow 0} \Eb{\mu_{\eta}^A (\Delta)} < \infty\,.
$$
This is straightforward by the definition of $\mu_{\eta}^A$. Indeed, let $d>0$ be the distance between $\p_1 A$ and $\p_2 A$;
a point in $\Delta$ has to be $d$-important in order to be $A$-important. Therefore, if $d'= d\wedge 1/2$:
\begin{equation}\label{e.tight}
\begin{aligned}
\Eb{\mu_{\eta}^A(\Delta)} & =  \sum_{x\in \Delta} \Pb{x\,\textrm{is}\,A\textrm{-important}}\, \eta^2\, \alpha_4^\eta(\eta,1)^{-1}\\
& \leq  \sum_{x\in \Delta} \alpha_4^\eta(\eta, d')\, \eta^2\, \alpha_4^\eta(\eta,1)^{-1}\\
&\asymp  \area(\Delta)\, \alpha_4(d',1)^{-1} < \infty.
\end{aligned}
\end{equation} 
The approximate equality in the last line follows from quasi-multiplicativity, from the fact that $\alpha_4^{\eta}(d',1)$ depends on $\eta$ but converges to the macroscopic probability $\alpha_4(d',1)$ when $\eta$ goes to zero, and from the observation that 
\begin{equation*}
\lim_{\eta\to 0} \big|\Delta\cap\eta\Tg\big| \, \eta^2  = \area(\Delta)\,.
\end{equation*} 
This observation can be proved by noticing that the area of points in $\Delta$ with distance at least $\eta$ from $\p_1 A$ converges to $\area(\Delta)$ as $\eta\to 0$, giving asymptotically $\eta^{-2}\area(\Delta)$ lattice points in $\Delta\cap\eta\Tg$, while the number of lattice points at distance at most $\eta$ from $\p_1A$ is at most $O(\eta^{-1}\length(\p_1 A))$. 

This proves tightness of $\{(\omega_{\eta}, \mu_{\eta}^A)\}_{\eta>0}$. Therefore, there exists some subsequential scaling limit $(\omega, \mu^A)$ along some subsequence $\{\eta_k\}_{k>0}$, where $\eta_k$ goes to 0. We will show that this $\mu^A$ can actually be recovered from $\omega$, which is the unique subsequential scaling limit of $\{\omega_{\eta}\}_{\eta>0}$, as we already know from Section~\ref{ss.topology}. Consequently, the pair $(\omega, \mu^A)$ will also be unique.
\medskip

Beyond the first moment, one can also easily bound higher moments. The following lemma will be useful later.

\begin{lemma}\label{l.Mom}
Let $A$ be a proper annulus with $\dist(\p_1 A,\p_2 A) \geq \eps$ and $\diam A<5\eps$. 
Assume that $U$ is a bounded open subset of the inner face of $A$, with a boundary $\p U$ that has finite length. Then
\begin{equation}\label{e.1stMom}
\Es{\mu_\eta^A(U)} \asymp \frac{\area(U)}{\alpha_4(\eps,1)} 
\asymp  \eps^{-5/4}\,\area(U)\,,
\end{equation}
with constant factors implied by the notation $\asymp$ that are independent of everything.
Assume now that $U$ is a ball or square of radius $r$, contained in the inner face of $A$. Then
\begin{equation}\label{e.2ndMom}
\Es{\mu_\eta^A(U)^2} 
\asymp \frac{1}{\alpha_4(\eps,1)} \frac{r^4}{\alpha_4(r,1)} 
\asymp  \eps^{-5/4}\,r^{11/4},
\end{equation}
for all small enough $\eta<\eta_0(r)$. Similarly,
\begin{equation}\label{e.3rdMom}
\Es{\mu_\eta^A(U)^3} 
\asymp \frac{1}{\alpha_4(\eps,1)} \frac{r^6}{\alpha_4(r,1)^2} 
\asymp  \eps^{-5/4}\,r^{7/2}.
\end{equation}
\end{lemma}

\proof The first moment bound~(\ref{e.1stMom}) is proved exactly the same way as~(\ref{e.tight}), just we have now equality up to constant factors  instead of just a one-sided bound due to the geometric properties of $A$ and quasi-multiplicativity.

For the second moment, note that $\Es{\mu_\eta^A(U)^2} = \sum_{x,y\in U}\Ps{x,y\textrm{ are $A$-important}}$. Now, for any $k\geq 0$ such that $2^k \eta < r $, there are $\Theta(r/\eta)^2 \, 2^{2k}$ pairs of points $(x,y)$ in $U\cap \eta\Tg$  satisfying $2^k \eta \leq d(x,y)<2^{k+1} \eta$. Furthermore, by quasi-multiplicativity of 4-arm probabilities, the probability that both $x$ and $y$ are $A$-important is comparable to the probability of the alternating 4-arm events from both points to distance about $2^k\eta$ and from a ball of radius about $2^k\eta$ containing both points to $\p_2 A$, which is at distance about $\eps$. Therefore, using quasi-multiplicativity several times,
\begin{align*}
\Es{\mu_\eta^A(U)^2} &\asymp \sum_{k=0}^{\log_2(r/\eta)} \frac{\alpha_4^\eta(\eta,2^k\eta)^2 \, \alpha_4^\eta(2^k\eta,\eps)}{(\alpha_4^\eta(\eta,1)\eta^{-2})^2} (r/\eta)^2\, 2^{2k}\\
&\asymp \frac{r^2}{\alpha_4^\eta(\eps,1)} \frac{\eta^2}{\alpha_4^\eta(\eta,1)} \sum_{k=0}^{\log_2(r/\eta)} \alpha_4^\eta(\eta,2^k\eta)2^{2k}\\
&\asymp \frac{r^2}{\alpha_4^\eta(\eps,1)} \frac{\eta^2}{\alpha_4^\eta(\eta,1)} (r/\eta)^2 \, \alpha_4^\eta(\eta,r)\\
&\asymp \frac{1}{\alpha_4(\eps,1)} \frac{r^4}{\alpha_4(r,1)}\,,
\end{align*}
where we used $R^2 \alpha_4^\eta(\eta,R\eta) = \Omega (R^{\zeta})$ with some $\zeta>0$ as $R\to\infty$  (uniformly in $\eta$) to get the third line. This proves~(\ref{e.2ndMom}).

The calculation for the third moment is very similar. We need to sum the probabilities $\Ps{x,y,z\textrm{ are all $A$-important}}$, which we can do by ordering the three points such that $2^k\eta\leq d(x,y) < 2^{k+1}\eta$ and $2^\ell\eta\leq d(y,z) <2^{\ell+1}\eta$, with $\ell\leq k$, and using quasi-multplicativity: 
\begin{align*}
\Es{\mu_\eta^A(U)^3} &\asymp \sum_{k=0}^{\log_2(r/\eta)} \sum_{\ell=0}^k \frac{\alpha_4^\eta(\eta,2^k\eta)\, \alpha_4^\eta(\eta,2^\ell\eta)^2 \, \alpha_4^\eta(2^\ell\eta,2^k\eta)\, \alpha_4^\eta(2^k\eta,\eps)}{(\alpha_4^\eta(\eta,1)\eta^{-2})^3} (r/\eta)^2\, 2^{2k}\, 2^{2\ell}\\
&\asymp \frac{r^2}{\alpha_4^\eta(\eps,1)} \frac{\eta^4}{\alpha_4^\eta(\eta,1)^2} \sum_{k=0}^{\log_2(r/\eta)}  \alpha_4^\eta(\eta,2^k\eta) \sum_{\ell=0}^k \alpha_4^\eta(\eta,2^\ell\eta)\,2^{2\ell}\\
&\asymp \frac{r^2}{\alpha_4^\eta(\eps,1)} \frac{\eta^4}{\alpha_4^\eta(\eta,1)^2}  \sum_{k=0}^{\log_2(r/\eta)} 
\alpha_4^\eta(\eta,2^k\eta)^2\,2^{2k}\\
& \asymp \frac{r^2}{\alpha_4^\eta(\eps,1)} \frac{\eta^4}{\alpha_4^\eta(\eta,1)^2} (r/\eta)^4 \alpha_4^\eta(\eta,r)^2\\
&\asymp \frac{1}{\alpha_4(\eps,1)} \frac{r^6}{\alpha_4(r,1)^2}\,,
\end{align*}
which proves~(\ref{e.3rdMom}).
\qed

\subsection{Strategy for the proof of uniqueness}

Since the subsequential scaling limit $\mu^A$ is a measure, we will not need to check sigma-additivity etc., we will only need to characterize uniquely the law of this random measure.
For this, it will be enough to determine, for any ball $B \subset \Delta$, the value of $\mu^A(B)$ as a function of $\omega$.
The strategy of the proof is as follows. If $\eta$ is a small mesh,
we are interested in $\mu_{\eta}^A (B)$, that is, in the number $X=X_\eta$ of points inside $B$ which are $A$-important (indeed,
by definition, $\mu_{\eta}^A (B) : = \frac{X} {\eta^{-2} \alpha_4^\eta(\eta,1)}$). We first need to show that $X$ can be ``guessed''
with arbitrarily good precision when one relies only on ``macroscopic'' information, since only the macroscopic information is preserved in the scaling limit $\omega$. Then we will also have to prove the actual convergence of the measures $\mu_\eta^A$. For the approximation from  macroscopic data, fix any grid of squares of radius $\eps$ (i.e., of side-length $2\eps$), and let $Y=Y^\eps_\eta$ be the number of squares $Q$ in this grid
which are contained in $B$ and which satisfy a 4-arm event from $2Q$ to $\p_2 A$. Why are we choosing $2Q$ instead of simply $Q$? This will save us from some boundary issues in our later conditionings. 

One would like to show that, knowing $Y=Y^\eps_\eta$, one can guess with good precision what $X=X_\eta$ is.
In other words, for some well-chosen factor of proportionality $\beta=\beta(\eps,\eta)$, one should have 
$X \approx \beta \, Y$ in some sense.
We will see that there is a natural candidate for this $\beta$-factor, and will show that 
for this $\beta$, one has $\Eb{|X -\beta Y|^2} = o(\Eb{X}^2)=o(\Eb{\beta Y}^2)$.
It will be convenient to think of this property as a kind of law of large numbers, even though we will in fact only achieve the above $L^2$-control (which fortunately will be sufficient for our ultimate goal).

Before going any further, let us make the above setup more precise.
First of all, let us fix the grid of $\eps$-squares to be the grid $G=G^\eps(a,\theta)$ of the $\eps$-squares centered at the points of the lattice 
$2\eps \,e^{i \theta} \Z^2 + a$, with $\theta\in\R$ and $a\in\C$. 
Why should one consider all possible orientations of the $\eps$-grid? In Section~\ref{s.covariance}, we will show that the asymptotic measure we are constructing 
satisfies some conformal covariance properties (in particular, rotational invariance). With this perspective in mind,
we should not restrict ourselves to a macroscopic lattice of $\eps$-squares conveniently 
chosen with respect to the orientation of our triangular grid $\eta \Tg$. This is the reason for this additional technicality.

The following proposition will be the main step in the proof of the uniqueness of the limiting measure.
\begin{proposition}\label{l.XY}
Let $A$ be some proper annulus, and let $B\subset \Delta$ be a ball included in the inside face of $A$.
Consider an $\eps$-grid $G$ corresponding to the lattice $2\eps\,e^{i \theta} \Z^2 + a$, with $\theta\in\R$ and $a\in\C$.

Recall that $X=X(\omega_\eta)$ denotes the number of sites in $\eta \Tg \cap B$ which are $A$-important and $Y=Y^\eps(\omega_\eta)$ is the number of $\eps$-squares $Q$ in $G\cap B$ which satisfy the four-arm event from $2Q$ to $\p_2 A$.

Then, there is a deterministic quantity $\beta=\beta(\eta,\eps,a,\theta)$ which satisfies
\begin{equation}\label{e.XY}
\Eb{(X-\beta Y)^2} \le \xi(\eta,\eps)\, \Eb{X}^2\,,
\end{equation}
where $\xi(\eta,\eps)$ goes to zero when $\eps$ and $\eta/\eps$ go to zero. Notice that the control we achieve is independent of 
the orientation $a,\theta$ of the grid (which will be relevant later) but depends on $B\subset \Delta$ and $A$. 
\end{proposition}

Let us sketch the proof of this second moment estimate. 
Let $Q_1, \ldots, Q_p$ denote the list of $\eps$-squares in the grid $G$ which are contained in $B$. For each $1\le i \le p$, define $y_i$
to be the indicator function of the event that there are 4 arms from $2Q_i$ to $\p_2 A$; define also $x_i$ to be the number of $A$-important points 
inside $Q_i$. Then one has $X:= x_1+\ldots+x_p$ and $Y=y_1+\ldots+y_p$. Formally, we neglected here some boundary issues,
since one should also consider the points lying in $B\setminus \bigcup Q_i$. We will come back to this minor issue in the actual proof.
We wish to obtain the estimate $\Eb{|X - \beta Y|^2} = o(\Eb{X}^2)$.
In particular, by Cauchy-Schwarz, $\beta$ needs to satisfy $|\Eb{X-\beta Y}|= o(\Eb{X})$. Let us then analyze the first moment
$\Eb{X-\beta Y}$. It can be rewritten as
\begin{eqnarray}
\Eb{X- \beta Y} &= & \sum_i \Eb{x_i - \beta y_i} \nonumber \\
& = & \sum_i \Pb{y_i=1} \bigl( \Eb{x_i \bigm| y_i=1} - \beta \bigr)\,, \nonumber
\end{eqnarray}
since $y_i=0$ implies $x_i=0$. 
Now, for any $i$, $\Eb{x_i\bigm| y_i = 1}$ is the expected number of $A$-important points inside $Q_i$ knowing that there are already four arms 
from $2Q_i$ to $\p_2 A$. The coupling result Proposition~\ref{pr.coupling2} essentially says that this quantity should depend very little on the shape of
$\p_2 A$ and on the position of $Q_i$ with respect to $\p_2 A$. In other words, by the coupling argument, we get that all these quantities 
$\Eb{x_i \bigm| y_i=1}$ are almost the same, and this common ``agreement'' is precisely what $\beta$ should represent. One could choose for $\beta$ precisely one of these quantities (say,
$\Eb{x_1 \bigm| y_1=1}$), but it will be convenient to define $\beta$ in a more universal way: i.e., based on some ``reference'' square $Q_0$ and some universal quad independent of $\p_2 A$
(in the spirit of Proposition~\ref{pr.couplingfaces}). The details are postponed to the actual proof.

   
   We sketched this first moment analysis only to highlight where the $\beta$ factor might come from. What we really need is the second moment.
It can be written as
\begin{align}\label{e.secmom}
\Eb{(X- \beta Y)^2} &=  \sum_{i,j} \Eb{\bigl( x_i - \beta y_i\bigr) \bigl(x_j - \beta y_j \bigr) } \nonumber \\
& =  \sum_{i,j} \Eb{y_i y_j}  \Eb{ \bigl( x_i - \beta \bigr) \bigl(x_j - \beta \bigr) \bigm| y_i=1,\, y_j=1}\,.
\end{align}

In this sum over squares $Q_i,Q_j$, there are relatively few pairs of nearby squares (say, whose distance $d(Q_i,Q_j)$ is less than $r$, with $\eps \ll r \ll 1$).
Hence one can neglect these pairs of squares, or in other words, we can neglect the near-diagonal term in~(\ref{e.secmom}) for a well-chosen 
scale $r$. Now, for any two squares $Q_i,Q_j$ such that $d(Q_i,Q_j)\geq r \gg \eps$,
if one conditions on $\{y_i=1, y_j=1\}$, i.e., on the event that both $2Q_i$ and $2Q_j$ are $A$-important, then, again by the coupling argument, one expects that the configuration seen
inside $Q_i$ should be almost independent of what is seen inside $Q_j$. This should lead to $\Eb{ \bigl( x_i - \beta \bigr) \bigl(x_j - \beta \bigr) \bigm| y_i=1,\, y_j=1} \approx
\Eb{x_i -\beta \bigm| y_i=1} \Eb{x_j -\beta \bigm| y_j=1}$. Now, since $\beta$ is precisely chosen to match well with these first moment quantities, one should indeed obtain a small
second moment in~(\ref{e.secmom}).

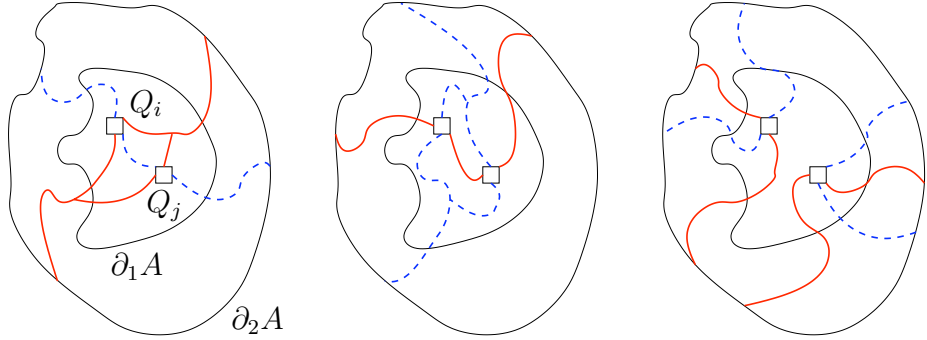
\begin{figure}[htbp]
\centerline{\SetLabels
(.15*.18)$\p_1 A$\\
(.28*.02)$\p_2 A$\\
(.16*.67)$Q_i$\\
(.18*.37)$Q_j$\\
\endSetLabels
\AffixLabels{
\includegraphics[height=1.75 in]{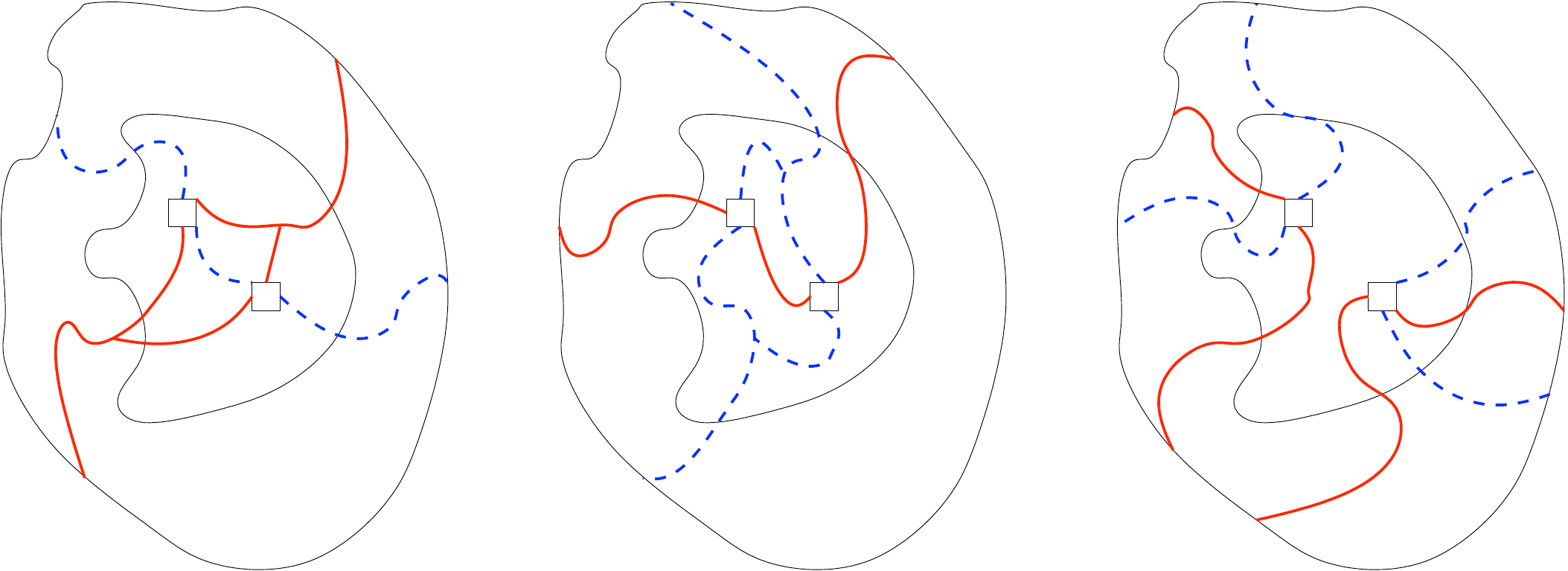}}}
\caption{Assume that all disjoint arms are shown between $Q_i$, $Q_j$ and $\p_2A$. Then, in the first example, $x_j>0$ requires a blue (dashed) connection within $Q_i$, while, in the second one, it requires a red (solid) connection. In the third example, $x_i$ and $x_j$ are independent.}
\label{f.insideeffect}
\end{figure}

However, in order to apply the coupling results from Section~\ref{s.coupling},
there are some issues about the above conditioning. 
Conditioned on the entire configuration outside $Q_i$ and $Q_j$ in a way that makes $\{y_i=y_j=1\}$ possible, the value of $x_j$ and $y_j$ might not at all be independent of the configuration inside $Q_i$. See Figure~\ref{f.insideeffect}.
This shows that somehow the configurations inside $Q_i$ and $Q_j$ interfere with each other in a nontrivial way, which is bad news for applying our coupling results from Section~\ref{s.coupling}. 

To handle this interference, first take some enlarged square $Q_i'$ around $Q_i$ so that $Q_j\subset (Q_i')^c$.
Under the conditioning $\{y_i=y_j=1\}$, one can further condition on all the percolation configuration seen in $(Q_i')^c$ (and thus, in particular, we know all the information inside $Q_j$). The hope is that if $Q_i'$ is much larger than $Q_i$, then the coupling results will imply that the configuration inside $Q_i$ is basically independent of everything we have conditioned on, except that we still have the same interference as before: even though we know all the information inside $Q_j$, the set of $A$-important points in $Q_j$ might still depend on how things are connected 
inside $Q_i'$. If the configuration around $Q_i'$ has many faces (with the definition of Section~\ref{s.coupling}), then this connectivity information could be quite complicated, hence it would be very convenient to assume that there is a configuration $\Theta$ of only four faces around $Q_i'$ (four being the minimum number compatible with the condition $\{y_i=1\}$). Then the only information inside $Q_i'$ which remains relevant outside $Q_i'$ is whether the open arcs of $\Theta$ are connected or not inside $Q_i'$ (or rather $\mathcal{D}_\Theta$).
This missing internal information is given by the random variable $U_\Theta$ introduced in Section~\ref{s.coupling}. If $\mathcal{F}$ denotes the filtration generated 
by the sites in $(Q_i')^c$ (or rather in $\mathcal{D}_\Theta^c$), summarizing the above discussion, we obtain that once conditioned on $\mathcal{F}$ and on the variable $U_\Theta$, 
the number $x_j$ of $A$-important points inside $Q_j$ is known. Hence one can write
\begin{align}
\Eb{(x_i - \beta y_i ) (x_j - \beta y_j)} & = \Eb{ (x_j - \beta y_j) \Eb{ x_i - \beta y_i \bigm| \mathcal{F}, U_{\Theta}}}, \nonumber
\end{align}
which allows us, by focusing on $\Eb{ x_i - \beta y_i \bigm| \mathcal{F}, U_{\Theta}}$, to ``decouple'' the behavior inside the $Q_i$ box from what is outside $Q_i'$.
Note that this strategy explains why our main coupling result 
in Section~\ref{s.coupling} dealt with this additional conditioning on $U_\Theta$.

But we still need to justify why one can assume that there are indeed 4 faces around $Q_i'$.
This is obviously not necessarily the case, especially if $Q_i'$ is of size comparable to $d(Q_i,Q_j)\ge r$. To make sure that with high conditional probability
(under $\{y_i=1\}$), there are only four faces around $Q_i'$, a solution is to fix the radius of our enlarged box $Q_i'$ at some intermediate scale $\gamma$ so that  
 $\eps\ll\gamma\ll r \ll 1$. Indeed, now it is unlikely (under $\{y_i =1\}$) that there will be 5 arms connecting $Q_i'$ and $Q_j$ (the five-arm event being much less likely than the 
 four-arm event).
\vskip 0.3 cm

The actual proof in the coming subsection will be rather technical: 
we will need to handle these various intermediate scales $\eps \ll \gamma \ll r \ll 1$ as well as some issues coming from the discrete lattice.
We will encounter some additional technicalities during the proof that we chose to hide in the above outline. Nevertheless, this sketch of proof should 
serve as a reading-guide for the next subsection.

\subsection{Proof of the main estimate ($L^2$-approximation)}


\proofof{Proposition \ref{l.XY}} In the following, a box $B(q,r)$ centered at $q$ (which will usually be the center of a tile) of radius $r$ will denote the union of tiles whose center is included inside $e^{i \theta} [-r,r)^2 + q$.
Let us consider in particular the square $Q_0 = B(a,\eps)$ in $G$ centered around $a$ (of course, $Q_0$ is not necessarily in $B\subset\Delta$).
Let $x_0$ be the number of sites inside $Q_0$ which are $(B(a,1)\setminus Q_0)$-important.
Following the notations of Section~\ref{s.coupling},
let $\mathcal{A}_0(2\eps,1)$ be the event that there are open arms from $B(a,2\eps)$ to the left and right sides of $\p B(a,1)$
and closed arms from $B(a,2\eps)$ to the lower and upper sides of $\p B(a,1)$
(here, the ``left side'' of $B(a,1)$ would be the image of the usual left side after multiplying by $e^{i \theta}$). Furthermore, let $U_{0}$ be the indicator function that there is a left-right crossing in $B(a,1)$.


We define
\begin{equation}\label{e.beta}
\beta:= \beta(\eta,\eps,a,\theta) =\Eb{ x_0 \bigm| \mathcal{A}_0(2\eps,1)}.
\end{equation}

We will start the proof by dividing the points in $\eta\Tg \cap B$ into squares in a more careful way than above, since there are some issues coming from the 
discrete lattice to deal with. 
Let $S_1,\ldots,S_p$ denote the list of $\eps$-squares in the grid $G$ which are contained in $B$.
For each square $S_i$, let $Q_i$ be the union of $\eta$-tiles whose center is included in $S_i$ (to avoid multiple allocation, each square
is considered to be a translate of $[-\eps,\eps)^2$, this convention being mapped by $z\mapsto e^{i\theta} z+ a$).
Therefore, the set of $\eta$-tiles inside $B$ is partitioned into the ``tiling'' $\{Q_i\}_i$ plus some boundary (or exterior) $\eta$-tiles which are at distance at most $2\eps$ from $\p B$. Let $Q_{\mathrm{ext}}$ be the set of these $\eta$-tiles.
Because of issues coming from the discrete lattice (mesh $\eta$),
we need to slightly change our definition of $Y$. For each $1\leq i \leq p$, let $q_i$ be the closest $\eta$-tile, in any reasonable sense, to the
center of $Q_i$, and define $y_i$ to be the indicator function of the event that there are 4 arms from $B(q_i, 2\eps)$ to $\p_2 A$.
Notice that the tiles $B(q_i, 2\eps)$ for $1\leq i \leq p$ are all translates of each other and are identical to $B(a,2\eps)$ which is used
in the definition of $\beta$. This will be relevant when we apply the coupling argument (Proposition~\ref{pr.couplingfaces}) to our
situation.
We now define $Y:= y_1+\ldots+y_p$. (Notice that, in the scaling limit, this  definition matches with our original one, which is the only relevant thing to us.) Similarly, for any $1\leq i \leq p$, let $x_i$ be the number of $A$-important points inside $Q_i$, and let $x_{\mathrm{ext}}$ be
the number of $A$-important points in $Q_{\mathrm{ext}}$. So, the number $X$ of $A$-important points in $B$ is
$X= x_1+\ldots+x_p+x_{\mathrm{ext}}$.

Let us recall at this point that there is some universal constant $C>0$ such  that $\Eb{Y^2} \le C \Eb{Y}^2$ and $\Eb{X^2} \le C \Eb{X}^2$
(these well-known second moment calculations can be found, for example, in \cite[Lemma 3.1]{\GPS1}). 
Furthermore, it will be helpful to keep in mind the following standard estimates:
\begin{equation}\label{e.stand}
\left\lbrace
\begin{array}{ccl} \beta &\asymp & \eps^2 \eta^{-2} \alpha_4^\eta(\eta,\eps)\\
\Eb{X} &\asymp& \eta^{-2} \alpha_4^\eta(\eta,1)\\
\Eb{X^2} &\asymp& \eta^{-4} \alpha_4^\eta(\eta,1)^2\,.\end{array}
\right.
\end{equation}

Let us now briefly explain why the additional term $x_{\mathrm{ext}}$ has a negligible contribution.
Notice that the estimate of $\Eb{x_{\mathrm{ext}}^2}$ is similar to the second moment of the number of points which are 1-important
inside a band of length 1 and width $\eps$. Writing this second moment as a sum over pairs of points, it is easy to check (using $\lim_{\eps\to 0} \eps^{-1} \alpha_4(\eps,1)= 0$, or the same for $\alpha_4^\eta$, uniformly in $\eta>0$) that in this case the main contribution comes from pairs of points which are about $\eps$-close (unlike the second
moment for the number of 1-important points in a square where the main contribution comes from ``distant'' points, because $\eps^{-2}\alpha_4(\eps,1)\to \infty$). Since we can divide the band into about $\eps^{-1}$ $\eps$-squares, each of which is 1-important with probability about $\alpha_4(\eps,1)$, we get
\begin{align*}
\Eb{x_{\mathrm{ext}}^2} &\leq O(1) \eps^{-1} (\eps/\eta)^{4} \alpha_4^\eta(\eta,\eps)^2 \alpha_4^\eta(\eps,1)\\
&\leq  O(1) \eps^3 \alpha_4^\eta(\eps,1)^{-1} \eta^{-4}\alpha_4^\eta(\eta,1)^2.
\end{align*}
Since $\eps^3 \alpha_4(\eps,1)^{-1} = o(1)$, when $\eps$ goes to 0, this implies $\Eb{x_{\mathrm{ext}}^2} = o(\Eb{X^2})$ by (\ref{e.stand}) (when $\eps$ goes to 0,
uniformly in $\eta\in (0,\eps)$). Another way to explain why $\Eb{x_{\mathrm{ext}}^2}$ is negligible (without going into the details) is to notice that $\eps^{-1} \alpha_4(\eps,1)\to 0$ implies that with probability going to 1 when $\eps$ goes to zero, $x_{\mathrm{ext}}=0$, and conditioning it to be non-zero, its second moment  should ``obviously'' be not much larger than the unconditioned $\Eb{X^2}$, since the pairwise correlations between possible 1-important points at a certain distance are similar in the two cases.

This allows us to restrict our attention to $\bar X = x_1+\ldots+x_p$, since
\begin{equation}
\|X-\beta Y \|_2 \leq \|\bar X -\beta Y\|_2 + \|x_{\mathrm{ext}}\|_2 \,,\nonumber
\end{equation}
and hence it is enough to prove that
\begin{equation}
\Eb{(\bar X - \beta Y)^2} = o(\Eb{ X^2}).\nonumber
\end{equation}

Now let us fix some $\delta>0$. We need to prove that for $\eps$ and $\eta/\eps$ small enough, we have
\begin{equation}\label{e.barX}
\Eb{(\bar X - \beta Y)^2} \leq \delta\, \Eb{X^2}.
\end{equation}
Let us write
\begin{equation}\label{e.correls}
\Eb{(\bar X - \beta Y)^2}= \sum_{i,j} \Eb{(x_i - \beta y_i)(x_j-\beta y_j)};
\end{equation}
so, in order to prove (\ref{e.barX}), we need to control the correlations between the number of $A$-important
points in squares $Q_i$ and $Q_j$. If the squares are close, then $x_i$ and $x_j$ are highly correlated, and some correlation is still there even if the squares are far away from each other, but at least in that case we will control their dependence well enough. So, if $r=r(\delta)$ is any distance that we will choose later depending on $\delta$, it will be
convenient to split the above sum (\ref{e.correls}) into a ``diagonal'' term corresponding to nearby squares, plus a
term corresponding to distant squares:
\begin{align}
\label{e.correls2}
\Eb{(\bar X - \beta Y)^2} = \sum_{d(Q_i,Q_j)\leq r} \Eb{(x_i - \beta y_i)(x_j-\beta y_j)} \\
+ \sum_{d(Q_i,Q_j)>r}
\Eb{(x_i - \beta y_i)(x_j-\beta y_j)}\,.\nonumber
\end{align}

First, we estimate from above the near-diagonal term. Take any $i,j$ such that $d(Q_i,Q_j)\leq r$. We want to bound from above
$\Eb{(x_i - \beta y_i)(x_j-\beta y_j)}$ (note that this might as well be negative, in which case it would ``help'' us). We have
\begin{equation}
\label{e.xixj}
\Eb{(x_i - \beta y_i)(x_j-\beta y_j)} \leq \Eb{x_i x_j +\beta^2 y_i y_j}
\end{equation}

We deal with the term $\sum_{d(Q_i,Q_j)\leq r} \Eb{x_i x_j}$, the other one being treated in a similar way.
There are $O(1) \eps^{-2}$ choices for the box $Q_i$ (where $O(1)$ depends on
$B\subset \Delta$). Choose one of the $Q_i$ boxes. For any $k\geq 0$ such that $2^k \eps < r $, there are $O(1) 2^{2k}$ boxes
$Q_j$ satisfying $2^k \eps\leq d(Q_i,Q_j)<2^{k+1} \eps$. For any of these boxes, we have (by the usual quasi-multiplicativity)
\begin{eqnarray}
\Eb{x_i x_j} & = & \sum_{x\in Q_i,\,y\in Q_j} \Pb{x,\,y \text{ are }\,A\text{-important}}\nonumber\\
&\leq & O(1) \eps^4 \eta^{-4} \frac{\alpha_4^\eta(\eta,1)^2}{\alpha_4^\eta(2^k \eps,1)}.\nonumber
\end{eqnarray}
This gives us
\begin{eqnarray}
\sum_{d(Q_i,Q_j)\leq r} \Eb{x_i x_j} &\leq& O(1) \eps^{-2} \sum_{k\leq \log_2(r/\eps)} 2^{2k} \alpha_4^\eta(2^k\eps,1)^{-1}
\eps^4\eta^{-4} \alpha_4^\eta(\eta,1)^2\nonumber\\
&\leq & O(1) r^2 \alpha_4^\eta(r,1)^{-1} \eta^{-4} \alpha_4^\eta(\eta,1)^2\nonumber\\
&\leq & O(1) r^2\alpha_4(r,1)^{-1} \Eb{X^2},\nonumber
\end{eqnarray}
where we used $r^2 \alpha_4(r,1)^{-1} = O(r^{\zeta})$ with some $\zeta>0$ as $r\to 0$  (uniformly in $\eta$) to get the second line. Using the same estimate $O(r^{\zeta})$, by choosing $r=r(\delta)$ small enough and by applying
the same argument to the other term of~(\ref{e.xixj}), we obtain that 
$$\sum_{d(Q_i,Q_j)\leq r}
\Eb{(x_i - \beta y_i)(x_j-\beta y_j)} \leq \frac{\delta}{2}\, \Eb{X^2}.$$
\vskip 0.3 cm

We now turn to the second term $\sum_{d(Q_i,Q_j)>r}
\Eb{(x_i - \beta y_i)(x_j-\beta y_j)}$ in (\ref{e.correls2}).
For the diagonal term, the strategy was to use the fact that there were few terms in the sum, and that each of the terms was of reasonable size; here we have many terms to deal with, so we need to proceed differently: we will prove that
if $\eps$ and $\eta/\eps$ are small enough, then for any $i,j$ such that  $l:=d(Q_i,Q_j)>r$ we have:
$$
\Eb{(x_i-\beta y_i)(x_j - \beta y_j)} \leq
\frac {\delta} 2 \, \Eb{x_i x_j}.
$$
Let $q_i$ and $q_j$ be the respective centers of these squares.
Let $\gamma  \in (2\eps,r/4)$ be some intermediate distance whose value will be fixed later.
Following the notations of Section~\ref{s.coupling}, let $\Upsilon$ be the set of all interfaces crossing
the annulus $B(q_i, l/2) \setminus B(q_i, \gamma)$. As previously, $\Qual(\Upsilon)$ will denote the least distance
between the endpoints of the interfaces on $\p B(q_i,\gamma)$ normalized by $\gamma$. Let $\mathcal{A}_4
=\mathcal{A}_4(\gamma, l/2)$ be the event that there are at least four arms of alternating colors in the annulus
$B(q_i, l/2) \setminus B(q_i, \gamma)$. Furthermore, let $\mathcal{A}_5$ be the event that there are at least five
arms in the same annulus, with four of them of alternating colors.
Recall that $\mathcal{G}=\mathcal{G}(\gamma,l/2)$ is the event that there are exactly four alternating arms (thus
$\mathcal{G}=\mathcal{A}_4 \setminus \mathcal{A}_5$), and $\mathcal{T}^{\alpha} =\mathcal{T}^{\alpha}(\gamma,l/2)$
is the event that $\{ \Qual(\Upsilon) >\alpha \}$, with the quality cut-off $\alpha$ to be fixed later. We now define the following disjoint events:
\begin{equation}
\left\lbrace\begin{array}{lcl}
\mathcal{W}_i & =& \mathcal{G} \cap \mathcal{T}^{\alpha}\nonumber\\
\mathcal{Z}_i &=& \mathcal{A}_5 \cup ( \mathcal{A}_4 \cap \neg\mathcal{T}^{\alpha} )\,.\nonumber
\end{array}\right.
\end{equation}
 In order for $(x_i-\beta y_i)(x_j - \beta y_j)$ to be non-zero, $\A_4$, hence either $\mathcal{W}_i$ or $\mathcal{Z}_i$,
 must hold. We can thus write
\begin{align}
\label{e.ZiWi}
\Eb{(x_i-\beta y_i)(x_j - \beta y_j)} &  \\
&  \hskip -3 cm = \Eb{(x_i-\beta y_i)(x_j - \beta y_j) 1_{\mathcal{Z}_i}} + \Eb{(x_i-\beta y_i)(x_j - \beta y_j) 1_{\mathcal{W}_i}}.\nonumber
\end{align}
We think of $\mathcal{W}_i$ as the ``good'' event, on which the coupling argument (Proposition~\ref{pr.couplingfaces}) can be used, while $\mathcal{Z}_i$ is the ``bad'' event that should have a negligible contribution by choosing $\alpha$ small enough.

We first want to bound the first term of (\ref{e.ZiWi}). By definition, $\Pb{\mathcal{Z}_i} \leq \Pb{\mathcal{A}_5} +
\Pb{\mathcal{A}_4 \cap \neg\mathcal{T}^{\alpha} }$.
For the second event, notice that
\begin{eqnarray}
\mathcal{A}_4 \cap \neg\mathcal{T}^{\alpha} &=& \mathcal{A}_4(\gamma,l/2) \cap \neg\mathcal{T}^{\alpha}(\gamma,l/2)\nonumber\\
&\subset& \mathcal{A}_4(2\gamma,l/2) \cap \neg\mathcal{T}^{\alpha}(\gamma, 2\gamma)\,,\nonumber
\end{eqnarray}
since there are more interfaces between radii $\gamma$ and $2\gamma$ than between $\gamma$ and $l/2$, therefore the
quality of the set of interfaces is smaller for the annulus $B(q_i,2\gamma) \setminus B(q_i,\gamma)$. 
It is known (see Lemma A.2 in \cite{\SchrammSteif}) that there is a function $h(\alpha)$ such that (uniformly in $\eta \in (0,\eps/10)$),
$\Pb{\neg \mathcal{T}^{\alpha}(\gamma, 2\gamma)} < h(\alpha)$, and furthermore $h(\alpha)$ goes to zero when $\alpha$
goes to zero. We then deduce (by independence on disjoint sets)
\begin{equation}
\label{e.Zibound}
\Pb{\mathcal{Z}_i} \leq  O(1) \alpha_5^\eta(\gamma, l) + O(1) h(\alpha) \alpha_4^\eta(\gamma , l).
\end{equation}
Furthermore, we have
\begin{eqnarray}
\Eb{(x_i-\beta y_i)(x_j - \beta y_j) 1_{\mathcal{Z}_i}} &\leq & \Eb{(x_i x_j + \beta^2 y_i y_j) 1_{\mathcal{Z}_i} } \nonumber \\
&\leq & \Pb{\mathcal{Z}_i} \Eb{x_i x_j +\beta^2 y_i y_j \bigm| \mathcal{Z}_i}. \nonumber
\end{eqnarray}
Let $\mathcal{A}_{i,j}$ be the event that there are four arms from $B(q_j,2\eps)$ to $B(q_j, l/2)$, four arms from
$B(q_i,2\eps)$ to $B(q_i, \gamma)$ and four arms from
$B(\frac {q_i+q_j} 2, l)$ to $\p_2 A$.
By independence on disjoint sets, we can write
\begin{eqnarray}
\Eb{x_i x_j + \beta^2 y_i y_j \bigm| \mathcal{Z}_i} &=& \Pb{\mathcal{A}_{i,j}\bigm| \mathcal{Z}_i}
\Eb{x_i x_j + \beta^2 y_i y_j \bigm| \mathcal{Z}_i, \mathcal{A}_{i,j}} \nonumber\\
&=& \Pb{\mathcal{A}_{i,j}} \Eb{ x_i x_j + \beta^2 y_i y_j \bigm| \mathcal{Z}_i, \mathcal{A}_{i,j}}\,.\nonumber
\end{eqnarray}
In order to have more independence, let us introduce the number $\tilde{x_i}$ of points in $Q_i$ which have
four arms to $\p B(q_i, 2\eps)$; it is clear that $x_i\leq \tilde{x_i}$. We define $\tilde{x_j}$ in the same way.
Then,
\begin{align}
\Eb{x_i x_j + \beta^2 y_i y_j \bigm| \mathcal{Z}_i} &\leq  \Pb{\mathcal{A}_{i,j}}
\Eb{\tilde{x_i}\tilde{x_j} +\beta^2 y_i y_j\bigm| \mathcal{Z}_i, \mathcal{A}_{i,j}} \nonumber\\
&\leq   O(1) \alpha_4^\eta(\eps, \gamma) \alpha_4^\eta(\eps, 1) (\Eb{\tilde{x_i} \tilde{x_j}} + \beta^2) \nonumber\\
&\leq  O(1) \alpha_4^\eta(\eps,\gamma) \alpha_4^\eta(\eps,1) (\eps/\eta)^{4} \alpha_4^\eta(\eta, \eps)^2. \label{e.Zicond}\\
&\leq O(1) \eps^4 \eta^{-4} \frac{\alpha_4^\eta(\eta,1)^2}{\alpha_4^\eta(\gamma,1)}\,\, \text{ by quasi-multiplicativity}.\nonumber
\end{align}
Combining (\ref{e.Zibound}) and (\ref{e.Zicond}) gives the following bound on the first term in~(\ref{e.ZiWi}):
\begin{eqnarray}
\Eb{(x_i-\beta y_i)(x_j - \beta y_j) 1_{\mathcal{Z}_i}} &\leq &
O(1) \left( \frac {\alpha_5^\eta(\gamma,l)}{\alpha_4^\eta(\gamma,l)}  + h(\alpha)  \right) \eps^4 \eta^{-4} \frac {\alpha_4^\eta(\eta,1)^2 }{\alpha_4^\eta(l,1)}.\nonumber
\end{eqnarray}
On the other hand, it is easy to check (using quasi-multiplicativity) that $\Eb{x_i x_j} \asymp \eps^4 \eta^{-4} \frac{\alpha_4^\eta(\eta,1)^2}{\alpha_4^\eta(l,1)}$,
hence
\begin{equation}
\label{e.termZi}
\Eb{(x_i-\beta y_i)(x_j - \beta y_j) 1_{\mathcal{Z}_i}} \leq O(1)
\left( \frac {\alpha_5^\eta(\gamma,l)}{\alpha_4^\eta(\gamma,l)}  + h(\alpha)  \right) \Eb{x_i x_j}.
\end{equation}

We need to bound now the second term $\Eb{(x_i-\beta y_i)(x_j - \beta y_j) 1_{\mathcal{W}_i}}$.
We will use for that purpose the coupling argument (Proposition \ref{pr.couplingfaces}). Indeed, on the
event $\mathcal{W}_i$, there are exactly 4 arms crossing the annulus $B(q_i, l/2) \setminus B(q_i, \gamma)$,
therefore there are exactly 4 interfaces crossing this annulus and, as we have seen in Section~\ref{s.coupling},
they induce a configuration of faces $\Theta=\{ \theta_1,\ldots,\theta_4\}$ at radius $\gamma$ around
$q_i$ (here $\theta_1,\theta_3$ are
the open faces).
As in Section~\ref{s.coupling}, let $\mathcal{D}_{\Theta}$ be the bounded component of $\C\setminus \Theta$
(which is a finite set of $\eta$-tiles) and let $U=U_{\Theta}$ be the indicator function that there is an open crossing
from $\theta_1$ to $\theta_3$ in $\mathcal{D}_{\Theta}$. Let $\mathcal{F}_{\Theta}$ be the $\sigma$-field generated by the tiles
in $\mathcal{D}_{\Theta}^c$. On the event $\mathcal{W}_i$, we may condition on $\mathcal{F}_{\Theta}$ in order to
``factorize'' the information in the $Q_i$ and $Q_j$ boxes, but notice that even if we condition on $\mathcal{F}_{\Theta}$
(and thus, in particular, we know all the information inside $Q_j$), the number $x_j$ of $A$-important points in $Q_j$
might still depend on the connectivities inside $Q_i$. That is why we also condition on $U_{\Theta}$ which gives the only information that is significant outside $\mathcal{D}_{\Theta}$ about what the connectivities are inside $\mathcal{D}_{\Theta}$.
We end up with
\begin{align*}
\Eb{(x_i-\beta y_i)(x_j - \beta y_j) 1_{\mathcal{W}_i}} &=
\Pb{\mathcal{W}_i} \Eb{(x_i-\beta y_i)(x_j - \beta y_j) \bigm| \mathcal{W}_i} \\
&= \Pb{\mathcal{W}_i} \Eb{ \Eb{(x_i-\beta y_i)(x_j - \beta y_j) \bigm| \mathcal{F}_{\Theta}, \, U_{\Theta}} \bigm| \mathcal{W}_i} \\
&= \Pb{\mathcal{W}_i} \Eb{ (x_j - \beta y_j) \Eb{x_i-\beta y_i\bigm| \mathcal{F}_{\Theta}, \, U_{\Theta}} \bigm| \mathcal{W}_i},
\end{align*}
since $x_j - \beta y_j$ is measurable with respect to the $\sigma$-field generated by $\mathcal{F}_{\Theta}$ and $U_{\Theta}$
(which allows us to ``factorize'' the $Q_i$ and $Q_j$ boxes). Therefore,
\begin{align}
\label{e.onWi}
\Big|\Eb{(x_i-\beta y_i)(x_j - \beta y_j) 1_{\mathcal{W}_i}}\Big| & \leq \\
& \hskip -3 cm \Pb{\mathcal{W}_i} \EB{ (x_j + \beta y_j) \big|\Eb{x_i-\beta y_i\bigm| \mathcal{F}_{\Theta}, \, U_{\Theta}}\big| \md \mathcal{W}_i}.\nonumber
\end{align}
As in Section \ref{s.coupling}, let $\mathcal{A}_{\Theta}=\mathcal{A}_{\Theta}(2\eps, \gamma)$ be the event that there are open arms from
$B(q_i, 2\eps)$ to the open faces $\theta_1, \theta_3$ and closed arms from $B(q_i,2\eps)$
to the closed faces $\theta_2,\theta_4$.
Let $\mathcal{X}_{\Theta}$ be the event that there are open arms in $\mathcal{D}^c_{\Theta}$ from $\theta_1, \theta_3$
to $\p_2 A$ and closed arms from $\theta_2,\theta_4$ to $\p_2 A$, so that on the event $\mathcal{W}_i$ we have
\begin{align}
\Eb{x_i-\beta y_i \bigm| \mathcal{F}_{\Theta}, U_{\Theta}}& = 1_{\mathcal{X}_{\Theta}}\,\Eb{x_i - \beta y_i \bigm|
\mathcal{F}_{\Theta}, U_{\Theta}}\nonumber\\
&= 1_{\mathcal{X}_{\Theta}}\,\Eb{x_i - \beta 1_{\mathcal{A}_{\Theta}} \bigm|
\mathcal{F}_{\Theta}, U_{\Theta}}\label{e.before} \\
&= 1_{\mathcal{X}_{\Theta}}\, \Pb{\mathcal{A}_{\Theta}\bigm| U_{\Theta}} \Eb{x_i - \beta \bigm|
\mathcal{F}_{\Theta}, \mathcal{A}_{\Theta}, U_{\Theta}} ,\nonumber
\end{align}
by independence on disjoint sets and by having $y_i=1_{\mathcal{A}_{\Theta}}$ on $\mathcal{W}_i \cap \mathcal{X}_{\Theta}$.

\vskip 0.3 cm

We now wish to bound $\big| \Eb{x_i - \beta \bigm|
\mathcal{F}_{\Theta}, \mathcal{A}_{\Theta}, U_{\Theta}}\big|$ using the coupling argument from Proposition~\ref{pr.couplingfaces}.
Recall that we defined $\beta:= \Eb{x_0 \bigm| \mathcal{A}_0(2\eps,1)}$. 
In order to apply Proposition~\ref{pr.couplingfaces}, it would be easier to deal with a modified version of $\beta$, namely:
\[
\hat \beta:= \Eb{x_0 \bigm| \mathcal{A}_0(2\eps,1), U_0=1}\,.
\]
Indeed, as we shall see below, the estimate $|\Eb{x_i - \hat \beta \bigm| \mathcal{F}_{\Theta}, \mathcal{A}_{\Theta}, U_{\Theta}}| = o(\hat \beta)$ follows easily from Proposition~\ref{pr.couplingfaces}. To connect this result with our goal, we will rely on the following extension of the coupling results of Section~\ref{s.coupling}:

\begin{lemma}\label{l.hatbeta}
As $\eps$ and $\eta/\eps$ go to zero, one has
\[
\beta = \hat \beta (1+o(1))\,.
\]
\end{lemma}

\proof By definition,
\begin{align}
\beta  = \Eb{x_0 \bigm| \mathcal{A}_0(2\eps , 1)} 
= &\; \Pb{U_0=1 \bigm| \mathcal{A}_0(2\eps,1)} \Eb{x_0 \bigm| \mathcal{A}_0(2\eps,1), U_0=1} \nonumber \\
&\; + \Pb{U_0=0 \bigm| \mathcal{A}_0(2\eps,1)} \Eb{x_0 \bigm| \mathcal{A}_0(2\eps,1), U_0=0}.\nonumber 
\end{align}

Let $\nu_0^{+}$ (resp.~$\nu_0^-$) be the law in $B(a,1)$ conditioned on the events $\mathcal{A}_0(2\eps,1)$ and $\{U_0=1\}$ (resp.~$\{U_0=0\}$).
Proposition~\ref{pr.couplingfaces} implies that with probability at least $1-(2\eps/\gamma)^k$, one can couple these two measures so that the induced faces 
at radius $2\eps$ are identical aside from their color which are reversed. When the coupling succeeds, even though the induced faces are of reversed color,  
it is clear (by switching the colors) that the conditional expected number of important points inside $Q_0=B(a,\eps)$ is the same for both measures
(here conditioned on the induced faces at radius $2\eps$, and on the event $\mathcal{S}$ that the coupling succeeds). This implies that 

\begin{align*}
\left|\Eb{x_0 \bigm| \mathcal{A}_0(2\eps,1), U_0=1}-\Eb{x_0 \bigm| \mathcal{A}_0(2\eps,1), U_0=0}\right|& \\
& \hskip -3 cm =  |\nu_0^+(x_0) - \nu_0^-(x_0)| \\
& \hskip -3 cm \le  \Bigl( \frac {2 \eps}{\gamma} \Bigr)^k \bigl| \nu_0^+(x_0 \bigm| \neg \mathcal{S}) - \nu_0^-(x_0 \bigm| \neg \mathcal{S})\bigr|\,. \\
& \hskip -3 cm \le  \Bigl( \frac {2 \eps}{\gamma} \Bigr)^k  \bigl( \nu_0^+(x_0 \bigm| \neg \mathcal{S}) + \nu_0^-(x_0 \bigm| \neg \mathcal{S}) \bigr) \,.
\end{align*}

Let us focus on the first term $\nu_0^+(x_0 \bigm| \neg \mathcal{S})$.
Call $\theta_\eps$ the configuration of faces around $B(a,2\eps)$ induced by
$\nu_0^+(\bullet \bigm| \neg \mathcal{S}) = \Eb{\bullet \bigm| \mathcal{A}_0(2\eps,1), U_0=1, \neg \mathcal{S}}$.
Let $U_\eps$ be the random variable equal to 1 if the open faces of $\theta_\eps$ are connected inside $\mathcal{D}_{\theta_\eps}$ and 0 else.
As such, conditioned on $\theta_\eps$, one has $U_0 = U_\eps$.
Furthermore, let $\mathcal{A}_\eps$ be the event that there are four arms from radius $\frac 3 2 \eps$ to radius 1, i.e., that the faces $\theta_\eps$ around $2\eps$ can be prolonged 
until radius $\frac 3 2 \eps$.
Now, if $\tilde{x_0}$ denotes the number of points in $Q_0$ which satisfy the four-arm event up to $\p B(a,\frac 3 2 \,\eps)$, 
it is clear that $x_0\le \tilde{x_0} \, 1_{\mathcal{A}_\eps}$. 

The advantage of $\tilde{x_0}$ over $x_0$ is that it is independent of both $\mathcal{A}_0(2\eps,1)$ and $\neg \mathcal{S}$. 
Nonetheless, it is not independent of  $\{U_0 = 1\}$, this is why the current argument requires some care.
To summarize, one has
\begin{align*}
\nu_0^+(x_0 \bigm| \neg \mathcal{S}) &\le \nu_0^+( \Eb{\tilde{x_0} \, 1_{\mathcal{A}_\eps} \bigm| \theta_\eps, U_\eps=1} \bigm| \neg \mathcal{S})\,.
\end{align*}

Now, for any configuration $\theta_\eps$ around $2\eps$ sampled according to $\nu_0^+(\bullet \bigm| \neg \mathcal{S})$, one has
\begin{align*}
\Eb{\tilde{x_0} \, 1_{\mathcal{A}_\eps} \bigm| \theta_\eps, U_\eps=1 } & \le \Eb{\tilde{x_0} \bigm| \theta_\eps, \mathcal{A}_\eps, U_\eps=1} \\
&\le \frac {1}{ \Pb{U_\eps =1 \bigm| \theta_\eps, \mathcal{A}_\eps}} \,\Eb{\tilde{x_0} \bigm| \theta_\eps, \mathcal{A}_\eps} \\
&= \frac {1}{ \Pb{U_\eps =1 \bigm| \theta_\eps, \mathcal{A}_\eps}} \,\Eb{\tilde{x_0}} \; \begin{array}{l} \text{by independence} \\ \text{on disjoint sets} \end{array}\\
&\asymp \frac {1}{ \Pb{U_\eps =1 \bigm| \theta_\eps, \mathcal{A}_\eps}} \, \beta\,.
\end{align*}
Hence, it remains to bound $\Pb{U_\eps = 1 \bigm| \theta_\eps, \mathcal{A}_\eps}$ from below. Whatever the separation quality $\Qual(\theta_\eps)$ at radius $2\eps$ is, the Separation Lemma~\ref{l.separation} says that, conditioned on $\mathcal{A}_\eps$, the interfaces around radius $\frac 3 2  \eps$ are well-separated ($\Qual > 1/4$) with positive conditional probability, and then, by the usual RSW gluing, $U_\eps=1$ also has a uniformly positive conditional probability.

Altogether (since the same argument works for $\nu_0^-(\bullet \bigm| \neg \mathcal{S})$), one ends up with 
\[
\Eb{x_0 \bigm| \mathcal{A}_0(2\eps,1), U_0=0} = \hat \beta + O(1)\Big(\frac{2\eps}{\gamma}\Big)^k \beta\,,
\]
which proves Lemma \ref{l.hatbeta}. \qed

\begin{remark}
Note that we used here Proposition~\ref{pr.couplingfaces} and Lemma~\ref{l.separation} in a slightly more general form, with boxes that are rotated by an angle $\theta$.
It is clear that the original proofs apply also here, with the exponent $k$ being independent of the direction $e^{i\theta}$.
\end{remark}
\vskip 0.3 cm

Recall from (\ref{e.before}) that we needed to bound $\big|\Eb{x_i - \beta \bigm|
\mathcal{F}_{\Theta}, \mathcal{A}_{\Theta}, U_{\Theta}}\big|$. Lemma \ref{l.hatbeta} allows us to focus on 
$\big|\Eb{x_i - \hat \beta \bigm|
\mathcal{F}_{\Theta}, \mathcal{A}_{\Theta}, U_{\Theta}}\big|$ instead.


Before applying the coupling argument from Proposition~\ref{pr.couplingfaces},  we need to be careful with the issues coming from the discrete lattice: indeed, $x_i$ is the number of points in $Q_i$ which are $A$-important, but $\hat \beta$ is defined as a (conditional) expected number of points in $Q_0=B(a,\eps)$, or, by translation invariance,
in $B(q_i,\eps)$. However, $Q_i$ and $B(q_i,\eps)$ do not exactly coincide (at the boundary points). Hence let us introduce $\hat x_i$ to be the number of $A$-important points in $B(q_i,\eps)$. We have
$$
x_i = \hat x_i + \sum_{x\in Q_i\setminus B(q_i,\eps)} 1_{\{x\text{ is $A$-important}\}} - \sum_{y\in B(q_i,\eps)\setminus Q_i}
 1_{\{y\text{ is $A$-important}\}}.
$$

There are $O(1) \eps \eta^{-1}$ such boundary points, each of them on the event
$\mathcal{W}_i \cap \mathcal{X}_{\Theta}$ and conditioned on $(\mathcal{A}_{\Theta},U_{\Theta})$ are $A$-important
with probability of order $O(1) \alpha_4^\eta(\eta, \eps)$. Hence,
on the event $\mathcal{W}_i \cap \mathcal{X}_{\Theta}$, we have
\begin{equation}
\label{e.boundaryQi}
\Eb{x_i | \mathcal{F}_{\Theta}, \mathcal{A}_{\Theta}, U_{\Theta}} = \Eb{\hat x_i |
\mathcal{F}_{\Theta}, \mathcal{A}_{\Theta}, U_{\Theta}} + O(1) \eps \eta^{-1}\alpha_4^\eta(\eta,\eps).
\end{equation}

In order to apply Proposition \ref{pr.couplingfaces}, one needs to consider both cases $U_{\Theta}=1$
and $U_{\Theta}=0$. On the event $\mathcal{W}_i \cap \mathcal{X}_{\Theta}$, if $\{ U_{\Theta} =1\}$ holds, we have
$$
\Eb{\hat x_i - \hat \beta \bigm| \mathcal{F}_{\Theta}, \mathcal{A}_{\Theta}, U_{\Theta}}  = \Eb{\hat x_i \bigm|
\mathcal{F}_{\Theta}, \mathcal{A}_{\Theta}, U_{\Theta} =1} - \Eb{ x_0 \bigm| \mathcal{A}_0(2\eps,1),\, U_0=1}.
$$
Proposition \ref{pr.couplingfaces} says that one can couple the probability measure conditioned on $\mathcal{A}_{\Theta}$ and $\{U_{\Theta}=1\}$
with the probability measure conditioned on $\mathcal{A}_0(2\eps, 1)$ and $\{U_0=1\}$ so that with probability at least
$1- (2\eps/ \gamma)^{k}$, we have $\hat x_i = x_0$. Let, as in Section~\ref{s.coupling}, $\mathcal{S}$ be the
event that the coupling succeeds. Hence, on the event $\mathcal{W}_i\cap\mathcal{X}_{\Theta}\cap \{U_{\Theta}=1\}$,
\begin{align}
\label{e.hatxi}
|\Eb{\hat x_i -\hat \beta \bigm| \mathcal{F}_{\Theta}, \mathcal{A}_{\Theta}, U_{\Theta}=1}| &   \\
&  \hskip -5 cm \leq  \left(\frac {2\eps}  {\gamma}\right)^{k}
\big(\Eb{\hat x_i  \bigm| \mathcal{F}_{\Theta}, \mathcal{A}_{\Theta}, U_{\Theta}=1, \neg \mathcal{S}}
+ \Eb{ x_0 \bigm| \mathcal{A}_0(2\eps, \gamma), U_0=1, \neg \mathcal{S}} \big) .\nonumber
\end{align}

Arguing exactly as in the proof of Lemma~\ref{l.hatbeta} (i.e., by introducing $\tilde{x_i}$, $\tilde{x_0}$, $\mathcal{A}_\eps$, and so on),
the inequality~\eqref{e.hatxi} becomes 
\begin{align}
\label{e.hatxi1}
|\Eb{\hat x_i - \hat \beta \bigm| \mathcal{F}_{\Theta}, \mathcal{A}_{\Theta}, U_{\Theta}=1}|& \le O(1)  \left(\frac {2\eps}  {\gamma}\right)^{k} (\eps/\eta)^2 \alpha_4^\eta(\eta,\eps) \nonumber \\
& \le O(1) \left(\frac {2\eps}  {\gamma}\right)^{k}\beta\,.
\end{align}

%

Now, on the event $\mathcal{W}_i \cap \mathcal{X}_{\Theta}$, if $\{ U_{\Theta} =0\}$ holds instead, we obtain similarly the 
following bound:
\begin{equation}
\label{e.hatxi2}
|\Eb{\hat x_i - \hat \beta \bigm| \mathcal{F}_{\Theta}, \mathcal{A}_{\Theta}, U_{\Theta}=0}|
\leq O(1) \left(\frac {2\eps}  {\gamma}\right)^{k} \beta \,.
\end{equation}
To obtain this bound, one can rely on the color-switching aspect of Proposition~\ref{pr.couplingfaces},
just as $\Eb{x_0 \bigm| \mathcal{A}_0(2\eps,1), U_0=0} = \hat \beta + o(\beta)$ was proved in Lemma~\ref{l.hatbeta}.
\vskip 0.3 cm

Summarizing: on the event $\mathcal{W}_i$, we have rewritten~(\ref{e.before}) as
\begin{align}
\label{e.summary}
\Eb{x_i-\beta y_i \bigm| \mathcal{F}_{\Theta}, U_{\Theta}} &  \\
&  \hskip -3 cm =
\left\lbrace\begin{array}{l} 1_{\mathcal{X}_{\Theta}, U_{\Theta}=1} \,\Pb{\mathcal{A}_{\Theta}\bigm| U_{\Theta}=1} (\Eb{x_i |
\mathcal{F}_{\Theta}, \mathcal{A}_{\Theta},U_{\Theta}=1} - \beta)  \\
 +\; 1_{\mathcal{X}_{\Theta}, U_{\Theta}=0} \,\Pb{\mathcal{A}_{\Theta}\bigm| U_{\Theta}=0} (\Eb{x_i |
\mathcal{F}_{\Theta}, \mathcal{A}_{\Theta},U_{\Theta}=0} - \beta) \,,\end{array}\right. \nonumber
\end{align}
and have bounded its different factors. The last ingredient is that
\begin{equation}\label{e.ATheta}
\Pb{\mathcal{A}_{\Theta} \bigm|
U_{\Theta}}< O_{\alpha}(1)\alpha_4^\eta(\eps,\gamma)\,,
\end{equation}
where the constant $O_\alpha(1)$ depends on our separation threshold $\alpha$ (to be fixed shortly).
In order to see why such an upper bound holds, one can write, for any $i\in\{0,1\}$,
\begin{align}
\Pb{\mathcal{A}_\Theta \bigm| U_\theta = i} & \le \frac {1}{ \Pb{U_\Theta = i}} \,\Pb{\mathcal{A}_\Theta} \nonumber \\
& \le \frac{O(1)}{\Pb{U_\Theta = i}}\, \alpha_4^\eta(\eps, \gamma)\,.\nonumber
\end{align}
Now, $\Qual(\Theta)>\alpha$, since we are on the event $\mathcal{W}_i$. Hence it is clear by RSW gluing (as in Figure~\ref{f.gluing} of Appendix~\ref{s.appendix}) that there is some constant $c(\alpha)>0$
such that for any $i\in\{0,1\}$, one has $\Pb{U_\Theta = i}>c(\alpha)$. This proves~\eqref{e.ATheta}.
\vskip 0.3 cm

Therefore, plugging (\ref{e.boundaryQi}), (\ref{e.hatxi1}, \ref{e.hatxi2}), Lemma \ref{l.hatbeta}, and~\eqref{e.ATheta}  into (\ref{e.summary}) gives (still on the event
$\mathcal{W}_i$):
\begin{align*}
\big|\Eb{x_i - \beta y_i \bigm| \mathcal{F}_{\Theta}, U_{\Theta}}\big| &\leq 1_{\mathcal{X}_{\Theta}} O_{\alpha}(1) \alpha_4^\eta(\eps,\gamma)
\left( \left(\frac{2\eps}{\gamma}\right)^{k} \eps^2 \eta^{-2}
+ \eps \eta^{-1}  \right) \alpha_4^\eta(\eta,\eps) \\
&\leq O_{\alpha}(1) \eps^2 \eta^{-2} \alpha_4^\eta(\eta,\gamma) \left( (2\eps/\gamma)^{k} + \frac {\eta}{\eps} \right) \; \begin{array}{l} \text{by quasimulti-} \\ \text{plicativity}\,.\end{array}
\end{align*}
(Here we used the quantitative estimate that $|\beta - \hat \beta| \le O(1) \left( \frac {2\eps} {\gamma} \right)^{k} \beta$, which follows from the proof of Lemma~\ref{l.hatbeta}).
\vskip 0.3 cm

It is clear that $\Pb{\mathcal{W}_i} \leq O(1) \alpha_4^\eta(\gamma, l)$, hence~(\ref{e.onWi}) becomes
\begin{eqnarray}
|\Eb{(x_i-\beta y_i)(x_j - \beta y_j) 1_{\mathcal{W}_i}}|  & & \nonumber \\
& & \hskip - 3.5 cm \leq \Pb{\mathcal{W}_i} \Eb{ (x_j + \beta y_j) \left|\Eb{x_i-\beta y_i\bigm| \mathcal{F}_{\Theta}, \, U_{\Theta}}\right| \bigm|
\mathcal{W}_i} \nonumber \\
& & \hskip - 3.5 cm \leq  O_{\alpha}(1) \eps^2 \eta^{-2} \alpha_4^\eta(\eta,l) \Eb{x_j+\beta y_j \bigm| \mathcal{W}_i} \left( (2\eps/\gamma)^{k} + \frac {\eta}{\eps} \right). \nonumber
\end{eqnarray}

In order to bound $\Eb{x_j | \mathcal{W}_i}$, we introduce $x^*_j$, the number of points in $Q_j$ which have four arms
to $\p B(q_j, l/2)$, and we let $\mathcal{G}^*$ be the event that there are four arms from $B(\frac{q_i+q_j} 2, l)$ to $\p_2 A$. By definition,
$x_j \leq x^*_j 1_{\mathcal{G}^*}$, therefore, by independence on disjoints sets ($x^*_j$ and $\mathcal{G}^*$ do not depend on $\mathcal{W}_i$),
we obtain $\Es{x_j \md \mathcal{W}_i} \leq O(1) \eps^2 \eta^{-2} \alpha_4^\eta(\eta,1)$. Similar estimates for $y_i$ lead to
$\Eb{x_j+\beta y_j \bigm| \mathcal{W}_i} \leq O(1) \eps^2 \eta^{-2} \alpha_4^\eta(\eta,1)$, hence
\begin{eqnarray}
|\Eb{(x_i-\beta y_i)(x_j - \beta y_j) 1_{\mathcal{W}_i}}|  &\leq & O_{\alpha}(1)
\eps^4 \eta^{-4} \frac{\alpha_4^\eta(\eta,1)^2}{\alpha_4^\eta(l,1)} \left( (2\eps/\gamma)^{k} + \frac {\eta}{\eps} \right)
\nonumber \\
&\leq & C(\alpha) \Eb{x_i x_j} \left( (2\eps/\gamma)^{k} + \frac {\eta}{\eps} \right), \nonumber
\end{eqnarray}
for some fixed constant $C(\alpha)>0$.
Now plugging this last expression (the ``good'' $\mathcal{W}_i$-part), together with the ``bad'' $\mathcal{Z}_i$-part'' (\ref{e.termZi}) into (\ref{e.ZiWi}) gives
\begin{eqnarray}
\Eb{(x_i-\beta y_i)(x_j - \beta y_j)} \nonumber\\
& & \hskip -3cm \leq  O(1) \Eb{x_i x_j} \left( \frac{\alpha_5^\eta(\gamma,l)}{\alpha_4^\eta(\gamma,l)}
+ h(\alpha) +  C(\alpha) \left(\frac{2\eps}{\gamma}\right)^{k}+ C(\alpha) \frac{\eta}{\eps} \right).\nonumber
\end{eqnarray}
Let us now fix the ``quality threshold'' $\alpha=\alpha(\delta)$ so that $O(1) h(\alpha) < \delta /8 $.
Recall we have already fixed $r=r(\delta)$ so that the diagonal term was less than $\delta/2 \, \Eb{X^2}$. So we have
$\gamma< r(\delta) < l$. It is a standard fact (proved by the BK inequality) that there is some $d>0$ such that, for any $\gamma<l$,
$\frac{\alpha_5^\eta(\gamma,l)}{\alpha_4^\eta(\gamma,l)} <  (\gamma/l)^d$. So, we fix $\gamma$ so that $O(1)
\frac{\alpha_5^\eta(\gamma,l)}{\alpha_4^\eta(\gamma,l)} <  \delta /8$. Now, by taking $\eps$ and $\eta/\eps$ small enough,
one obtains $\Eb{(x_i-\beta y_i)(x_j - \beta y_j)} \leq  \delta/2 \, \Eb{x_i x_j}$, which
ends the proof of Proposition~\ref{l.XY}.
\QED

\subsection{Ratio limit theorems}\label{ss.ratio}

In this subsection, we state and prove some limit theorems that will be used in the sequel and are also interesting in their own right. Their relevance to the existence and conformal covariance of the limit of the measures $\mu^A_\eta$ will be explained at the beginning of Subsection~\ref{ss.limitpivomeas}.

$\A_\square^{a,\theta}(r,R)$ will denote the four-arm event from the rotated square $B(a,r) = a+ e^{i\theta}[-r,r)^2$ to the rotated square $\p B(a,R)$ with the additional requirement that open arms are connected to the (rotated) right-left boundaries and closed arms are connected to the upper-lower boundaries. The probability of this event is denoted by $\alpha_\square^{\eta,a,\theta}(r,R)$. By $\A^{a,\theta}_\square(\eta,R)$ and $\alpha_\square^{\eta,a,\theta}(\eta,R)$, we mean the four-arm event from the hexagon containing $a$. Events and probabilities with the subscript $4$ always refer to the standard four-arm events between the standard unrotated squares. Finally, omitting the superscript $\eta$ from an $\alpha$-probability means it is understood in the scaling limit.

\begin{proposition}\label{pr.ratio}
For any fixed $r>0$,
\begin{equation}\label{e.limitlambda}
\lim_{\eta\rightarrow 0}\frac{\alpha_4^\eta(\eta,r)}{ \alpha_4^\eta(\eta, 1)} = \lim_{\eps\to 0} \frac{\alpha_4(\eps,r)}{\alpha_4(\eps,1)} = r^{-5/4}\,.
\end{equation}
Furthermore, there is an absolute constant $C>0$, such that, uniformly in the shift and orientation, $a, e^{i\theta}$,
\begin{equation}\label{e.ratio}
\lim_{\eta\to 0} \frac{\alpha_\square^{\eta,a,\theta}(\eta,1)}{\alpha_4^\eta(\eta,1)} =
\lim_{\eps\to 0} \frac{\alpha_\square^{a,\theta}(\eps,1)}{\alpha_4(\eps,1)} = 
C\,.
\end{equation}
\end{proposition}

\begin{remark}\label{r.ratioexplain}
This proposition might appear obvious at the first sight, and highly non-trivial at the second, but the truth is in between. The issue is that although we know the critical exponent 5/4, the probability $\alpha_4^\eta(\eta, 1)$ is only known to be $\eta^{5/4+o(1)}$, so there could be large but sub-polynomial factors, while the scaling limit quantity $\alpha_4(\eps,1)$ is known to be $\eps^{5/4}$ only up to constant factors. Therefore, the proposition is not a direct consequence of the determination of the exponent, and, to our knowledge, is a new result in itself. On the other hand, we do not have to improve the exponent results to prove our results: somehow, the possible sub-polynomial corrections cancel out in the quotients under interest.
\end{remark}

\proof Let us start by proving~(\ref{e.ratio}), using an extension of the coupling arguments from Section~\ref{s.coupling}. 
In order to show that the limit in $\eta$ exists, we will prove that the sequence $\big\{\alpha_\square(\eta,1)/ \alpha_4(\eta,1)\big\}_{\eta>0}$ satisfies the Cauchy criterion. Fix an intermediate radius $\gamma$ with $10 \eta < \gamma < \frac 1 {10}$, and let $\mathcal{A}_{\circ\square}^{a,\theta}(\gamma,1)$ be the event that there are four arms from the {\bf circle} of radius $\gamma$ around $a$ to the four arcs of the {\bf rotated box} $B(a,1)= a + e^{i\theta}[-1,1)^2$. Also, let $\mathcal{A}_{\circ4}(\gamma,1)$ be the event that there are four arms from the circle of radius $\gamma$ around $a$ to the {\bf un-rotated box} $a+ [-1,1)^2$ (without prescribed arcs).  We now write
\begin{equation}\label{e.circbox}
\frac{\alpha_\square^{\eta,a,\theta}(\eta,1)}{\alpha^\eta_4(\eta,1)} = 
\frac{ \P\big[ \mathcal{A}_\square^{a,\theta}(\eta,1) \bigm| \mathcal{A}_{\circ\square}^{a,\theta}(\gamma,1) \big] } {\P\big[ \mathcal{A}_4(\eta,1) \bigm| \mathcal{A}_{\circ 4}(\gamma,1) \big]} \, 
\frac {\alpha_{\circ\square}^{\eta,a,\theta}(\gamma,1)} {\alpha_{\circ 4}^\eta(\gamma,1)}. 
\end{equation}
As the mesh $\eta$ goes to zero, since we are on the triangular lattice, we can use the existence and conformal invariance of the scaling limit between the fixed radii $\gamma$ and 1 to get that 
$$
\lim_{\eta\to 0} \frac{\alpha_{\circ\square}^{\eta,a,\theta}(\gamma,1) }{ \alpha_{\circ 4}^\eta(\gamma,1)} = 
\frac{ \alpha_{\circ\square}^\theta(\gamma,1)}{\alpha_{\circ 4}(\gamma,1)} =: C(\gamma)
$$
exists and is independent of $\theta$. On the other hand, the first factor on the right hand side of~(\ref{e.circbox}) can be handled using the coupling results from Section~\ref{s.coupling}:
$$
\frac{ \P\big[ \mathcal{A}_\square^{a,\theta}(\eta,1) \bigm| \mathcal{A}_{\circ\square}^{a,\theta}(\gamma,1) \big] } {\P\big[ \mathcal{A}_4(\eta,1) \bigm| \mathcal{A}_{\circ 4}(\gamma,1) \big]} = 1 + O(\gamma^k),
$$
uniformly in $\eta < \gamma/10$, since the inside circles in the numerator and the denominator are identical from the discrete lattice perspective. In fact, we need here the coupling results in a slightly generalized way: firstly, the ``targets'' are circles now, not squares; secondly, we have here $\alpha_4$ together with $\alpha_\square$. The first issue requires a trivial modification; the second issue can be handled in a way very similar to the proof of Lemma~\ref{l.hatbeta}, hence the details are left to the reader. 

Altogether, (\ref{e.circbox}) becomes
\begin{equation}\label{e.circboxlim}
\frac{\alpha_\square^{\eta,a,\theta}(\eta,1)}{\alpha^\eta_4(\eta,1)} =(1+o(1))\, C(\gamma)\,,
\end{equation}
where $o(1)$ goes to zero as $\eta$ and $\gamma$ go to zero in such a way that $\eta < \gamma/10$. This proves the needed Cauchy criterion as $\eta\to 0$. Moreover, the resulting limit $C^{a,\theta}$ must equal $\lim_{\gamma\to 0} C(\gamma)$ (which thus needs to exist), implying that $C^{a,\theta}$ is in fact independent of both $a$ and $\theta$. The finiteness and positivity of this limit $C$ is clear from quasi-multiplicativity arguments.

For the $\eps\to 0$ limit in~(\ref{e.ratio}), we can do the same proof, arriving at the right hand side of~(\ref{e.circboxlim}) again, which shows that we get the same limit $C$. 
\medskip

The proof for the existence of the limits in~(\ref{e.limitlambda}) is almost identical to the above, except that we can simply use the coupling result of Proposition~\ref{pr.coupling2}. That is, assuming $r<1$ (the case $r \geq 1$ being symmetric), for $10 \eta < \gamma < \frac r {10}$,
\begin{align*}
\frac{\alpha_4^\eta(\eta,r)}{\alpha_4^\eta(\eta,1)} &= \frac{ \P\big[\mathcal{A}_4(\eta,r) \bigm| \mathcal{A}_4(\gamma,r)\big]}{\P\big[\mathcal{A}_4(\eta,1) \bigm|
\mathcal{A}_4(\gamma,1)\big]} \frac {\alpha_4^{\eta}(\gamma, r)}{ \alpha_4^{\eta}(\gamma,1)}  \\
&= \big(1+ O(1)\,(\gamma/ r)^k \big)\, \big(1+o(1)\big)\, \frac{ \alpha_4(\gamma, r)}{\alpha_4(\gamma,1)}\,,
\end{align*}
where the $O(1)$ term is uniform in $\eta$, while the $o(1)$ factor goes to 0 as $\eta$ goes to zero, not necessarily uniformly in $\gamma$. This is enough for the Cauchy criterion for $\big\{\alpha_4^\eta(\eta,r)/\alpha_4^\eta(\eta,1)\big\}_{\eta>0}$.
%
So, we have some limit $\ell$ (which is easily seen
to be positive and finite), which is also the limit of the ratios of the macroscopic probabilities $\alpha_4(\gamma,r) /\alpha_4(\gamma,1)$ when $\gamma\to 0$. It is easier to work with events at the scaling limit (since at the scaling limit we can use scale invariance), so we identify $\ell$ using the $\gamma\to 0$ limit.

On the triangular grid it is known that
\begin{equation}
\lim_{n \rightarrow \infty} \frac {\log \alpha_4(r^n, 1)}{n} =  \log(r^{5/4}).
\label{e.5/4}
\end{equation}

But one can write $\alpha_4(r^n,1)$ in the following way:
$$
\alpha_4(r^n,1) = \frac {\alpha_4(r^n,1)}{\alpha_4(r^n,r)} \frac{\alpha_4(r^{n-1},1)}{\alpha_4(r^{n-1},r)}\ldots
\frac{\alpha_4(r,1)}{1}\,.
$$
Therefore,
\begin{equation}
 \label{e.cesarol}
 \frac {\log \alpha_4(r^n,1)}{n} = \frac 1 n \sum_{j=1}^n \log \frac{\alpha_4(r^j,1)}{\alpha_4(r^j,r)}\,.
 \end{equation}
But, since $r^j$ goes to zero with $j$, we have that $\lim_{j\rightarrow \infty} \log\frac{\alpha_4(r^j,1)}{\alpha_4(r^j,r)} = \log \frac 1 \ell $.
By the convergence of the Ces\`aro mean, the right hand side of~(\ref{e.cesarol}) converges to $\log \frac 1 \ell$, hence comparing with~(\ref{e.5/4}) gives that $\log \frac 1 \ell = \log(r^{5/4})$, which concludes the proof. \qed
\medskip

Of course, similar results hold for the $\alpha_1$ and $\alpha_2$ probabilities, and we will prove another somewhat similar limit theorem, the rotational invariance of the two-point function $\Pb{x \mathop{\llra}\limits^{\omega_\eta} y}$. See Section~\ref{s.others}.

\subsection{The limit of the counting measures}\label{ss.limitpivomeas}

Proposition \ref{l.XY} shows that for any ball $B\subset \Delta$, the number $X$ of $A$-important points in $\eta\Tg\cap B$ is well approximated 
(in the $L^2$ sense) by the macroscopic quantity $\beta Y$. However, instead of $X$, what we really want to control is the measure of the ball  $B$, i.e.,
$\mu_\eta^A(B):= \frac X {\eta^{-2} \alpha_4^\eta(\eta,1)}$: we want to approximate it by some macroscopic quantity, and want to show its convergence as $\eta\to 0$. We will also want to prove scale and rotational invariance properties for the limit measure. 

Recall that the quantity $\beta$, which we have been working with so far, is defined as
\[
\beta= \beta(\eta,\eps,a,\theta):= \E\big[ x_0 \bigm| \mathcal{A}_0(2\eps, 1)\big]\,,
\]
where $\mathcal{A}_0(2\eps,1)=\A_\square^{a,\theta}(2\eps,1)$, as defined in the previous subsection. We have already used several times that $\alpha_\square^{\eta,a,\theta}(2\eps,1)$ is the same,  up to constant factors,  as $\alpha_4(\eps,1)$, a slightly more canonical quantity that does not depend on some given quad and some angle $\theta$; however, to achieve the goals we summarized above, we need to be more precise than ``up to constant factors''. This is the need fulfilled by Proposition~\ref{pr.ratio} above.

So, based on Proposition~\ref{l.XY} and Proposition~\ref{pr.ratio}, we shall first of all prove the following:

\begin{proposition}\label{pr.muB}
There exists an absolute constant $c>0$ (in fact, it is $4\cdot 2^{-5/4}$) such that for any $B\subset \Delta$
\begin{equation}
\EB{\bigl(\mu_\eta^A(B) - \frac c {\eps^{-2} \alpha_4(\eps,1)} Y \bigr)^2} \le \zeta(\eta,\eps)\,,
\end{equation}
where $\zeta(\eta,\eps)$ goes to zero when $\eps$ and $\eta/\eps$ go to zero, and $\alpha_4(\eps,1)$ may equally mean the four-arm probability on $\eta\Tg$ or in the scaling limit. As in Proposition~\ref{l.XY}, this control is independent of the chosen grid for $Y=Y_\eta^{\eps,a,\theta}$, but might depend on the ball $B$.
\end{proposition}

The proposition easily follows from Proposition~\ref{l.XY} plus the following claim:
\begin{lemma}
\label{l.betalimit}
For $c=4\cdot 2^{-5/4}$, uniformly in the orientation $a, e^{i\theta}$,
\[
\frac {\beta}{\eta^{-2} \alpha_4^\eta(\eta,1)} = \frac {\beta(\eta,\eps,a,\theta)} {\eta^{-2} \alpha_4^\eta(\eta,1)} = (1+o(1)) \frac c {\eps^{-2} \alpha_4(\eps,1)}\,,
\]
as $\eps$ and $\eta/\eps$ go to zero.
Recall here that the four-arm probabilities $\alpha_4(\eta,1)$ and $\alpha_4(\eps,1)$ are independent of the orientation of the grid $G$.
\end{lemma}

\proofof{Proposition \ref{pr.muB} assuming Lemma~\ref{l.betalimit}} 
By the triangle inequality, one has
\begin{align}
\Big\| \mu_\eta^A(B) - \frac {c} {\eps^{-2} \alpha_4(\eps,1)} Y \Big\|_2 & \nonumber \\
& \hskip -2.2 cm \le \Big\| \mu_\eta^A(B) - \frac{\beta}{\eta^{-2} \alpha_4^\eta(\eta,1) } Y\Big\|_2 + \Big\|\frac {\beta}{\eta^{-2}\alpha_4^\eta(\eta,1)} Y - \frac c {\eps^{-2}\alpha_4(\eps,1)} Y\Big\|_2 \nonumber\\
& \hskip -2.2 cm = \frac {\Eb{(X-\beta Y)^2}^{1/2}}{\eta^{-2} \alpha_4^\eta(\eta,1)} + \Bigl|\frac {\beta} {\eta^{-2} \alpha_4^\eta(\eta,1)} - \frac c {\eps^{-2} \alpha_4(\eps,1)} \Bigr| \, \Eb{Y^2}^{1/2} \nonumber \\
& \hskip -2.2 cm \le O(1) \xi(\eta,\eps)^{1/2} + o(1)\,, \nonumber
\end{align}
where, in the last line, the first term comes from Proposition~\ref{l.XY} and $\Eb{X}\asymp \eta^{-2} \alpha_4^\eta(\eta,1)$, while the second term (the $o(1)$) comes from Lemma~\ref{l.betalimit} and $\Eb{Y}\asymp \eps^{-2} \alpha_4(\eps,1)$.
So, assuming we have a uniform control over $a,\theta$ on the $o(1)$ in the above lemma, this indeed proves Proposition~\ref{pr.muB}.\qed

\proofof{Lemma \ref{l.betalimit}}
%
For any quad $\Quad\subset \C$ and any $\eta$-tile $x\in\Quad$, we denote by $\mathcal{A}_\square(x,\p \Quad)$ the event of having four alternating arms from $x$ to the four prescribed arcs of $\p \Quad$. So, one may rewrite $\beta$ as 
\begin{equation}
\beta  = \frac {\sum_{x\in \eta \Tg \cap B(a,\eps)}  \Pb{\mathcal{A}_\square(x, \p B(a,1))}  }{  \Pb{\mathcal{A}_0(2\eps, 1)}} \nonumber 
\end{equation}

Now, it is not hard to check that for all $x\in B(a,\eps)$, one has 
\begin{align}
\Pb{\mathcal{A}_\square (x, \p B(a,1))} &= \Pb{\mathcal{A}_\square(x, \p B(a,1))}  \nn \\
&= \Pb{\mathcal{A}_\square(x, \p B(x,1))} (1+o(1)) \nn \\
&= \Pb{\mathcal{A}_\square(a, \p B(a,1))}(1+o(1)) \; \text{by translation invariance} \nn \\
&= \alpha_\square^{\eta,a,\theta}(\eta,1)(1+o(1))\,, \nn
\end{align} 
where $o(1)$ goes to zero uniformly as $\eps$ goes to zero (independently of $a,\theta$). 
Only the second equality requires an argument: 

\begin{figure}[htbp]
\SetLabels
(.47*.47)$x$\\
(1.1*.24)$B(a,1)$\\
(.82*.3)$B(x,1)$\\
\endSetLabels
\centerline{
\AffixLabels{
\includegraphics[height=2 in]{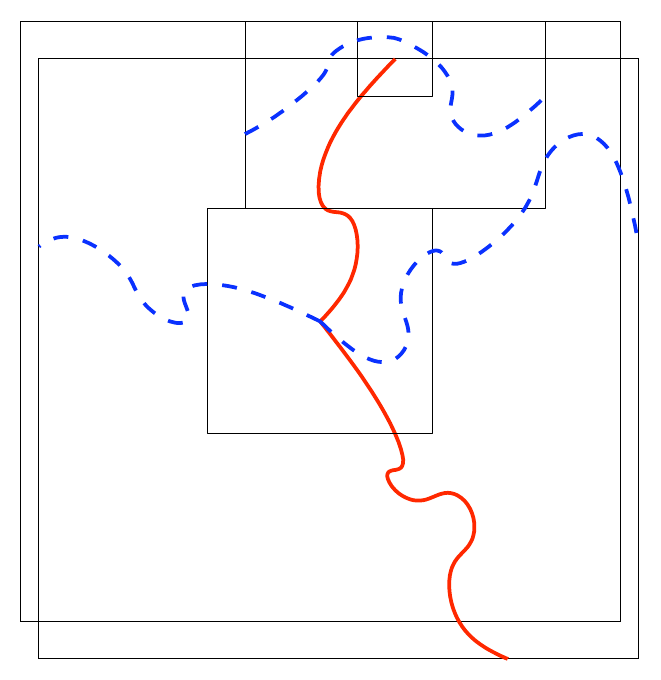}}}
\caption{The effect of shifting the box on the four-arm event.}
\label{f.shiftB}
\end{figure}

The square $B(x,1)$ is at least $\eps$-close to $B(a,1)$ when 
$x\in B(a,\eps)$, which will imply that 
\begin{equation}
\label{e.shiftB}
\Pb{\mathcal{A}_\square(x, \p B(x,1)) \Delta \mathcal{A}_\square(x, \p B(a,1))} =o_{\eps\to 0} \bigl(\Pb{\mathcal{A}_\square(x, \p B(x,1))} \bigr).
\end{equation}
Indeed, the symmetric difference $\mathcal{A}_\square(x, \p B(x,1)) \Delta \mathcal{A}_\square(x, \p B(a,1))$ holds only if there are four arms 
from $x$ to $B(x,\frac 1 2)$ and if there is some ball of radius $\eps$ on $\p B(a,1)$ which has three arms in $B(a,1)\cup B(x,1)$ up to distance $1/2-\eps$.
See Figure~\ref{f.shiftB}.
Using the fact that the three arms exponent in $\H$ is 2 (also known for $\Z^2$), and using independence on disjoint sets, one easily obtains~\eqref{e.shiftB}.

Summarizing, since there are $(4+o(1))\eps^2 \eta^{-2}$ points $x\in B(a,\eps)$, one ends up with
\begin{align*}
\frac{\beta}{ \eta^{-2} \alpha_4^\eta(\eta,1)} &=  \frac{\alpha_\square^{\eta,a,\theta}(\eta,1)}{ \alpha_4^\eta(\eta,1)} \frac{4}{\eps^{-2} \alpha_\square^{\eta,a,\theta}(2\eps,1)} (1+o(1))\\
& = \frac{4C}{\eps^{-2}\alpha_4(2\eps,1)} \frac{\alpha_4(2\eps,1)}{\alpha_\square^{\eta,a,\theta}(2\eps,1)}(1+o(1))\quad \textrm{by (\ref{e.ratio})}\\
& =  \frac{4}{\eps^{-2}\alpha_4(\eps,1)} \frac{\alpha_4(\eps,1)}{\alpha_4(2\eps,1)} (1+o(1))\quad \textrm{by (\ref{e.ratio}) again}\\
& =  \frac{4\cdot 2^{-5/4}}{\eps^{-2}\alpha_4(\eps,1)}(1+o(1)) \quad \textrm{by scale invariance and (\ref{e.limitlambda})}\,,
\end{align*}
which ends the proof of Lemma~\ref{l.betalimit} and Proposition~\ref{pr.muB}. \qed
\medskip

We now conclude the proof of Theorem~\ref{th.measurable}, then explain what needs to be modified to prove Theorem~\ref{th.4meas}.

\proofof{Theorem~\ref{th.measurable}}
Let us summarize what the current state of the proof is.
We first showed that $\{(\omega_\eta, \mu^A_\eta)\}_{\eta>0}$ is tight. Moreover we know that $\omega_\eta$
converges in law to the continuum percolation $\omega$. Using the tightness of $\{(\omega_\eta, \mu^A_\eta)\}_{\eta}$,
one can thus
find a subsequence $\{\eta_k\}_{k\geq 1}$ and some limiting random measure $\mu^\star$ so that 
\begin{equation}\label{e.weakconv}
(\omega_{\eta_k}, \mu_{\eta_k}^A) \overset{d}{\longrightarrow} (\omega, \mu^\star) 
\end{equation}

In order to prove Theorem~\ref{th.measurable}, we need to show that there is a unique possible choice for the random measure $\mu^\star$
and that this random measure is a measurable function of the first coordinate $\omega$.

This requires some topological considerations.
First of all, recall that all $\omega_{\eta_k}$ and $\omega$ live in the space $\HH = \HH_{\hat \C}\,$. As we have seen in Subsection~\ref{ss.topology}, 
$(\HH,\T)$ is a
compact, metrizable, separable space.
Furthermore, all the (random) measures $\mu_{\eta_k}^A$ and $\mu^\star$ live in the space $\mathfrak{M}(\closure{\Delta})$ of finite measures defined on the compact set $\closure{\Delta}$.
This space, equipped with the weak* topology (which is called ``weak convergence of measures'' by probabilists) is a complete, metrizable, separable space. 
Here is an example of a metric on $\mathfrak{M}(\closure{\Delta})$ which induces the desired weak* topology: let $\{\phi_j\}_{j\geq 1}$
be a countable basis of the separable space ($\mathcal{C}^0(\closure{\Delta}), \| \cdot \|_\infty)$ of continuous functions on $\closure{\Delta}$ (assume also that $\| \phi_j \|_\infty \le 1$ for all $j$).
Then
\begin{align}
d(\nu^1, \nu^2):= \sum_{j \geq 1} \frac {| \nu^1(\phi_j) - \nu^2(\phi_j) |}{ 2^j}
\end{align}
defines a metric on $\mathfrak{M}(\closure{\Delta})$ for the weak* topology.

Hence all $(\omega_{\eta_k}, \mu_{\eta_k}^A)$ and $(\omega, \mu^\star)$ live on the metric space $\HH\times \mathfrak{M}(\closure{\Delta})$,
which is complete and separable. In particular, one can apply the Skorokhod representation theorem: since we assumed the convergence in law of $(\omega_{\eta_k}, \mu^A_{\eta_k})$
towards $(\omega, \mu^\star)$, there is a coupling (on the same probability space)
of all the $(\omega_{\eta_k}, \mu_{\eta_k}^A)$ and $(\omega,\mu^*)$ such that 
\begin{enumerate}
\item $\omega_{\eta_k}  \longrightarrow \omega$ almost surely (for the quad-crossing topology $\T$) 
\item $\mu_{\eta_k}^A \longrightarrow \mu^*$ almost surely (for the topology induced by the metric $d$).
\end{enumerate}
 
Let us now characterize uniquely the limit measure $\mu^\star$ with the help of the previous subsections.
\medskip

\noindent{\bf Characterisation of $\mu^\star$}. Proposition \ref{pr.muB} says that for any ball $B\subset \Delta$ (or in fact $\closure{\Delta}$ would work equally well) the random measure $\mu_\eta^A(B)$ is well approximated by $Y=Y_\eta^{\eps,a,\theta}$ (times a deterministic factor) which counts 
how many of certain macroscopic 4-arm events hold. 
In Subsection~\ref{ss.limitarms}, we have seen (see Remark~\ref{r.joint4arm}) that as $\eta$ goes to zero, $Y_\eta^{\eps,a,\theta}(\omega_\eta)$ converges in law 
to a random variable $Y^{\eps,a,\theta}(\omega)$ which is measurable with respect to the scaling limit $\omega$ of critical percolation. Moreover, in the above joint coupling, where $\omega_{\eta_k} \to \omega$ almost surely, we have, for any $\eps>0$ and any $B\subset\Delta$,  the convergence $Y^{\eps,a,\theta}_\eta \to Y^{\eps,a,\theta}$ in probability. Of course, we also have the convergence of 4-arm probabilities, so we obtain that, in probability,
 \[
 \frac {Y_\eta^{\eps,a,\theta}(\omega_\eta)} {\eps^{-2} \alpha_4^\eta(\eps,1)} \longrightarrow \frac {Y^{\eps,a,\theta}(\omega)} {\eps^{-2} \alpha_4(\eps,1)} \,.
 \]
Now, we have the following bound on the third moments:
\begin{equation}\label{e.3rdMomY}
\E\left( \frac {Y_\eta^{\eps,a,\theta}(\omega_\eta)} {\eps^{-2} \alpha_4^\eta(\eps,1)}\right)^3 < C_{A,B} < \infty\,,
\end{equation}
uniformly in $\eta$ and $\eps$, with the constant $C_{A,B}$ depending only on the annulus $A$ and the ball $B$. This is proved exactly as~(\ref{e.3rdMom}) in Lemma~\ref{l.Mom}, just summing over boxes on the $\eps$-scale instead of points on the $\eta$-scale. This~(\ref{e.3rdMomY}) implies that the sequence 
$$
\left\{ \left( \frac{Y_\eta^{\eps,a,\theta}(\omega_\eta) }{ \eps^{-2} \alpha_4^\eta(\eps,1)) }\right)^2\right\}_{\eta>0}
$$ is uniformly integrable, and hence the convergence in probability also implies convergence in $L^2$:
$$
\lim_{\eta\to 0} \left\|\frac{Y^{\eps,a,\theta}_\eta}{\eps^{-2}\alpha_4^\eta(\eps,1)}-\frac{Y^{\eps,a,\theta}}{\eps^{-2}\alpha_4(\eps,1)}\right\|_2=0\,.
$$
 Combining this with the $L^2$ control provided by Proposition~\ref{pr.muB} easily leads to the fact that $\left\{ \omega \mapsto \frac {c \,Y^{\eps,a,\theta}(\omega)}{\eps^{-2}\alpha_4(\eps,1)}\right\}_{\eps>0}$ 
is a Cauchy sequence in $L^2$ as $\eps\to 0$. In particular, it has a limit:
\begin{equation}
\mu^A(B) := \lim_{\eps\to 0} \frac{c \, Y^{\eps,a,\theta}}{\eps^{-2} \alpha_4(\eps ,1)}
\qquad\hbox{in $L^2$}.
\label{e.muY}
\end{equation}
This means that for any ball $B\subset \closure{\Delta}$, there is an $L^2$ limit $\mu^A(B)(\omega)$ which is a good candidate for the desired limit measure of the ball $B$.
\vskip 0.3 cm

One would like at this stage to conclude that almost surely, $\mu^\star(B) = \mu^A(B)(\omega)$; this would characterize our limit measure $\mu^\star$ as a measurable function 
of $\omega$. Unfortunately, indicator functions of balls do not behave well with respect to the weak* topology, therefore we need to modify~(\ref{e.muY}) to the case of continuous functions on $\closure{\Delta}$. 
(See Remark~\ref{r.avoid} below regarding the issues that might appear.)
It is a standard technique to use Riemann sums for such an extension, and we get the following proposition:

\begin{proposition}\label{pr.Zeps}
For any $j\geq 1$, and any $0< \eta < \eps$, there exist functions $Z_j^{\eta,\eps}$ of the random discrete configurations $\omega_\eta$ such that
as $\eps \to 0$ and $\eta/\eps \to 0$, one has
\begin{equation}\label{e.Zmu}
\| \mu_\eta^A(\phi_j) - Z_j^{\eta,\eps} \|_2 \longrightarrow 0\,,
\end{equation}
where the set $\{\phi_j\}$ is our countable basis of continuous functions.  
Furthermore, these functions $Z_j^{\eta,\eps}$ are ``macroscopic'' in the sense that they converge in law as $\eta\to 0$ towards a measurable function of $\omega$: $Z_j^{\eps} = Z_j^\eps(\omega)$. Moreover, in the coupling with the pointwise convergence of $\omega_{\eta_k}$ towards $\omega$, we have that $Z_j^{\eta_k,\eps}= Z_j^{\eta_k,\eps}(\omega_{\eta_k})$ converges in probability towards $Z_j^\eps=Z_j^\eps(\omega)$.

Finally, the sequence of observables $\{Z_j^\eps\}_{\eps>0}$ is a Cauchy sequence in $L^2$, hence it has an $L^2$-limit that we may write as
\begin{equation}\label{e.phiZ}
\mu^A(\phi_j):= \mu^A(\phi_j)(\omega) = \lim_{\eps\to 0}  Z_j^\eps(\omega)
\; \hbox{in $L^2$}.
\end{equation}
\end{proposition}

Just like in~(\ref{e.3rdMomY}), we have a uniform control on the $L^3$ norms of $Z_j^{\eta_k,\eps}$ and $Z_j^\eps$, hence from the convergence in probability under the above coupling we obtain that, for any $j\geq 1$:
\begin{equation}
\| Z_j^{\eta_k,\eps}(\omega_{\eta_k}) - Z_j^\eps(\omega) \|_2 \longrightarrow 0\,. \nn
\end{equation}
Combining this with equations (\ref{e.Zmu}) and (\ref{e.phiZ}) from the proposition, one obtains using the triangle inequality:
\begin{equation}\label{e.phieta}
\| \mu_{\eta_k}^A(\phi_j) - \mu^A(\phi_j) \|_2 \longrightarrow 0\,.
\end{equation}

Recall now that we also assumed the pointwise weak convergence of the measures $\mu_{\eta_k}^A$ towards $\mu^\star$.
Since from (\ref{e.3rdMom}) we have a uniform control on $\Eb {\mu_{\eta_k}^A (\closure{\Delta})^3}$, this implies that for any $j\geq 1$
\begin{equation}\label{e.starZ}
\| \mu_{\eta_k}^A(\phi_j) - \mu^\star(\phi_j)\|_2 \longrightarrow 0\,.
\end{equation}

Combining equations (\ref{e.phieta}) and (\ref{e.starZ}), we obtain that, almost surely, for any $j\geq 1$, $\mu^\star(\phi_j) = \mu^A(\phi_j)(\omega)$.
This characterizes uniquely $\mu^\star$ and furthermore, $\mu^\star$ is indeed measurable with respect to $\omega$. \qed



\begin{remark}\label{r.avoid} Proposition~\ref{pr.Zeps} is actually weaker than the convergence of the measures for a countable basis of balls $B_j\subset \closure{\Delta},\, j\geq 1$. If, for example, $\mu^\star$ happened to have Dirac point masses, then the convergence of integrals of continuous functions might not imply the convergence of the measures of balls, hence one would need some additional regularity property for the possible weak limits $\mu^\star$. This also means that the weak*-convergence to $\mu^*$ does not alone imply the convergence for balls, i.e., the $L^2$-limit $\mu_A(B)(\omega)$ might not be equal to $\mu^*(B)$ a.s.  For any given ball $B$, one could easily prove equality (see Corollary~\ref{c.ballconv} below), but for all balls $B$ simultaneously the situation is more involved, so we chose not to address this question.
\end{remark}

\begin{remark}\label{r.as}
In~(\ref{e.muY}) and (\ref{e.phieta}), we have convergence only in $L^2$, not almost surely. We expect that the stronger versions also hold, but our strategy does not seem to yield this, or at least not without a lot of additional work. For comparison, one might consider  the convergence of the random Riemann-Stieltjes sums to the It\^{o} integral. This holds only in $L^2$ in general, not almost surely; however, it is not very hard to see that for reasonably regular functions (say, H\"older$(\eps)$ for some $\eps>0$) that do not depend on the Brownian motion we are integrating against, the convergence does hold a.s.
\end{remark}

\begin{corollary}\label{c.ballconv}
For a proper annulus $A\subset\C$ with inner face $\Delta$, let $\mu^A$ be the limit measure defined by (\ref{e.phiZ}) in the coupling of all the $\omega_\eta$ and $\omega$. Then, for any given open set $U\subset\Delta$ with boundary $\p U$ having finite length, we have $\mu^A_\eta(U)\to \mu^A(U)$ in $L^2$. In particular, Proposition~\ref{pr.Mom} follows immediately from Lemma~\ref{l.Mom}.
\end{corollary}

\proof For any $\delta>0$, let $N_\delta(\p U)\subseteq\Delta$ be the set of points at distance at most $\delta$ from $\p U$, and for $i=1,2$ let $\phi_{\delta,i}:\Delta\longrightarrow [0,1]$ be continuous functions with $\phi_{\delta,1} \leq \1_U \leq \phi_{\delta,2}$ and $\mathrm{supp} (\phi_{\delta,i}-\1_U)\subseteq N_\delta(\p U)$. We obviously have 
\begin{equation*}
\mu^A(\phi_{\delta,1}) \leq \mu^A(\1_U) \leq \mu^A(\phi_{\delta,2}) \qquad\textrm{a.s.} 
\end{equation*}
On the other hand, note that $|\mu^A_\eta(\phi_{\delta,i})-\mu^A_\eta(\1_U)|^2 \leq \mu^A_\eta(N_\delta(\p U))$, hence
\begin{equation*}
\big\|\mu^A_\eta(\phi_{\delta,i})-\mu^A_\eta(\1_U)\big\|_2 \to 0 \quad\textrm{ as }\delta\to 0\,, 
\end{equation*}
uniformly in $\eta<\eta_0(\delta)$, by~(\ref{e.1stMom}) of Lemma~\ref{l.Mom}. Finally, for any $\delta>0$ and $i=1,2$, by~(\ref{e.starZ}), we have 
\begin{equation*}
\big\|\mu^A_\eta(\phi_{\delta,i})-\mu^A(\phi_{\delta,i})\big\|_2 \to 0 \quad\textrm{ as }\eta\to 0\,.
\end{equation*}
Combining the three displayed formulas gives the desired result. \qed

\proofof{Theorem~\ref{th.4meas}} The proof is basically the same as for Theorem~\ref{th.measurable}, with the following  modifications. We define the variables $x_i,y_i,X,Y$ using not $A$-importance, but $\Quad$-pivotality. The definition of $\beta$ does not have to be changed. The measurability of $Y$ in the scaling limit of percolation now relies on Remark~\ref{r.quadarm}. For the application of the coupling arguments to the $L^2$-approximation, there is the following issue: for those $Q_i$ boxes in the $\eps$-grid that are close to $\p\Quad$, the enlarged $\gamma$-box $Q_i'$ does not lie in $\Quad$, hence there is no room to apply the coupling lemma. The solution is to notice that there are actually no $\Quad$-pivotals close to $\p\Quad$, so, for small enough $\eps$, we can ignore the $\eps$-boxes close to $\p\Quad$. What do we mean by ``no pivotals close to $\p\Quad$''? If an $\eps$-box $Q_i$ is such that the $\gamma$-box $Q_i'$ intersects $\p\Quad$, then $Q_i$ can contain pivotals only if $Q_i'$ has the 3-arm event to the remaining three boundary arcs of $\p\Quad$. (Or, if $Q_i'$ intersects two boundary arcs, then it needs to satisfy a 2-arm event.) But the boundary of $\p\Quad$ is assumed to be piecewise smooth, the half-plane 3-arm exponent is 2, and the number of possible essentially different such $\gamma$-boxes is $O(1/\gamma)$, hence, by Markov's inequality, the probability that there exists such a $\gamma$-box goes to 0 as $\gamma\to 0$, uniformly in $\eta$. Therefore, in the joint coupling of $\omega_\eta$ and $\omega$, there is a random positive distance between $\p\Quad$ and the closest pivotal, and the coupling lemma can be applied.
These estimates also show that the measures $\mu^\Quad_\eta$ are still tight, despite the presence of the boundary. Hence all the arguments work fine.\qed

\subsection{Filtering}\label{ss.filtering}

It was convenient to introduce the measures $\mu_{\eta}^A$ for some proper annulus $A\subset \C$ instead of working directly with $\bar\mu_{\eta}^{\rho}$, the counting measure on $\rho$-important points.
Indeed, suppose that in the previous proof we were working with $\rho$-important points instead of $A$-important points. Then, for different points $x,y$ inside some $\eps$-square $Q_i$, given some configuration of faces around $Q_i$, we might need quite different information about where these faces are connected outside $Q_i$ if we want to know how $x$ and how $y$ has to be connected to these faces from the inside in order to be $\rho$-important. Say, if one of the four arms emanating from the faces around $Q_i$ goes to distance $\rho-\eps$ but not to $\rho$, then it might happen that some points in $Q_i$ that are connected ``pivotally'' to the four arms will be $\rho$-important while others will not.
$A$-important points are simpler to handle in this respect.

If one wants to use the more natural concept of $\rho$-important points, then needs to relate in some way the measures $\mu^A_\eta$ to the measures $\bar\mu^{\rho}_\eta$. The following partial order between enhanced tilings will handle this issue.

\begin{definition}\label{d.refine}
Let us say that an enhanced tiling $\mathcal{H}$ {\bf refines} another, $\mathcal{H}'$, denoted by $\mathcal{H} \leq \mathcal{H'}$, if the following holds: for any pair of annuli $A=B_2\setminus B_1 \in\mathcal{H}$ and $A'=B_2'\setminus B_1'\in\mathcal{H}'$, if the inner faces $B_1$ and $B_1'$ intersect each other, then $B_2\subset B_2'$.
\end{definition}

For example, one can consider the enhanced tiling $$\mathcal{H^\rho_\eta}:=\big\{B(x,\rho)\setminus \{x\}: \text{all $\eta$-tiles $x$ of }D\big\}$$ in a domain $D$. Now, if $\mathcal{H}$ is an enhanced tiling with $\diam\,\mathcal{H} < \rho$, then $\mathcal{H} \leq \mathcal{H^\rho_\eta}$. On the other hand, if $d(\p_1 A,\p_2 A)> \rho$ for all $A\in\mathcal{H}$, then $\mathcal{H} \geq \mathcal{H}^\rho_\eta$.

The point of this definition is that if $\mathcal{H} \leq \mathcal{H'}$, then we have the reversed domination between the associated annulus-pivotal measures: $\mu^\mathcal{H}_\eta \geq \mu^{\mathcal{H}'}_\eta$. Therefore, there is a coupled pair $(\mathcal{P},\mathcal{P}')$ of Poisson samples from these measures such that $\mathcal{P}\supset \mathcal{P}'$. In particular, if $\big\{\mathcal{H}(\eps) : \eps\in I \big\}$ is an ordered family of enhanced tilings, i.e., $\mathcal{H}(\delta) \le \mathcal{H}(\eps)$ whenever $\delta < \eps$, for $\delta,\,\eps\in I \subset \R_+$, then for each $\eta>0$ we get a family $\big\{\mu^{\mathcal{H}(\eps)}_\eta : \eps\in I\big\}$ of increasing measures (as $\eps$ decreases), called a {\bf filtered measure}, and there is an associated increasing family of Poisson samples $\big\{\mathcal{P}(\mathcal{H}(\eps)) : \eps\in I\big\}$.

If $\big\{\mathcal{H}(\eps) : \eps\in I\big\}$ is an ordered family of enhanced tilings satisfying $\diam \,\mathcal{H}(\eps) < \eps$ and $d(\p_1 A,\p_2 A) > \eps/2$ for all $A\in\mathcal{H}(\eps)$, then we have the following dominations between the associated filtered measure and the $\rho$-pivotal measures. For any $\rho>0$, if $\eta>0$ is small enough, then
\begin{equation}
\mu^{\mathcal{H}(\rho)}_\eta
\geq \bar\mu^\rho_\eta = \mu^{\mathcal{H}^\rho_\eta}_\eta
\geq \mu^{\mathcal{H}(2\rho)}_\eta\,.
\label{e.domi}
\end{equation}
This is a useful tool for comparing annulus-pivotality to $\rho$-pivotality.

%
%

\section{Limit results for other special points}\label{s.others}

We will now explain how to extend the above techniques to the case of counting measures on other geometrical objects of interest in critical percolation, namely:
\begin{itemize}
\item An ``area''-measure on the clusters. 
\item A ``length''-measure for interfaces ($\mathrm{SLE}_6$ curves in the scaling limit).
\item A ``length''-measure on exterior boundaries ($\mathrm{SLE}_{8/3}$ curves).
\end{itemize}
These three families of measures will be called {\it cluster measures}, {\it interface measures}
and {\it frontier measures}, respectively.
\medskip

One difference from the case of the {\it pivotal measure} is that the $\sigma$-algebra in which the limit measures will be measurable will not always be the Borel $\sigma$-algebra of the quad-crossing space, but, for instance, the {\it interface measure} will be defined as a function of the $\mathrm{SLE}_6$ curve.
\medskip

The proofs are in large part very similar to the proof given for the {\it pivotal measure}, hence we will only highlight the places that need modifications. The first case is much easier to establish due to the monotonicity of the one-arm event. The second one is very similar to the construction of the {\it pivotal measure} (except that four-arm events are replaced by two-arm events), while we will see that the third one requires some more care. In general, the analog of Section~\ref{s.measurable} carries through easily in each case, while the coupling part (Section~\ref{s.coupling}) needs some fixing that depends on the number of arms.
\medskip

Finally, after the construction of the {\it cluster measure}, we will prove a corollary of independent interest on the so-called
{\bf two-point function} in percolation, namely that the two-point connectivity function is asymptotically rotationally invariant.

\subsection{The {\it cluster measure} (or {\it area measure})}\label{ss.clusters}

\begin{minipage}{\textwidth}
\begin{center}
\includegraphics[width=0.5 \textwidth]{specialWWcluster}
\end{center}
\end{minipage}
\vskip 0.5 cm

In this subsection, we wish to construct an ``area''-measure on the clusters of critical percolation. 
Similarly to the rest of the paper, a natural way to proceed is to consider the properly rescaled counting measures on 
``macroscopic'' clusters. Let $\rho>0$ be some cut-off length and consider the following rescaled measure on the clusters 
of an $\omega_\eta$ percolation configuration in $\eta\Tg$ whose diameters are greater than $\rho$:
\begin{equation}
\aream_\eta^\rho = \aream_\eta^\rho(\omega_\eta) := \sum_{x:\; \diam(\mathcal{C}_x)> \rho} \delta_x \, \frac {1} {\eta^{-2}\alpha_1^\eta(\eta,1)}\,,
\end{equation}
where for each $x\in\eta \Tg$, $\mathcal{C}_x$ denotes the cluster attached to $x$.

For the same reasons that we encountered in the case of pivotals, it turns out to be more convenient to work with the following related measures: for any proper annulus $A$, define
\begin{equation}
\aream_\eta^A = \aream_\eta^A(\omega_\eta) := \sum_{x:\; \{x \leftrightarrow \p_2 A\} } \delta_x \, \frac {1} {\eta^{-2}\alpha_1^\eta(\eta,1)}\,.
\end{equation}

With techniques similar to the ones used for the {\it pivotal measure}, one can prove the following analog of Theorem~\ref{th.measurable}. We will only explain below what adaptations are needed.

\begin{theorem}\label{th.areameasure}
Let $A$ be a fixed proper annulus of the plane. When $\eta\to 0$, the random variable $(\omega_{\eta}, \aream_{\eta}^A)$ converges in law
to some $(\omega, \aream^A)$, where $\omega$ is the scaling limit of critical percolation, and the measure $\aream^A= \aream^A(\omega)$ is a measurable function of $\omega$. The measure thus defined may be thought of as an {\bf area measure} on macroscopic clusters of percolation.
\end{theorem}

A filtering argument as in Subsection~\ref{ss.filtering} allows us, if needed, to deal with the initial $\{\aream^\rho\}_{\rho>0}$ measures.
\medskip

{\bf Rough sketch of proof.}
Due to the monotonicity of the one-arm event, the proof is in many ways much simpler than the one for the {\it pivotal measure}.
The reason for this is that all the separation of arms type of results needed in Sections~\ref{s.coupling} and~\ref{s.measurable} are not needed here. RSW is the only necessary tool in this case.
\vskip 0.2 cm

Let us first sketch what the analog of the coupling part (Section~\ref{s.coupling}) is.
We start with some notations. If $\eta<r <R$, let $\mathcal{A}_1(r,R)$ be the one-arm event 
from radius $r$ to radius $R$. If $\Omega$ is any smooth simply connected domain and $\eta <r $ is such that $B(0,r)\cap \p \Omega =\emptyset$,
define $\mathcal{A}_1(r,\p \Omega)$ to be the event that there is an open path from $\p B(0,r)$ to $\p \Omega$.
Note that, in this definition, $B(0,r)$ need not be included in $\Omega$. Also if $\Omega= B(0,R)$, this event coincides with $\mathcal{A}_1(r,R)$, so we will only state the coupling result in this general setup.

Now if $\eta<r < u$, let $\Gamma = \Gamma(r,u)$ be the outermost open circuit in the annulus $A(r,u)$, i.e., the simple open circuit which is the closest 
to $\p B(0,u)$. If there are no such circuits, then let $\Gamma:= \emptyset$. It is a standard fact that such a circuit can be revealed in a ``Markovian'' way,
only using bits from the exterior of $\Gamma$ in such a way that  the region surrounded by $\Gamma$ is left untouched.

The analog of the coupling result that we need is the following:

\begin{proposition}\label{pr.couplingonearm}
Let $\Omega$ be some piecewise smooth simply connected domain with $0\in \Omega\setminus \p \Omega$. Let $d>0$ be the distance from 0 to the boundary $\p \Omega$ and $d' = d\wedge 1$. 

For any $10 \eta < r < \frac {d'} {100}$, let $\nu$ be the law on percolation configurations conditioned on the event $\mathcal{A}_1(r,\p \Omega)$,
and let $\nu_0$ be the law conditioned on the event $\mathcal{A}_1(r,1)$.

Let $h:= \lfloor \log_2(\frac {d'} r) \rfloor$ and $u:= 2^{h/2} r$. Let $\Gamma=\Gamma(r,u)$ and $\Gamma^0=\Gamma^0(r,u)$ be the outermost contours under $\nu$ and $\nu^0$.

Then, there is a coupling of the conditional laws $\nu$ and $\nu^0$ such that with (conditional) probability at least $1-(r/d')^k$, both $\Gamma$
and $\Gamma^0$ are non-empty and are identical, where $k>0$ is some absolute exponent.
\end{proposition}

The proof of this coupling is essentially in Kesten's work on the IIC \cite{\KestenIIC}, but let us briefly explain how to adapt the strategy of Section~\ref{s.coupling} to this case. Namely, one tries to couple the two measures in each dyadic annulus ranging from $d'$ to $u= 2^{h/2} r$, where ``dyadic'' will actually correspond to $\log_8$, for convenience. Roughly speaking, in each of these annuli, say $A(s, 8s)$, with positive probability there are circuits in each of $A(2s,4s)$, $A(4s,6s)$ and $A(6s,8s)$. Let $\Gamma^1$, $\Gamma^2$, $\Gamma^3$ be respectively the outermost circuit in $A(2s,4s)$, the innermost one in $A(4s,6s)$,
 and the outermost one in $A(6s,8s)$. Following the notations of Section~\ref{s.coupling}, let $S=S(\omega_\eta)$ be the union of all hexagons which are determined 
 by $\Gamma^1$ and $\Gamma^2$ plus all hexagons which lie between $\Gamma^1$ and $\Gamma^2$. Let also $m_S$ be the number of hexagons 
 that are in $S$.
 
Whenever there is an open circuit in $A(6s,8s)$, which happens with positive conditional probability under $\A_1(r,1)$ or $\A_1(r,\p\Omega)$, we can explore $A(6s,8s)$ from the outside in a way that detects the outermost circuit $\Gamma^3$ and does not touch anything inside $\Gamma^3$. Then, using RSW, one easily obtains that there is some constant $0<c<\infty$ such that, for any admissible $S$,
 \[
 2^{-m_S} / c \le \Pb{S(\omega_\eta)=S \bigm| \exists\,\Gamma_3,\ \mathcal{A}_1(r, \Gamma^3)} \le c 2^{-m_S}\,.
 \]

Now, following the same technology as in Section \ref{s.coupling}, one has that with high probability ($1- (u/d')^{k_1}$ for some $k_1>0$),
the coupling succeeds before reaching the intermediate radius $u=2^{h/2} r$, i.e., we have identical circuits $\Gamma^1$ and $\Gamma^2$ in the two configurations, given by $\nu$ and $\nu^0$, respectively. Since the exploration was Markovian so far, one can now couple these two measures so that they match inside $\Gamma^1$. Then, with high (conditional) probability $1-(r/u)^{k_2}$, there is an (outermost) open circuit $\Gamma= \Gamma(r,u)$ in $A(r,u)$. By construction, this contour is identical for the two measures and can be discovered in a Markovian 
way from the outside (so, this contour could be thought of as a ``stopping time''). \qed


Now, as with the pivotal points, it remains to use this coupling property in order to prove Theorem~\ref{th.areameasure}. It is straightforward to adapt Section~\ref{s.measurable} to our present case. Again, the present situation is much simpler in many ways.  For example, if one defines $X$ and $Y$ as before, then, in the second moment estimate, two $\eps$-squares $Q_i$ and $Q_j$ in general interact very little (while, with the pivotal points, this interaction was a non-negligible issue).
The reason for this is that if one conditions $Q_i$ and $Q_j$ to be part of a large cluster, then with high (conditional) probability, there are open circuits around $Q_i$ and $Q_j$ disconnecting one from the other. This implies that with high conditional probability the behavior in one box is decorrelated from the behavior in the other one.

So, we get the $L^2$-approximation $\Eb{(X-\beta Y)^2}=o((\E X)^2)$, using 
$$
\beta=\beta(\eta,\eps,a,\theta):=\Pb{ x_0 \md \A_1^{a,\theta}(2\eps,1)},
$$ 
where $ \A_1^{a,\theta}(2\eps,1)$ denotes the one-arm event from the rotated $2\eps$-square to the rotated $1$-square, and $x_0$ is the number of $\eta$-tiles in the rotated $\eps$-square $Q_0=B(a,\eps)$ that are connected to $\p B(a,1)$. Then we need to prove the analogs of the ratio limit results of Subsection~\ref{ss.ratio}: for any fixed $r>0$,
\begin{equation*}
\lim_{\eta\rightarrow 0}\frac{\alpha_1^\eta(\eta,r)}{ \alpha_1^\eta(\eta, 1)} = \lim_{\eps\to 0} \frac{\alpha_1(\eps,r)}{\alpha_1(\eps,1)} = r^{-5/48}\,,
\end{equation*}
and there is an absolute constant $C>0$, such that, uniformly in the orientation $a,e^{i\theta}$,
\begin{equation*}
\lim_{\eta\to 0} \frac{\alpha_1^{\eta,a,\theta}(\eta,1)}{\alpha_1^\eta(\eta,1)} =
\lim_{\eps\to 0} \frac{\alpha_1^{a,\theta}(\eps,1)}{\alpha_1(\eps,1)} = 
C\,.
\end{equation*}
Using the above coupling property, these again require just straightforward modifications.
   
 Finally, to prove that the limiting measure is indeed measurable with respect to the continuum percolation, one needs to argue that the macroscopic observable $Y=Y^{\eps,\eta}$ (which in this case would count roughly the number of $\eps$-squares which are connected to $\p_2 A$) has a limit when the mesh $\eta\to 0$ and that this limit is measurable with respect to $\omega$ (i.e., $Y^{\eps,\eta} \to Y^\eps=Y^\eps(\omega)$). This is precisely the purpose of Lemma~\ref{l.meas4arm}. \QED

\subsection{Rotational invariance of the two-point function}\label{ss.twopoint}

The result of this subsection uses rotational invariance of the scaling limit plus the coupling technology for one arm, Proposition~\ref{pr.couplingonearm}. Since, as we mentioned in the previous subsection, the coupling idea was already a key ingredient in \cite{\KestenIIC}, this result is not fundamentally new; however, we do not know of an explicit appearance of it so far, so we decided to include it here to highlight the power of the coupling ideas. 

The asymptotic rotational invariance of the two-point function in critical percolation on $\Z^d$ for $d\geq 19$ was proved by very different methods by Hara \cite{\Hara}.

\begin{proposition}\label{pr.twopoint} There exists a universal constant $C>0$ such that if $x$ and $y$ are two sites of $\eta\Tg$, with Euclidean distance $d=d(x,y)>0$ (hence graph distance $\sim g(\theta)\,d/\eta$, with the constant $g$ depending on the angle $\theta$ of the vector $y-x$), then 
$$\Pb{x\mathop{\llra}\limits^{\omega_\eta} y} = (C+o(1))\, \alpha_1^\eta(\eta,d/2)^2$$ as $\eta\to 0$, uniformly in the angle $\theta$.
\end{proposition}

\proof
The proof will be very similar to the one we used for the {ratio limit theorem} in Subsection \ref{ss.ratio}.
Namely, we will show that $\Pb{x \overset{\omega_\eta}{\longleftrightarrow} y}/ \alpha_1^\eta(\eta,d/2)^2$ is a Cauchy sequence 
whose limit is a universal constant $C>0$ independent of $x$ and $y$.

Let us start by 
\begin{equation}\label{e.xy}
\Pb{x \overset{\omega_\eta}{\longleftrightarrow} y}  =  \Pb{x \leftrightarrow y \bigm| B(x,\eps)  \leftrightarrow  B(y, \eps) } \, \Pb{ B(x,\eps)  \overset{\omega_\eta}{\longleftrightarrow} B(y, \eps)} \,.
\end{equation}

Now observe that 
\begin{align}
\Pb{ x \leftrightarrow y \bigm| B(x,\eps)  \leftrightarrow  B(y, \eps)} & \label{e.xyeps} \\
& \hskip -4 cm = \Pb{ x \lra \p B(x,d/2)\,, \, y \lra \p B(y,d/2) \bigm|  B(x,\eps)  \leftrightarrow  B(y, \eps) } (1+o(1))\,, \nonumber
\end{align}
where $o(1)$ goes to zero polynomially in $\eps/d$. This follows from the fact that 
the symmetric difference of the event 
$B_1:= \mathcal{A}_1(x,d/2)\cap \mathcal{A}_1(y,d/2)\cap\{  B(x,\eps)  \leftrightarrow  B(y, \eps) \}$ with  $B_2:=\{ x \llra y \}$ 
is included in the event
\[
\mathcal{A}_1(x,\eps)\cap \mathcal{A}_1(y,\eps)\cap \Bigl( \mathcal{A}_4(x, \eps, d/2) \cup \mathcal{A}_4(y,\eps,d/2) \Bigr) \,,
\]
whose probability is much smaller than the probabilities of $B_1$ or $B_2$ (this can be seen using quasimultiplicativity and the fact that 
$\alpha_4(r,R) \ll \alpha_1(r,R)^2$).

Therefore, one can now work instead with the event 
$$
\Pb{ \A_1(x,d/2),\; \A_1(y,d/2) \bigm|  B(x,\eps)  \leftrightarrow  B(y, \eps)} \,.
$$

Now, if $\nu_{x,y}$ denotes the conditional law $\Pb{\, \bullet \bigm|  B(x,\eps) \lra  B(y,\eps) }$ and if 
$\nu_x,\, \nu_y$ denote the conditional laws $\Pb{\, \bullet \bigm|   B(x,\eps) \lra  B(x, d/2)}$ and $\Pb{\, \bullet \bigm|   B(y,\eps) \lra  B(y, d/2)}$ respectively, then, similarly to Subsection~\ref{ss.clusters}, one can couple $\nu_{x,y}$ with $\nu_x \otimes \nu_y$ (product measure, or in other words independent copies).
What we mean by coupling here is the straightforward analog of Theorem~\ref{th.areameasure}, except here one couples ``simultaneously'' around $B(x,\eps)$ and $B(y,\eps)$. If $\SS$ is the event that the coupling succeeds (which has probability at least $1-(\eps/d)^k$ for some $k>0$), then, in $\nu_{x,y}$,
\begin{align}
\Pb{ \A_1(x,d/2),\; \A_1(y,d/2) \md  B(x,\eps)  \leftrightarrow  B(y, \eps)} & \label{e.joint}\\
& \hskip -6.5 cm = \Pb{  \A_1(x,d/2),\; \A_1(y,d/2) \md \SS,\; B(x,\eps)  \leftrightarrow  B(y, \eps)} \, \Pb{\SS \md B(x,\eps)  \leftrightarrow  B(y, \eps)} \nonumber \\
& \hskip -6.5 cm + \Pb{ \A_1(x,d/2),\; \A_1(y,d/2) \md \neg\SS,\;  B(x,\eps)  \leftrightarrow  B(y, \eps)} \, \Pb{\neg\SS \md B(x,\eps)  \leftrightarrow  B(y, \eps)}\,.  \nonumber
\end{align}
In the second term (with $\neg\SS$), the first factor is at most $\Pb{\A_1(x,\eps),\; \A_1(y,\eps)} = \alpha_1(\eta,\eps)^2$, while the second factor is $o(1)$ as $\eps/d \to 0$. In the first term (with $\SS$), the first factor, by the definition of $\SS$, equals 
\begin{equation}\label{e.bothS}
\Pb{  \A_1(x,d/2),\; \A_1(y,d/2) \md \SS,\; B(x,\eps) \lra  B(x, d/2),\; B(y,\eps) \lra  B(y, d/2) }\,,
\end{equation}
while the second factor is $1-o(1)$. Now, in  $\nu_x \otimes \nu_y$, on one hand we directly have 
\begin{align}
\nu_x \otimes \nu_y[\A_1(x,d/2),\; \A_1(y,d/2)] &= \Pb{ x \lra B(x,d/2) \bigm| B(x,\eps) \lra B(x,d/2)}^2\nonumber\\
&= \frac{\alpha_1^\eta(\eta,d/2)^2}{\alpha_1^\eta(\eps,d/2)^2} \asymp \alpha_1^\eta(\eta,\eps)^2\,;\label{e.separ}
\end{align}
on the other hand, we can do the decomposition according to $\SS$ and $\neg\SS$, where the latter term is again $o(1)\, \alpha_1^\eta(\eta,\eps)^2$. This is negligible compared to (\ref{e.separ}), making the first term dominant:
\begin{align*}
\nu_x \otimes \nu_y \big[ \A_1(x,d/2),\; \A_1(y,d/2) \big] & \\
&\hskip - 5.5 cm = (1-o(1)) \, \Pb{  \A_1(x,d/2),\; \A_1(y,d/2) \md \SS,\; B(x,\eps) \lra  B(x, d/2),\; B(y,\eps) \lra  B(y, d/2) }\,.
\end{align*}
This is asymptotically the same as (\ref{e.bothS}), hence also in~(\ref{e.joint}) the $\neg\SS$-term is negligible compared to the $\SS$-term,  and altogether we get
\begin{align*}
\nu_{x,y}\big[ \A_1(x,d/2) ,\; \A_1(y,d/2) \big] &= (1+o(1)) \, \nu_x \otimes \nu_y \big[ \A_1(x,d/2) ,\; \A_1(y,d/2) \big] \\
& = (1+o(1)) \, \frac{\alpha_1^\eta(\eta,d/2)^2}{\alpha_1^\eta(\eps,d/2)^2}\,.
\end{align*}

Combining this with (\ref{e.xyeps}), we can write~(\ref{e.xy}) as 
\[
\Pb{x \overset{\omega_\eta}{\longleftrightarrow}  y} = (1+o(1)) \, \frac{\alpha_1^\eta(\eta,d/2)^2}{\alpha_1^\eta(\eps,d/2)^2} \, \Pb{B(x,\eps)  \overset{\omega_\eta}{\longleftrightarrow} B(y,\eps)} \,,
\]
where $o(1)$ goes to zero as $\eps$ goes to zero, uniformly in $\eta < \eps/10$.

Now, from conformal invariance on the triangular lattice and $\mathrm{SLE}$ technology, we get that for fixed $\eps>0$ and as the mesh $\eta$ goes to zero,
\[
\lim_{\eta\to 0}  \frac{\Pb{B(x,\eps)  \overset{\omega_\eta}{\longleftrightarrow} B(y,\eps)}}{\alpha_1^\eta(\eps, \frac d 2)^2} = 
\frac{\Pb{B(x,\eps) \longleftrightarrow B(y,\eps)}}{\alpha_1(\eps, \frac d 2)^2}  =: C(\eps/d)
\]
exists and depends only on the ratio $\eps/d$ (by scaling invariance and rotational invariance of critical percolation on $\Tg$).

Altogether, this gives us
\[
\frac{\Pb{x \overset{\omega_\eta}{\longleftrightarrow}  y} }{\alpha_1^\eta(\eta, \frac d 2)^2} = C(\eps/d)(1+o(1))\,,
\]
where $o(1)$ goes to zero as $\eta$ and $\eps$ go to zero with $\eta < \eps/10$. This proves the needed Cauchy criterion as $\eta\to 0$. Moreover, the resulting limit $C=C^{x,y}$ must equal $\lim_{\eps\to 0} C(\eps/d)$ (which thus needs to exist), implying that $C^{x,y}$ is in fact independent of the points $x$ and $y$. The finiteness and positivity of this limit $C$ is clear from quasi-multiplicativity arguments.  \qed

\begin{remark}
Using Proposition \ref{pr.ratio}, one may rewrite this result as 
\[
\Pb{x\mathop{\llra}\limits^{\omega_\eta} y} = \alpha_1^\eta(\eta,|x-y|)^2 (\tilde C+o(1)) \,,
\]
for some universal $\tilde C>0$.
\end{remark}

\subsection[The {\it interface measure} on $\mathrm{SLE}_6$ curves]{The {\it interface measure} (related to the natural parametrization of $\mathrm{SLE}_6$ curves)}\label{ss.interface}

\begin{minipage}{\textwidth}
\begin{center}
\includegraphics[width=0.5 \textwidth]{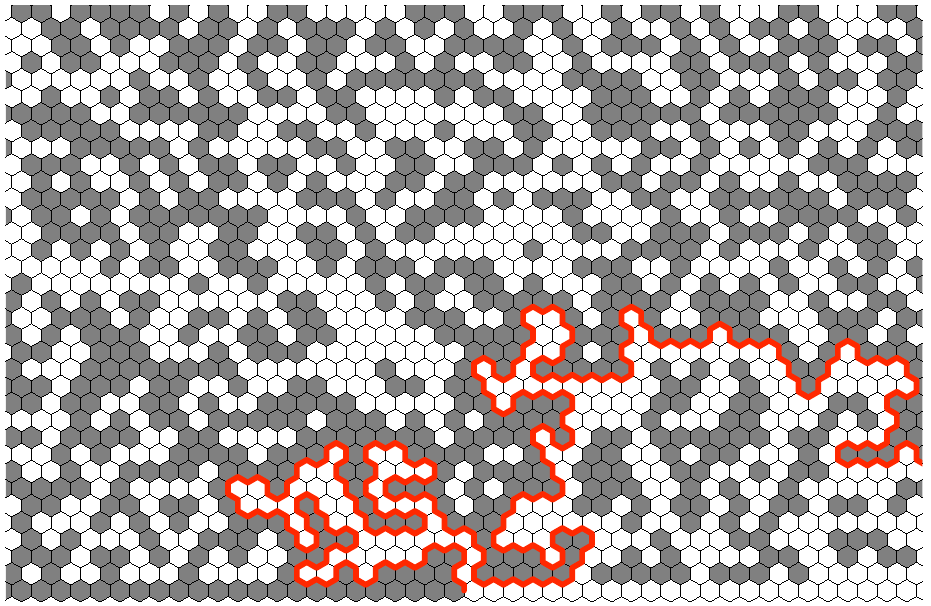}
\end{center}
\end{minipage}
\vskip 0.5 cm

Interfaces between open and closed clusters have been much studied since the discovery of the $\mathrm{SLE}$ processes.
When the mesh $\eta$ goes to zero, these curves are known to converge towards $\mathrm{SLE}_6$ curves. At the scaling limit, these curves are 
classically parametrized by their ``capacity''.  This capacity parametrization also makes sense on the discrete level for the exploration processes $\gamma_\eta$ at mesh $\eta$, 
but it is very different from the natural parametrization on $\gamma_\eta$ which consists in counting the number of iterated steps done by the exploration process so far. Whence the question of defining a parametrization of the $\mathrm{SLE}_6$ curve at the scaling limit, which would be, in an appropriate sense, the scaling limit of the 
natural discrete parametrizations.
\medskip

In this subsection, we prove a result in this direction. Let us define some notation before stating the main result.
Let $(\Omega,a,b)$ be a smooth domain with two marked points on its boundary. For any mesh $\eta>0$, let $(\Omega_\eta, a_\eta, b_\eta)$ be the discrete domain 
approximating $(\Omega,a,b)$ (in any ``reasonable'' way, say from inside for example). Let $\gamma_\eta$ denote the chordal exploration process (with Dobrushin's boundary conditions) from $a_\eta$ to $b_\eta$. By definition, $\gamma_\eta$ consists of all edges $e$ of the hexagonal lattice which lie between an open hexagon connected to the open arc and a closed hexagon connected to the closed arc. Each edge $e$ corresponds to one iteration of the exploration path, hence each edge contributes one unit of time to the natural parametrization. We thus define the following counting measure supported on $\gamma_\eta$:
\begin{equation}\label{e.param}
\param_\eta = \param_\eta(\omega_\eta) = \param_\eta(\gamma_\eta):= \sum_{\mathrm{edges} \, e \in \gamma_\eta}  \delta_e \, \frac {\eta^2} {\alpha_2^\eta(\eta,1)}\,.
\end{equation}

As in the case of the {\it pivotal measure} or the {\it cluster measure}, we would like to show that this counting measure $\param_\eta$, which we will call
the {\it interface measure}, has a scaling limit when $\eta\to 0$. More precisely, we would like the joint $(\omega_\eta, \param_\eta)$ to converge towards 
a coupling $(\omega,\param)$, where $\param=\param(\omega)$ would be measurable with respect to the scaling limit.

One could prove such a statement (following the setup of Sections~\ref{s.coupling} and~\ref{s.measurable}), but this would have the disadvantage that the exact relation of the resulting $\param=\param(\omega)$ to the $\mathrm{SLE}_6$ curve $\gamma$ (the scaling limit of the exploration paths $\gamma_\eta$) would be unclear. 

To start with, even the measurability of $\gamma$ w.r.t.~$\omega$ is unresolved, see Question~\ref{q.SLEmeas} at the end of Subsection~\ref{ss.limitarms}. But, assuming an affirmative answer to that question, $\gamma$ contains much less information than $\omega$. In definition~(\ref{e.param}) of the discrete interface measure, we wrote $\param_\eta = \param_\eta(\omega_\eta) = \param_\eta(\gamma_\eta)$: the amount of information in $\omega_\eta$ and $\gamma_\eta$ relevant to $\param_\eta$ is obviously the same. But it is not clear what happens in the scaling limit. Hence, at least assuming Question~\ref{q.SLEmeas}, having a scaling limit measure $\param=\param(\gamma)$ measurable w.r.t.~the SLE curve $\gamma$ would be stronger and more satisfying than measurability w.r.t.~the entire quad-crossing configuration $\omega$. Therefore, this is the path we are taking.
\medskip

One more issue to mention is whether the interface measure $\param$ can be considered as a time parametrization  of the $\mathrm{SLE}_6$ curve $\gamma$. We expect so, but as explained in Subsection~\ref{ss.time}, the measure does not see the order of points in which the curve is traversed, hence additional work is needed. Even if we assume that the support of $\param$ is the entire range of $\gamma$, which will be proved in \cite{\MetricProp}, we still have the following questions. For any subdomain $F$ that contains a neighbourhood of the starting point $a$, let $\gamma_F\subset \gamma$ be the curve stopped upon first exiting $F$, and let $\param_F$ be the measure collected so far. Then we clearly have $\param_F(\gamma_F)\leq \param(\gamma_F)$, but do we have equality? If not, then we clearly do not have a time parametrization by $\param$. Moreover, we need almost sure equality not only for a fixed $F$, but for all $F$ simultaneously. We will try to address these questions in \cite{\MetricProp}. 
\medskip

Now, let us briefly recall what the setup of convergence of $\gamma_\eta$ towards $\mathrm{SLE}_6$ is. It is a convergence in law under the Hausdorff topology on simple curves.
See \cite{\WWperc} for more detail on this.
For simplicity, let us restrict ourselves to the case where the domain $\Omega$ (with two marked points $a,b \in \p \Omega$) is smooth and bounded.
We consider our exploration processes $\gamma_\eta$ as continuous paths $\gamma_\eta: [0,1] \to \closure{\Omega}$, such that $\gamma_\eta(0)=a$ and $\gamma_\eta(1)=b$.
The space $\CC$ of such curves $\gamma$ (which lie in $\closure{\Omega}$) is equipped with the (pseudo)-distance
\[
d(\gamma,\gamma') = \inf_{\phi} \sup_{u\in[0,1]} |\gamma(u) - \gamma'(\phi(u)) |\,,
\] 
where the infimum is taken over all continuous increasing bijections of $[0,1]$ onto itself. Let $\TC$ be the topology on $\CC$ generated by this distance.
The important theorem by Smirnov \cite{\SmirnovPerc,\SmirnovICM}
is that the discrete exploration paths $\gamma_\eta$ converge in 
law, in the sense of the Borel $\sigma$-field of $(\CC,\TC)$, towards an $\mathrm{SLE}_6$ curve in $\closure{\Omega}$ from $a$ to $b$ (see also \cite{\CamiaNewmanConv}).

We wish to prove the following theorem (an analog of Theorem \ref{th.4meas}):

\begin{theorem}\label{th.lengthmeasure}
Let $\Omega$ be a smooth bounded domain of the plane with two marked points $a,b \in \p \Omega$.
When $\eta\to 0$, the random variable $(\gamma_{\eta}, \param_{\eta})$ converges in law,
under the topology of $(\CC,\TC)$,
to $(\gamma, \param)$, where $\gamma$ is a chordal $\mathrm{SLE}_6$ in $\closure{\Omega}$ from $a$ to $b$, and the measure $\param= \param(\gamma)$ is a {\bf measurable function} of $\gamma$. We call this random measure the {\bf length measure} of $\gamma$.
\end{theorem}

\begin{remark}
The same theorem holds for non-bounded smooth domains ($\neq \C$), but for simplicity we stick to this case (essentially this would boil
down to looking at larger and larger bounded regions in $\closure{\Omega}$).
\end{remark}

\proof As in Section~\ref{s.measurable}, we first discuss tightness and boundary issues, then we characterize uniquely the subsequential scaling limits.
\medskip

\ni {\bf Tightness and treating the boundary issues.} As in Subsection~\ref{ss.tightness}, in order to prove the tightness of the measures $\param_\eta$ (as the mesh $\eta$ goes to zero), one would have to show that in expectation, the total mass $\param_\eta(\closure{\Omega})$ remains bounded as $\eta\to 0$.
In the bulk, this is done very similarly to Subsection~\ref{ss.tightness}, but near the boundary, this would already require some care. Still, even if we succeeded in showing this, it could still be that at the scaling limit, a positive fraction of the measure
would be concentrated on the boundary $\p \Omega$. This is a scenario that we definitively want to exclude --- otherwise, the measure could look quite strange. Therefore, we need to 
prove that the measure gets smaller near the boundary. More precisely, for any $\delta>0$, if $\Omega^{(\delta)}$ denotes the points in $\Omega$
at distance at least $\delta$ from the boundary, then we would like to show that 
\[
\param_\eta(\closure{\Omega} \setminus \Omega^{(\delta)}) = o(1)\,,
\]
where $o(1)$ goes to zero when $\delta\to 0$, uniformly in $\eta<\delta/10$. Note that this implies tightness in $\closure{\Omega}$, as well.

The proof of this estimate, though not conceptually hard, needs some rather tedious computations. We will here only highlight briefly how it works.
(See \cite{\GPS1} where many such boundary issues were treated carefully.) 
Divide the set $\closure{\Omega}\setminus \Omega^{(\delta)}$ into $\asymp 1/\delta$ squares of size $\delta$. Roughly speaking (ignoring the effects of the possible
curvature of $\p \Omega$), group these squares in a dyadic fashion depending on their distance from the set $\{a,b\}$. This singles out the two squares which touch
the tips $a$ and $b$. Group further all the sites in these two squares dyadically depending on their distance from $a$ and $b$, respectively. In all the other $\delta$-boxes, group the sites dyadically according to their distance from $\p\Omega$. In the $\delta$-box at $a$, at distance about $2^k\eta$ from $a$, the probability of a site to be on $\gamma_\eta$ is roughly $\alpha_2^\eta(\eta,2^k\eta)$, and there are roughly $2^{2k}$ such sites. So, their contribution to the un-normalized expected length measure is 
\begin{equation*}
\asymp \sum_{k=1}^{O(\log (\delta/\eta))} \alpha_2^\eta(\eta,2^k\eta) \, 2^{2k} \asymp \eta^{-2}\alpha_2^\eta(\eta,1) \, \delta^{7/4+o(1)}\,.
\end{equation*}
In the $\delta$-boxes at unit order distance from $a$ and $b$, for a site at distance about $2^k\eta$ from $\p\Omega$, the probability of being on $\gamma_\eta$ is roughly $\alpha_2^\eta(\eta,2^k\eta)\,\alpha_1^{+,\eta}(2^k\eta,1)$. In all these $\delta$-boxes, there are about $2^k/\eta$ such sites. So, their total contribution to the expectation is 
\begin{align*}
\asymp \sum_{k=1}^{O(\log (\delta/\eta))} \alpha_2^\eta(\eta,2^k\eta) \, \alpha_1^{+,\eta}(2^k\eta,1) \, 2^k / \eta 
& \asymp \alpha_2^\eta(\eta,\delta) \, \alpha_1^{+,\eta}(\delta,1) \, \delta/\eta^2 \\
&\asymp \eta^{-2}\alpha_2^\eta(\eta,1)\, \delta^{-1/4+1/3+1+o(1)}\\
& = \eta^{-2}\alpha_2^\eta(\eta,1)\, \delta^{13/12+o(1)}\,.
\end{align*}
After normalization by $\eta^{-2}\alpha_2^\eta(\eta,1)$, both contributions go to 0 with $\delta\to 0$. (For instance, the rate $\delta^{13/12}$ corresponds to fact that the $7/4$-dimensional curve $\gamma$ has a $2/3$-dimensional intersection with the boundary, so in this 2/3-dimensional subset the ``$7/4$-dimensional measure'' should scale with exponent  $7/4-2/3=13/12$.) The contributions of the $\delta$-boxes whose distance from $\{a,b\}$ is less than order 1 but more than zero is handled by a dyadic grouping again, and we are done.
\medskip

\ni {\bf Characterization of the subsequential scaling limits.}
Here, the whole setup carries through easily. Let $B$ be a ball inside $\Omega$; we wish to approximate $\param_\eta(B)$ using macroscopic information. Let us then define $Y=Y^\eps_\eta$ to be the number of $\eps$-squares $Q$ of $B$ (in the same grid $G$ as in Section~\ref{s.measurable}, but let us not worry here about the orientation of the grid) which are such that the square $2Q$ is connected via an open arm to the arc $ab$ and is connected via a closed arm to the opposite arc $ba$. Note that this is the same as asking that $\gamma_\eta \cap 2Q \neq \emptyset$.
The number we wish to approximate is 
\[\param_\eta(B) = \frac{X}{\eta^{-2} \alpha_4^\eta(\eta,1)}\,,
\]
where $X=X^\eps_\eta$ is the number of edges $e\in \eta \Tg \cap B$ which lie on the interface $\gamma_\eta$.

It is straightforward to adapt the coupling argument (Section~\ref{s.coupling}) to the case of the two-arm event. Basically the only change is that one deals with configurations of two faces, instead of four. Using this coupling result and following the details of Section~\ref{s.measurable}, one ends up 
with
\[
\Eb{(X- \beta_{\mathrm{two-arm}}\, Y)^2} = o(\Eb{X^2})\,,
\]
where the proportionality factor $\beta_{\mathrm{two-arm}}$ is the analog of $\beta$ in Section~\ref{s.measurable}.
\medskip

There is one part though which needs some more words of explanation: 
for the pivotal measure, we used the fact that $Y^\eps_\eta$ converged to $Y^\eps=Y^\eps(\omega)$. This was a direct consequence of Lemma~\ref{l.meas4arm}. Since, we do not rely on the quad-crossing topology in this subsection, we need an analog of Lemma~\ref{l.meas4arm} in our present setup. Namely, we need the following result:

\begin{proposition}
If $Q$ is any square included in $(\Omega,a,b)$, then 
\[
\Pb{\gamma_\eta \cap Q \neq \emptyset} \underset{\eta\to 0}{\longrightarrow} \Pb{\mathrm{SLE}_6 \cap Q \neq \emptyset}\,.
\]
Moreover, in any coupling with $\gamma_n\to\mathrm{SLE}_6$ a.s., we have the convergence in probability of the indicator variables of the events in question.
\end{proposition}

This ``continuity result'' about the convergence of discrete exploration paths towards $\mathrm{SLE}_6$ is somewhat {folklore} in the community. A sketch of this proposition can be found in \cite{\SmirnovWerner}. The complete proof is not hard but a bit technical: the main idea is to use the $3$-arm event in $\H$ in order to justify that it is very unlikely for an interface $\gamma_\eta$ or an $\mathrm{SLE}_6$ to come close to $\p Q$, without intersecting $Q$.
\QED

\subsection[Hints for the {\it frontier measure} and the $\mathrm{SLE}_{8/3}$ curve]{Hints for the {\it frontier measure} and a natural parametrization of $\mathrm{SLE}_{8/3}$ curves}\label{ss.frontier}

\begin{minipage}{\textwidth}
\begin{center}
\includegraphics[width=0.5 \textwidth]{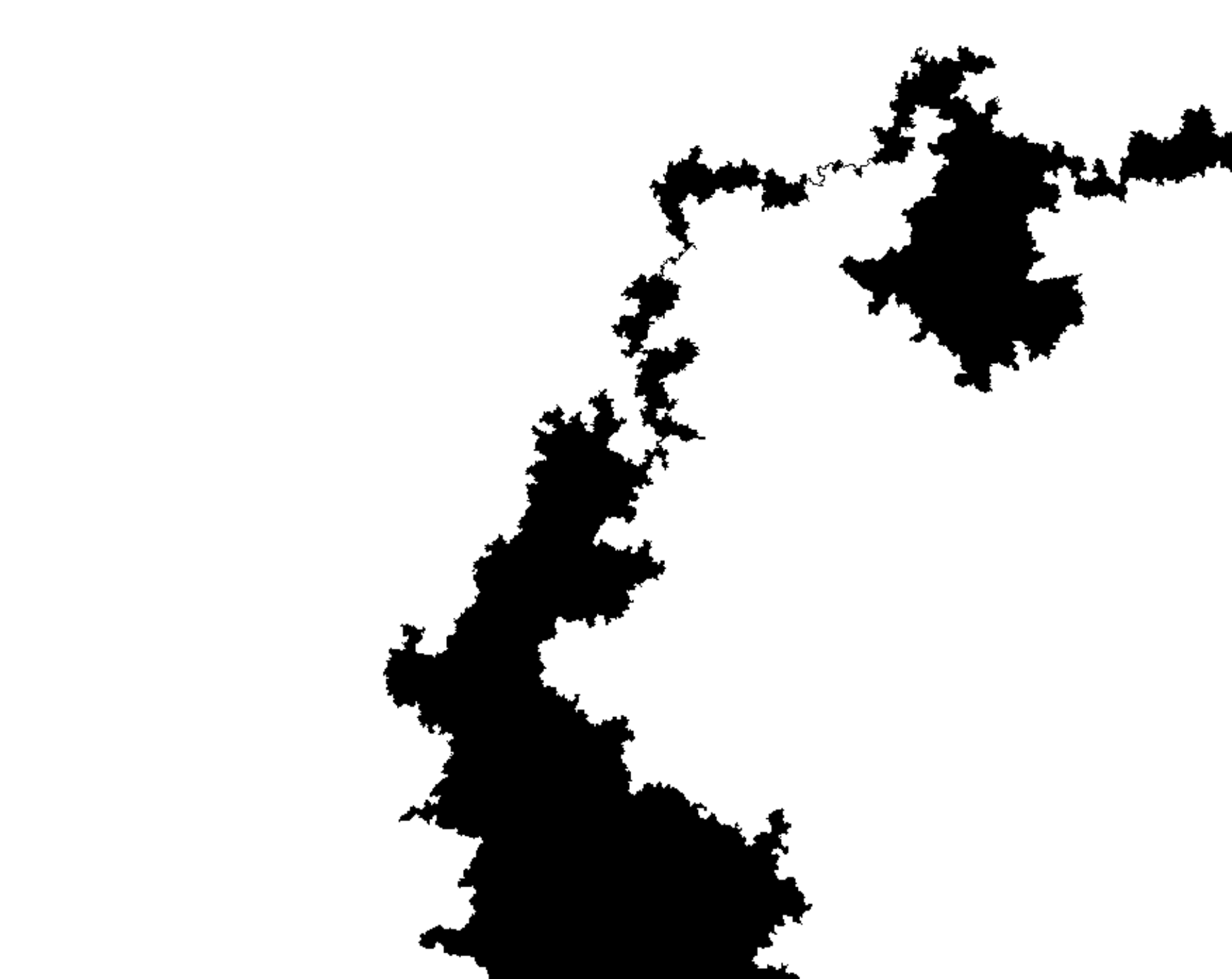}
\end{center}
\end{minipage}
\vskip 0.5 cm

Here we would want the normalized counting measure on the {\bf percolation frontier}, or {\bf exterior boundary} of large clusters. Namely, we say that an open hexagon neighbouring the interface $\gamma_\eta$ is on the open frontier $\Psi^+$ if it has not just one, but two disjoint open arms to the open boundary arc. The closed frontier $\Psi^-$ is  defined analogously. The right normalization for these points is clearly the polychromatic 3-arm event:
\begin{equation}\label{e.frontier}
\psi^{\pm}_\eta = \psi^{\pm}_\eta(\omega_\eta) = \psi^{\pm}_\eta(\gamma_\eta) =  \psi^{\pm}_\eta(\Psi^{\pm}_\eta):= \sum_{x \in \Psi^{\pm}}  \delta_x \, \frac {\eta^2} {\alpha_3^\eta(\eta,1)}\,.
\end{equation}

The scaling limit for this measure would need a bit more work than the other cases. The issue is that three arms is much less appropriate for our setup of faces $\Theta$ in Section~\ref{s.coupling}: having two arms of the same colour next to each other means that the 3-arm event is not characterized purely by interfaces that would produce faces and would cut information between past (a percolation configuration outside a large box) and future (the arrival of interfaces or arms at the small target box). The solution is to take a mixture of the approaches in the 2-arm and 1-arm cases, and use two interfaces plus RSW circuits to cut the information between the two interfaces on that side where we need the two arms of the same colour. Another small issue is that the separation of interfaces results of Appendix~\ref{s.appendix} cannot be used automatically: instead of an even number of interfaces, we would need to work with three well-chosen arms.

We do not go into the details, mainly because we do not have good applications in mind. In what $\sigma$-algebra would one want to have the limit measure be measurable? The quad-crossing configuration $\omega$ is one possibility, for which the 3-arm case of Lemma~\ref{l.meas4arm} can be used. This would be useful to measure the size of the exterior boundary of a large percolation cluster, which is a reasonable but much less important goal than the pivotal measure. On the other hand, measurability w.r.t.~to the frontier itself (defined as a curve) would be less nice than the case of interface measure, because the frontier has a scaling limit that is closely related to SLE$_{8/3}$, but is not that exactly. Getting a natural parametrization of SLE$_{8/3}$ by tweaking the frontier measure would certainly be much less natural and exciting than getting it, say, from Self-Avoiding Walk paths. Finally, measurability w.r.t.~the interface scaling limit $\mathrm{SLE}_6$ could possibly be interesting from some SLE duality point of view, but this again does not seem to be such a natural choice. 

\section[Conformal covariance of the limit measures]{Conformal covariance properties of the limit measures}\label{s.covariance}



We first discuss the annulus-important and quad-pivotal measures $\mu^A$ and $\mu^\Quad$, then the cluster and interface measures $\aream^A$ and $\param^{\Omega,a,b}$ in Subsection~\ref{ss.cothers}.

Let $\Omega, \widetilde\Omega$ be two simply connected domains of the plane, and let $f:\Omega\longrightarrow \widetilde\Omega$ be some conformal
map. By conformal invariance, the image $\tilde\omega:= f(\omega)$ is also a realization of continuum percolation in $\widetilde\Omega$; see Subsection~\ref{ss.black}. Consider some proper annulus $A \subset \bar A\subset \Omega$. Since $f$ is conformal on $\Omega$, we have that $f(A)$ is again a proper annulus, and $\mu^{f(A)}$ is the scaling limit of $\mu^{f(A)}_\eta$. We will prove the following:

\begin{theorem}\label{pr.covariance}
Let $f_*(\mu^A(\omega))$ be the pushforward measure of $\mu^A$. Then, for almost all $\omega$, the Borel measures $\mu^{f(A)}(\tilde\omega)$ and $f_*(\mu^A(\omega))$  on $f(\Delta)$ are absolutely continuous w.r.t.~each other, and their Radon-Nikodym
derivative satisfies, for any $w=f(z) \in \widetilde\Omega$,
\begin{equation}
\frac {d \mu^{f(A)}(\tilde\omega)}{ d f_*(\mu^A(\omega))}\,(w) = |f'(z))|^{3/4}, \nonumber
\end{equation}
or equivalently, for any Borel set $B\subset \Delta$,
\begin{equation}
\mu^{f(A)}(f(B))(\tilde\omega) = \int_{B} |f'|^{3/4} d\mu^A(\omega) \,. \nonumber
\end{equation}
\end{theorem}

\begin{remark}\label{r.inlaw}
Using the conformal invariance of the scaling limit, from the almost sure equality in the theorem we get that $\mu^{f(A)}(f(B))$ has the same law as $\int_{B} |f'|^{3/4} d\mu^A$.
\end{remark}

\medskip

Since any conformal map is locally a rotation times a dilatation, we will, as a warm-up, first check the theorem on these particular cases. This will be easier than the general case, mainly because the grid of $\eps$-squares that we used in defining the approximating macroscopic quantities $Y^\eps$ is preserved quite nicely under rotations and dilatations, while distorted by a general conformal map $f$.

\subsection{Rotational invariance}\label{ss.rot}

Let us consider some proper annulus $A$ of the plane, some ball $B\subset \Delta$, and the rotation $T:\, z\mapsto \exp(i\theta)\,z$ by an angle $\theta$. We need to show that
\begin{equation}
\mu^A(B,\omega)=\mu^{TA}(TB,T\omega)\,.
\label{e.rotinv}
\end{equation}

By~(\ref{e.muY}), the right hand side equals
$\lim_\eps Y^{G(\eps),TA}(TB,T\omega)/ F(\eps)$, an $L^2$-limit,
where $G(\eps)$ is our usual grid of $\eps$-squares and $F(\eps):=\eps^{-2}\alpha_4(\eps,1) / c$. Or, by using a rotated grid for the approximation in~(\ref{e.muY}), it is also equal to $\lim_\eps Y^{TG(\eps),TA}(TB,T\omega) / F(\eps)$, which (by rotating back the entire universe) is trivially the same as $\lim_\eps Y^{G(\eps),A}(B,\omega)/F(\eps)$, giving the left hand side of~(\ref{e.rotinv}), as desired.

\begin{remark}
One might speculate that this type of rotational invariance should hold even if we have a scaling limit $\omega$ without (or unproved) rotational invariance (such as a subsequential limit of critical percolation on $\Z^2$), since the backbone of the above argument seems to be the following: $\mu^A$ is a function of $\omega$, and the definition of this function in~(\ref{e.muY}) does not depend on any special orientation $\theta$, hence if one rotates $\omega$, then $\mu^A$ should get rotated, as well. (Of course, rotational invariance of the law of $\omega$ would still be essential for the equality in law discussed in Remark~\ref{r.inlaw}.)

However, the rotational invariance of the scaling limit is in fact used here, in a somewhat implicit way, through the fact that the normalization factor $F(\eps)$ in~(\ref{e.muY}) does not depend on $\theta$, which was one of the implications of Proposition~\ref{pr.ratio}.


We will see the same phenomenon in the proofs below: the conformal covariance comes in some sense from the fact that the normalization factor cannot be changed when applying a conformal map.
\end{remark}

\subsection{Scaling covariance}\label{ss.scaling}

We show here the following special case of Theorem~\ref{pr.covariance}:

\begin{proposition}
\label{pr.scalingcov}
Let $A$ be some proper annulus of the plane and $\lambda>0$ some scaling factor. Then, for any $B\subset \Delta$:
\begin{equation}
\label{e.scaling}
 \mu^{\lambda A}(\lambda B,\lambda\omega) = \lambda^{3/4} \mu^A(B,\omega)\,.
 \end{equation}
\end{proposition}

\proof
By (\ref{e.muY}), we have the $L^2$-limits
$$
\mu^{\lambda A}(\lambda B,\lambda\omega)=\lim_{\eps\to 0}
\frac{Y^{\lambda\eps,\lambda A}(\lambda B,\lambda\omega)}{F(\lambda\eps)}
= \lim_{\eps\to 0} \frac{Y^{\eps,A}(B,\omega)}{F(\lambda\eps)}\,,
$$
where the second equality is just a tautology even for fixed $\eps>0$. So, in order to get the right side of~(\ref{e.scaling}), we need that
$\lim_{\eps\to 0} \alpha_4(\lambda\eps,1)/\alpha_4(\eps,1)=\lambda^{5/4}$. By the scale invariance of the scaling limit,
$\alpha_4(\lambda\eps,1) = \alpha_4(\eps,1/\lambda)$, so the ratio limit result~(\ref{e.limitlambda}) finishes the proof.\qed

\subsection{Proof of conformal covariance}


The key step in the proof of Theorem \ref{pr.covariance} will be the following:

\begin{lemma}\label{l.covariance}
There are absolute constants $\delta_0=\delta_0(A,\Omega,\widetilde\Omega)>0$ and $K= K(A,\Omega,\widetilde\Omega)>0$
so that almost surely (with respect to $\omega$), for any ball or square $B\subset \Delta$ centered at $z$ of radius
$\delta \leq \delta_0$, we have
\begin{equation}
(1 - K \delta) \mu^A (B,\omega) |f'(z)|^{3/4} \leq \mu^{f(A)} (f(B),\tilde\omega) \leq (1+ K\delta) \mu^A (B,\omega) |f'(z)|^{3/4}.
\nonumber
\end{equation}
\end{lemma}

Indeed, let us briefly show that Lemma  \ref{l.covariance} implies that a.s.~$\mu^{f(A)} \ll f_*(\mu^A)$. Take any Borel set $U$ such that
$f_*(\mu^A)(U) = 0$. For any $\eps>0$ there is some finite cover $\bigcup_i B_i$ of $f^{-1} (U)$ (also a Borel set) by disjoint open squares $B_i$, each of radius less than $\delta_0$, so that
$\mu^A(\bigcup_i B_i) \leq \eps$. Therefore, by Lemma~\ref{l.covariance},
\begin{align*}
\mu^{f(A)}(U)  & \leq  ( 1+ K \delta_0)\, \sup_{z\in \Delta} |f'(z)|^{3/4}\, \mu^A\big(\bigcup_i B_i\big) \\
 & \leq  (1+ K \delta_0) \, H^{3/4} \, \eps\,,
\end{align*}
where $H:= \sup_{\Delta} |f'|$ is finite, since $\bar \Delta$ is a compact set inside $\Omega$ where $f$ is conformal. By letting $\eps$
go to zero, this proves $\mu^{f(A)} \ll f_*(\mu^A)$. The other direction is proved in the same way. Therefore, the two measures
are absolutely continuous, and it is straightforward from Lemma~\ref{l.covariance} that their Radon-Nikodym derivative is indeed as in
Theorem~\ref{pr.covariance}.

\proofof{Lemma \ref{l.covariance}}
Since $\bar A$ is a compact subset of $\Omega$, we can define
\begin{equation} \label{e.2derivativeF}
\begin{array}{lcr}
H_f &:=& \sup_{z\in A} |f'(z)| <\infty \\
L_f &:=& \sup_{z\in A} |f'' (z)| <\infty .
\end{array}
\end{equation}
Since also $\overline{f(A)}$ is a compact subset of $\widetilde\Omega$, we have
\begin{equation} \label{e.2derivativeG}
\begin{array}{lcr}
H_g &:=& \sup_{w\in f(A)} |g'(w)| <\infty \\
L_g &:=& \sup_{x\in f(A)} |g'' (w)| <\infty\,,
\end{array}
\end{equation}
where $g:\widetilde\Omega \longrightarrow \Omega$ is the inverse of $f$.

We will fix the value of $\delta_0>0$ later on.
Let $B=B(z_0,\delta)$ some ball of radius $\delta \leq \delta_0$ centered at $z_0$ and satisfying $B\subset \Delta$.

For any parameters $a\in \C,\theta\in[0,2\pi)$, let $G$ be the grid of $\eps$-squares centered at $a$ and rotated by $e^{i\theta}$.
As in the previous sections, $Y^{\eps,a,\theta} = Y^{\eps,a,\theta}(\omega)$
will be the random variable corresponding to the number of $G$-squares $Q$ inside $B$ for which $2Q$ is $A$-important
for the configuration $\omega$. Recall that (\ref{e.muY}) gives an approximation to $\mu^A(B)$ using $Y^{\eps,a,\theta}$, with a speed of convergence that is independent of $a,\theta$ and also of the ball $B$.

In particular, for any $\eps>0$, if $\pi^{\eps}$ is any probability measure on the parameters $a,\theta$, we obtain the $L^2$-limit
\begin{equation}\label{e.measurepi}
\mu^A(B) = \lim_{\eps\rightarrow 0} \int_{a,\theta} \frac{Y^{\eps,a,\theta}} {(2\eps)^{-2} \alpha_4(2\eps,1)}\, d\pi^{\eps}(a,\theta).
\end{equation}

In our setup we make the natural choice to define $\pi^{\eps}$ as the product measure of
the normalized uniform measure on $[-\eps,\eps]^2$ for $a$, times the normalized uniform measure on $[0,2\pi]$ for $\theta$.
With this particular choice, it turns out that one can rewrite (\ref{e.measurepi}) in a nicer way.
First let us define for any $z\in B$ and any $\eps,\theta$ the random variable $X_{\theta}^{\eps}(z)$ to be the indicator
function of the event that the square of radius $2\eps$ centered at $z$ and rotated by $e^{i\theta}$ is $A$-important.
We will show the following lemma:


\begin{lemma}\label{l.rewriteY} We have the $L^2$-limit
\begin{equation}
\mu^A(B) = \lim_{\eps\rightarrow 0} \int_{B\times [0,2\pi)}  \frac{X_{\theta}^{\eps}(z)}{\alpha_4(2\eps,1)}\, d\A(z) \,d\mathcal{L}(\theta)\,,
\nonumber
\end{equation}
where $d\A$ is the (non-renormalized) area measure on $B$ and $d\mathcal{L}$ is the normalized Lebesgue measure on $[0,2\pi]$.
\end{lemma}

\proof Recall that $B=B(z_0,\delta)$. It is straightforward to check that, by the definition of $Y^{\eps,a,\theta}$, we have
\begin{multline*}
\int_{\scriptsize\begin{array}{l} B(z_0, \delta -4 \eps) \\ (0,2\pi)\end{array}} X^{\eps}_{\theta}(z)\, d\A(z)\, d\mathcal{L}(\theta) \leq
 \int_{\scriptsize\begin{array}{l}[-\eps,\eps]^2 \\ (0,2\pi)\end{array}}   Y^{\eps,a,\theta}\, d\A(a) \,
d\mathcal{L}(\theta) \\
\leq  \int_{\scriptsize\begin{array}{l} B(z_0, \delta) \\ (0,2\pi)\end{array}} X^{\eps}_{\theta}(z)\, d\A(z) \, d\mathcal{L}(\theta).\nonumber
\end{multline*}
Since we have $d\pi^{\eps}(a,\theta) = \frac{1}{4\eps^4}\, d\A(a) d\mathcal{L}(\theta)$ with the above choice
of the measure $\pi^{\eps}$, the above inequalities can be rewritten as
\begin{multline*}
\int_{\scriptsize\begin{array}{l} B(z_0, \delta -4 \eps) \\ (0,2\pi)\end{array}} \frac {X^{\eps}_{\theta}(z)}{\alpha_4(2\eps,1)}\, d\A(z) \, d\mathcal{L}(\theta) \leq
 \int_{\scriptsize\begin{array}{l}[-\eps,\eps]^2 \\ (0,2\pi)\end{array}}   \frac{Y^{\eps,a,\theta}}{(2\eps)^{-2} \alpha_4(2\eps,1)}\, d\pi^{\eps}(a,\theta) \\
\leq \int_{\scriptsize\begin{array}{l} B(z_0, \delta) \\ (0,2\pi)\end{array}} \frac{X^{\eps}_{\theta}(z)}{\alpha_4(2\eps,1)}\, d\A(z) \, d\mathcal{L}(\theta).
\end{multline*}

So, it is enough to prove that the boundary effect
$$
\mathcal{E}:=  \int_{\scriptsize\begin{array}{l} B(z_0, \delta)\setminus B(z_0,\delta-4\eps) \\ (0,2\pi)\end{array}}
 \frac {X^{\eps}_{\theta}(z)}{\alpha_4(2\eps,1)}\, d\A(z) \, d\mathcal{L}(\theta)
$$
is negligible when $\eps$ goes to zero. For each $z\in \Delta$,
the probability that $X_{\theta}^{\eps}(z)$ equals 1
is of order $O(1)\alpha_4(\eps,1)$ (where $O(1)$ only depends on $A$). Since the area of $B(z_0,\delta)\setminus B(z_0,\delta- 4\eps)$
is of order $\delta \eps$, altogether we have
\[
\Eb{\mathcal{E}} \leq O(1) \delta \eps,
\]
which completes the proof of Lemma \ref{l.rewriteY}.
\QED

For any $z,\theta$ and $\eps>0$, $B^{\theta}(z,\eps)$ will denote the square of radius $2\eps$, centered at $z$ and rotated
by $e^{i\theta}$; in particular $X_{\theta}^{\eps}(z) = 1_{B^{\theta}(z,\eps)\textrm{ is $A$-important}}$.

Lemma \ref{l.rewriteY} says that
$$
\mu^A(B) = 
\lim_{\eps\to 0} \frac{1}{\alpha_4(2\eps,1)} \int_{B\times [0,2\pi)}  1_{B^{\theta}(z,\eps)\textrm{ is $A$-important}} \,d\A(z)\,d\mathcal{L}(\theta).
$$
Let us change variables in the following way:
\[
\begin{array}{lcl}
\left\lbrace \begin{array}{lcl}
\tilde z &=& f(z)\\
\tilde \theta & = & \theta + \Im(\log f'(z))
\end{array}\right.
& \text{ or equivalently } &
\left\lbrace \begin{array}{lcl}
z &=& g(\tilde z)\\
\theta & = & \tilde \theta + \Im(\log g'( \tilde z)).
\end{array}\right.
\end{array}
\]
The Jacobian of the change of variables $(\tilde z, \tilde \theta) \mapsto ( z, \theta)$, from $\R^3$ to $\R^3$, is $|g'(\tilde z)| ^2$, one therefore has
\begin{align}\label{e.changeofvariable}
\int_{B\times [0,2\pi)}  1_{B^{\theta}(z,\eps)\textrm{ is $A$-important for $\omega$ }} \,d\A(z)\,d\mathcal{L}(\theta) \\
&  \hskip -6 cm  = \int_{f(B) \times [0.2\pi)}  1_{f(B^{\theta}(z,\eps)) \textrm{ is $f(A)$-important for $\tilde{\omega}$ }} \, |g'(\tilde z)|^2 
\, d\A(\tilde z) \, d\mathcal{L}(\tilde \theta) \,, \nonumber
\end{align}
since $\tilde{\omega}$ is the continuum percolation satisfying $\tilde{\omega} =f(\omega)$.

Now, for any $z\in B= B(z_0, \delta)$, by the definition of $H_f$, we have that $|f(z)- f(z_0)|  < H_f |z-z_0| \leq H_f \delta$;
hence, if $\tilde{z_0}:= f(z_0)$, then $f(B)\subset B(\tilde{z_0}, H_f \delta)$.
Now, by the definition of $L_g$, for any $\tilde z \in f(B)$, we have that $|g'(\tilde z)| \leq |g'(\tilde{z_0})| + L_g H_f \delta$;
here one needs to take $\delta$ small enough so that $B(\tilde{z_0}, H_f \delta)$ is still included in $f(A)$.
This gives, for any $\tilde z \in f(B)$,
\begin{align}
|g'(\tilde z)|^2  &\leq  |g'(\tilde{z_0})|^2 + 2 L_g H_g H_f \delta + L_g^2 H_f^2 \delta^2  \nonumber\\
& \leq  |g'(\tilde{z_0})|^2 + O(1) \delta \, . \label{e.gprime}
\end{align}
Similarly, we have that $|g'(\tilde z)|^2 \geq |g'(\tilde{z_0})|^2 - O(1)\delta$.

Now notice that the $2\eps$-squares are very little distorted by $f$. Indeed, consider some square $B^{\theta}(z,\eps)$ (recall that it is the square of radius $2\eps$, centered at $z$, rotated by $\theta$); for any point $u\in B^{\theta}(z,\eps)$ we have
\[
|f(u)- f(z) - f'(z)(z-u) | \leq L_f \frac{(4\eps)^2}{2},
\]
since $|u-z|\leq \sqrt{2}\, 2 \eps \leq 4\eps$. Therefore, if
\[
\left\lbrace\begin{array}{lclcl} \tilde{\eps_1} &=&\tilde{\eps_1}(z)& =& |f'(z)| \eps - 4\eps^2 L_f \\
\tilde{\eps_2} &=&\tilde{\eps_2}(z)& =& |f'(z)| \eps + 4\eps^2 L_f \, ,\end{array} \right.
\]
then
\begin{equation}
B^{\tilde{\theta}}(\tilde{z}, \tilde{\eps_1}) \subset f(B^{\theta}(z,\eps)) \subset B^{\tilde{\theta}} (\tilde z, \tilde{\eps_2}).
\label{e.BfBalls}
\end{equation}
This and (\ref{e.gprime}) imply the following upper bound (and the lower bound would work in a similar way)
\begin{eqnarray}
\int_{f(B) \times [0.2\pi)}  1_{f(B^{\theta}(z,\eps)) \textrm{ is $f(A)$-important for $\tilde{\omega}$ }} \, |g'(\tilde z)|^2 
\, d\A(\tilde z) \, d\mathcal{L}(\tilde \theta) \nonumber\\
& & \hskip -8cm \leq \int_{f(B)\times [0,2\pi) } 1_{B^{\tilde{\theta}}(\tilde z,\tilde{\eps_2}) \textrm{ is $f(A)$-important for $\tilde{\omega}$ }} \, |g'(\tilde z)|^2 
\, d\A(\tilde z) \, d\mathcal{L}(\tilde \theta)  \nonumber \\
& & \hskip -8cm \leq \int_{f(B)\times [0,2\pi) }X_{\tilde{\theta}}^{\tilde{\eps_2}}(\tilde z) \,  \big( |g'(\tilde{z_0})|^2 + O(1) \delta \big)\, d\A(\tilde z)\, d\mathcal{L}(\tilde{\theta}). \nonumber
\end{eqnarray}
Combined with (\ref{e.changeofvariable}), and using that there is a uniform lower bound on $|g'(\tilde{z_0})|^2$,
this leads to
\begin{align}
\label{e.eps2bound}
\int_{ B\times [0,2\pi) } \frac {X_{\theta}^{\eps}(z)}{\alpha_4(2\eps,1)}\, \, d\A(z) \, d\mathcal{L}(\theta)\\
& \hskip - 4cm \leq |g'(\tilde{z_0})|^2 \big(1+ O(1)\delta\big) \,\frac{\alpha_4(2 \tilde{\eps_2},1)}{\alpha_4(2\eps,1)}\,
\int_{f(B) \times [0,2\pi)} \frac {X_{\tilde{\theta}}^{\tilde{\eps_2}}(\tilde{z})}{\alpha_4(2\tilde{\eps_2},1)} \,d\A(\tilde{z})\,d\mathcal{L}(\tilde{\theta})\nonumber.
\end{align}

The ratio limit result (\ref{e.limitlambda}) and the scale invariance of $\alpha_4(\cdot,\cdot)$ imply that
\[
\frac{\alpha_4(2\tilde{\eps_2},1)}{\alpha_4(2\eps,1)} = \frac{\alpha_4\big( 2 |f'(z)|\eps + 8 \eps^2 L_f, 1\big) }{\alpha_4(2\eps, 1)}
\overset{\eps\to 0}{\longrightarrow} |f'(z)|^{5/4}.
\]

Therefore, by letting the mesh $\eps$ (hence also the mesh $\tilde{\eps_2}$) go to 0 in (\ref{e.eps2bound}),
and using Lemma \ref{l.rewriteY} in both domains $\Omega, \widetilde\Omega$, one ends up with
\begin{eqnarray}
\mu^A(B) &\leq &  |g'(\tilde{z_0})|^2 |f'(z)|^{5/4} \big(1+ O(\delta)\big) \, \mu^{f(A)}(f(B))  \nonumber\\
 & \leq & |g'(\tilde{z_0})|^{3/4} \big(1+ O(\delta)\big)\, \mu^{f(A)}(f(B)) \nonumber,
\end{eqnarray}
since we have a uniform control on how $|f'(z)|$ is close to $|f'(z_0)|= |g'(\tilde{z_0})|^{-1}$ on the ball $B=B(z_0,\delta)$.

This, together with the lower bound that is proved in the same way, completes the proof of Lemma~\ref{l.covariance}: the last thing to do is to choose the threshold radius $\delta_0$ to be small enough so that for any $z\in \Delta$, the ball $B(z, H_f \delta_0)$ is still inside $f(A)$. \QED
\medskip

Given Theorem~\ref{th.4meas}, the above proof for Theorem~\ref{pr.covariance} works also for Theorem~\ref{th.4cov}. 
A detail that might be worth pointing out is why it is that the covariance exponent is given by the 4-arm exponent also close to $\p\Quad$ and not by the half-plane 3-arm event (and the 2-arm event at the corners). The reason is that although the $\eps$-boxes counted in the macroscopic approximation $Y^\eps$ that are close to $\p\Quad$ indeed need to have more-or-less the half-plane 3-arm event, but, on this macroscopic scale, we can use conformal {\it in\/}variance combined with the effect of the normalization constant being still the whole plane 4-arm probability. This is the same effect that was pointed out in the last paragraph of Subsection~\ref{ss.rot} and was seen in Subsection~\ref{ss.scaling}. This works fine because the normalization constant comes from the microscopic scale, and there, inside the $\eps$-box, again the whole-plane 4-arm event is the relevant.

\subsection{Conformal covariance for the cluster and interface measures}\label{ss.cothers}

Here are the results concerning the cluster (area) measure $\aream^A$ and interface (length) measure $\param^{\Omega,a,b}$. 

Let $f: \Omega \longrightarrow  \widetilde\Omega$ be a conformal isomorphism of two simply connected domains of the plane, with continuum percolation configurations $\omega$ and $\tilde\omega=f(\omega)$. Consider some proper annulus $A \subset \bar A\subset \Omega$. Then $f(A)$ is again a proper annulus, and $\aream^{f(A)}$ is the scaling limit of $\aream^{f(A)}_\eta$.

\begin{theorem}\label{th.1cov}
Let $f_*(\aream^A(\omega))$ be the pushforward measure of $\aream^A$. Then, for almost all $\omega$, the Borel measures $\aream^{f(A)}(\tilde\omega)$ and $f_*(\aream^A(\omega))$  on $f(\Delta)$ are absolutely continuous w.r.t.~each other, and their Radon-Nikodym
derivative satisfies, for any $w=f(z) \in \widetilde\Omega$,
\begin{equation}
\frac {d \aream^{f(A)}(\tilde\omega)}{ d f_*(\aream^A(\omega))}\,(w) = |f'(z))|^{91/48}, \nonumber
\end{equation}
or equivalently, for any Borel set $B\subset \Delta$,
\begin{equation}
\aream^{f(A)}(f(B))(\tilde\omega) = \int_{B} |f'|^{91/48} d\aream^A(\omega) \,. \nonumber
\end{equation}
\end{theorem}

For the case of interfaces, let $\Omega$ be a smooth simply connected domain, $a,b \in\p\Omega$ two distinct points, $\gamma$ the continuum interface (a chordal $\mathrm{SLE}_6$) from $a$ to $b$, and $\param^{\Omega,a,b}$ be the corresponding interface (length) limit measure. Let $f: \Omega \longrightarrow  \widetilde\Omega$ be a conformal isomorphism extending continuously to $\p\Omega$ and $|f'|$ being bounded on $\Omega$. (For instance, if $\p\Omega$ is not only smooth, but Dini-smooth, then $f'(z)$ extends continuously to $\p\Omega$, see \cite[Theorem 10.2]{\PommerenkeUniv}.) Let $\tilde a:=f(a)$ and $\tilde b:=f(b)$, and $\tilde\gamma:=f(\gamma)$, which has the law of a continuum interface from $\tilde a$ to $\tilde b$, due to the conformal invariance of $\mathrm{SLE}_6$.

\begin{theorem}\label{th.2cov}
Let $f_*(\param^{\Omega,a,b}(\gamma))$ be the pushforward measure of $\param^{\Omega,a,b}$. Then, for almost all $\gamma$, the Borel measures $\param^{\widetilde\Omega,\tilde a,\tilde b}(\tilde\gamma)$ and $f_*(\param^{\Omega,a,b}(\gamma))$  on $f(\Omega)$ are absolutely continuous w.r.t.~each other, and their Radon-Nikodym
derivative satisfies, for any $w=f(z) \in \widetilde\Omega$,
\begin{equation}
\frac {d \param^{\widetilde\Omega,\tilde a,\tilde b}(\tilde\gamma)}{ d f_*(\param^{\Omega,a,b}(\gamma))}\,(w) = |f'(z))|^{7/4}, \nonumber
\end{equation}
or equivalently, for any Borel set $B\subset \Omega$,
\begin{equation}
\param^{\widetilde\Omega,\tilde a,\tilde b}(f(B))(\tilde\gamma) = \int_{B} |f'|^{7/4} d\param^{\Omega,a,b}(\gamma) \,. \nonumber
\end{equation}
\end{theorem}

Given the results and proofs of Section~\ref{s.others}, the proofs are straightforward modifications of Theorems~\ref{pr.covariance} and~\ref{th.4cov}.

\appendix

\section{Appendix: The strong separation phenomenon}\label{s.appendix}

In many planar statistical physics models, it is a very useful phenomenon that if one conditions on certain curves (trajectories or interfaces) to reach to a long distance without colliding with each other, then with a probability uniformly bounded away from zero, they arrive at their target with a good separation, i.e., with a distance from each other that is as large as it can be (up to a constant factor). Such a result was probably first proved and used in Kesten's fundamental paper \cite{\KestenScaling}: conditioned on the alternating four-arm event in $A(r,R)$, with positive probability there are four interfaces well-separated at distance $R$. He used this  to prove the quasi-multiplicativity of multi-arms probabilities. Cleaner proofs appear in \cite[Appendix]{\SchrammSteif} and \cite{\NolinKesten}. However, this is weaker than what we need here: we want the conclusion under conditioning on arbitrarily bad starting points of the four given interfaces at radius $r$. Such a strong separation lemma was proved for Brownian motions by Greg Lawler in \cite[Lemma 4.2]{\LawlerStrict}, which is a fundamental tool in computing Brownian intersection exponents \cite{\LSWan}. Then, using the Lawler-strategy, an analogous result was proved for a random walk approaching a level line of the discrete Gaussian free field in \cite[Lemma 3.15]{\SchShDGFF}, and for the meeting of a random walk and a loop-erased random walk in \cite[Theorem 4.7]{\Masson}, used again in \cite{\BarlowMassonTail}. This strategy could also be adapted to percolation interfaces; however, as was kindly explained to us by Art\"em Sapozhnikov, it is simpler to modify Kesten's argument. That approach has appeared in \cite{\Sapozhnikov}, but, for completeness, roughly following their argument but in a form more suitable for us, we include it here. The main ingredients of Kesten's strategy are the same as those of Lawler's, just the organization of the ideas is different; e.g., the latter uses random stopping times that the former does not. Consequently, in Kesten's strategy it is easier to separate the input that is needed from the model itself in order for the argument to work. To emphasize its universal aspect, in the second part of this Appendix we will use Kesten's strategy to prove Lawler's separation lemma for Brownian motions. We hope that this mostly self-contained Appendix will be useful for other models, as well.

The Lawler-strategy was streamlined and used in \cite{\LawlerVermesi} to show strong separation and coupling results for three-dimensional Brownian motions. This shows that the separation phenomenon is not restricted to planar models, even if proving the separation of interface surfaces in the three-dimensional critical percolation is presently just a dream, especially that the RSW technology does not really exist yet.

\subsection{Percolation interfaces}

Let us stress that the results in this section hold for critical percolation on any planar lattice where the RSW estimates are known; in particular, for bond percolation on $\Z^2$ and site percolation on $\Tg$.

Recall the notations from the beginning of Section~\ref{s.coupling}. Our main goal is a proof of Lemma~\ref{l.separation}, but we will first prove a more traditional version, where there is no ``forward conditioning'' (i.e., the conditioned target for the interfaces is at the same radius where we want the good separation), and there is no $U_\Theta=\tau$ conditioning.

\begin{proposition}[Strong Separation Lemma, \cite{\Sapozhnikov}]\label{p.separation}
Fix $N\in 2\Z_+$ for the number of interfaces considered, and let $A=A(R_1,R_2)$ be the annulus centered around 0 with radii $100 N \eta <  2 R_1 < R_2$. Assume we are given some
configuration of faces $\Theta=\{\theta_1,\ldots,\theta_N\}$ around the square of radius $R_2$, with some arbitrary initial quality $\Qual(\Theta)$.

Let $\mathcal{A}_\Theta (R_1,R_2)$ be the event that the $N$-tuple of interfaces in $\mathcal{D}_\Theta$ that start at the $N$ endpoints of $\Theta$ reach radius $R_1$, and let $\Gamma$ be the set of these interfaces. Then, there is some constant $c(N)>0$ such that 
\[
\PB{ \Qual(\Gamma) >\frac{1}{N} \md \mathcal{A}_\Theta (R_1,R_2)} > c(N)\,. 
\]
\end{proposition}

\proof
Let us start with a sketch of the main ideas. In the unconditional measure, it is unlikely that in a fixed dyadic annulus $A(s,2s)$ there are interfaces very close to each other; this is formalized in the well-known Lemma~\ref{l.nobadsep} below, proved by RSW estimates. Hence, still in the unconditional measure, the probability of fulfilling the conditioning but having very low quality configurations on $k$ consecutive dyadic scales  is exponentially small in $k$, with a base that can be made arbitrarily small if the ``low quality'' is chosen low enough. On the other hand, from a dyadic scale with a not very low quality, with some small probability the situation can be improved to a good quality in the next dyadic scale, and this good quality then can be preserved along several scales, with an exponentially small probability but with a uniformly positive base (by the also standard~(\ref{e.moveon}) and~(\ref{e.goon}) below). These suggest that conditioned on $\A_\Theta(R_1,R_2)$, if the starting quality at $R_2$ was not very bad, then good quality scales would appear with a positive frequency, and having a good quality at the end would have positive probability. However, we start with an arbitrary starting configuration, so we need to improve its quality little by little, with each step having a size comparable to the quality at that moment. 

Since the starting configuration might not be uniformly bad, the step sizes should also vary locally. Therefore, we build a simple hierarchical structure on the interface endpoints $V=\{x_1,\dots,x_N\}\subset B(R_2)$ of the starting configuration: we first try to improve the quality in small steps (in annuli on dyadic scales) for the interfaces that are very close, then other interfaces join in at larger scales given by the distances in the starting configuration, and so on, resulting in a tree structure. Here is the precise procedure:

Let the minimum distance between the points $x_i$ be $s$. For $j\in\N$, let $G_j$ be the graph on $V$ where two points are joined if their $\ell^\infty$-distance is less than $s2^j$; in particular, $G_0$ is the graph on $V$ with
no edges. If $j\in\N_+$ and $C\subseteq V$ is a connected component
of $G_j$, but is not connected in $G_{j-1}$, then $C$ is called a {\bf level $j$ cluster} of $V$. The level of a cluster is denoted by $j(C)$. Note that if two clusters intersect, then one contains the other.
There is an obvious tree structure on the clusters, with $V$ (as a single cluster) being the root, and the singletons $\{x_i\}$ being the leaves. The parent of a cluster $C$ will be denoted by $C^p$.

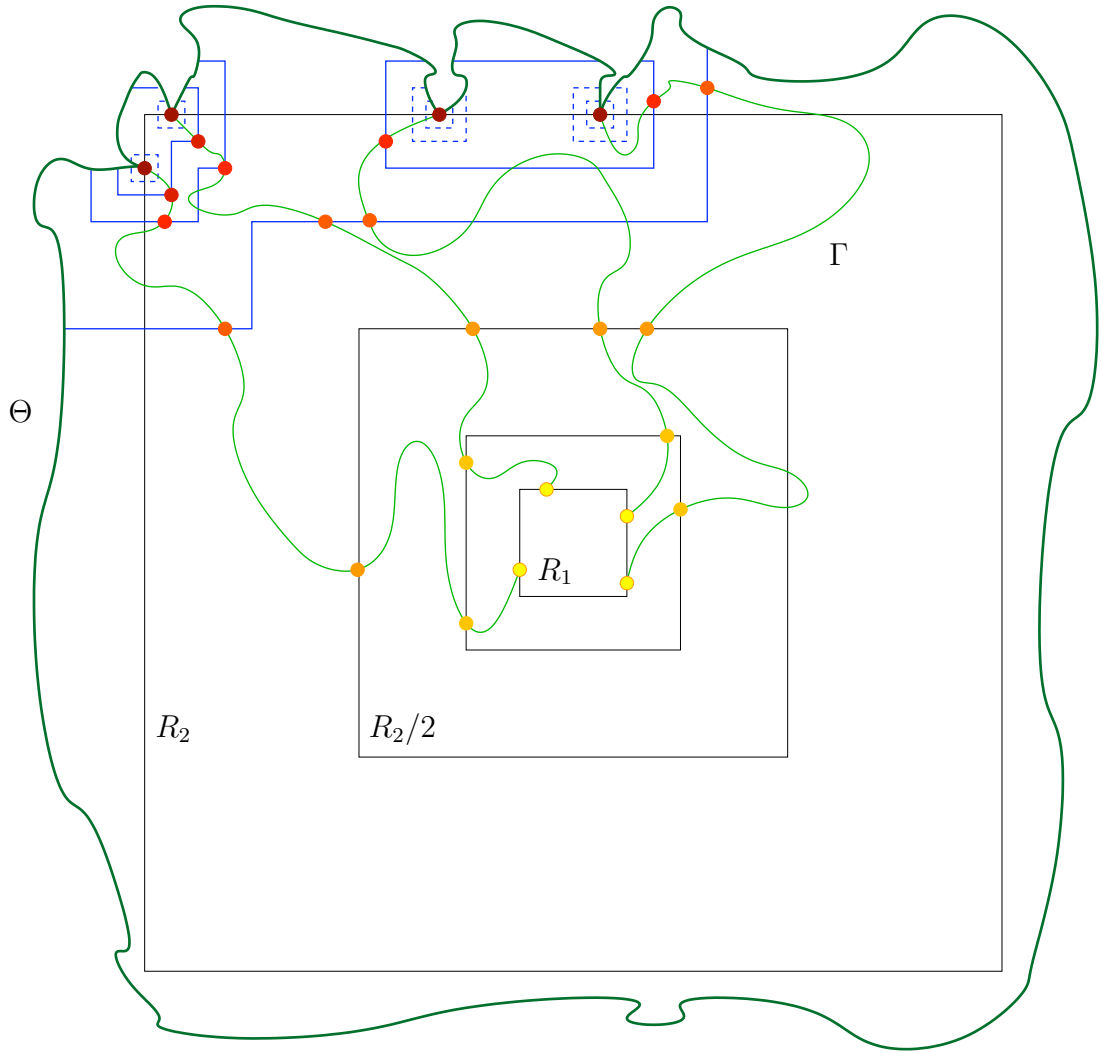
\begin{figure}[htbp]
\SetLabels
(0*.6)$\Theta$\\
(.75*.75)$\Gamma$\\
(.14*.3)$R_2$\\
(.35*.3)$R_2/2$\\
(.49*.45)$R_1$\\
\endSetLabels
\centerline{
\AffixLabels{
\includegraphics[height=5.5 in]{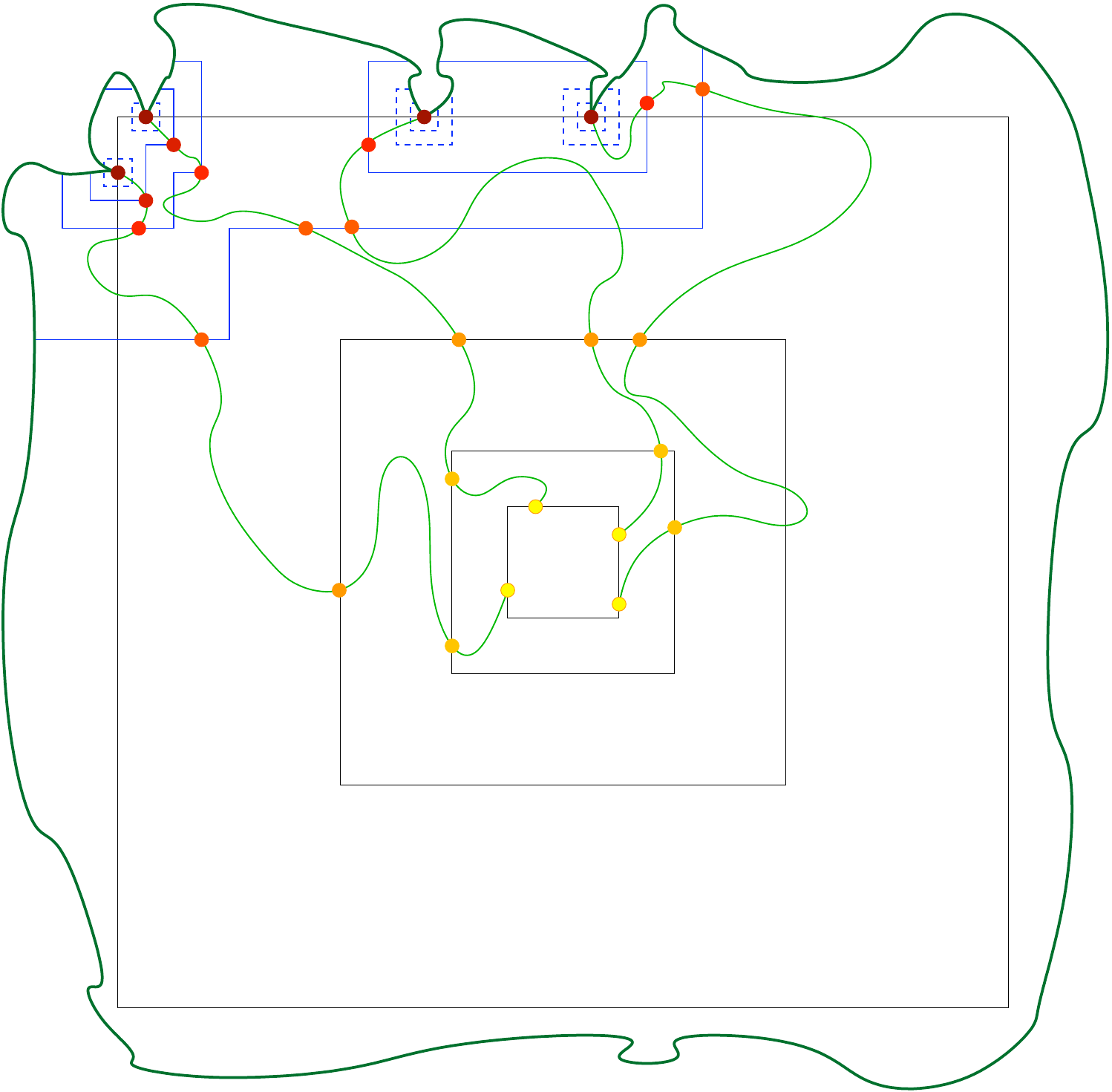}
}}
\caption{The bounding boxes associated to a starting configuration $\Theta$ with four faces, and a possible continuation $\Gamma$ of the four interfaces (with good separation at the end, radius $R_1$, but not everywhere on the way). The collection $\CC$ of bounding boxes has $M=3$ layers; there is one maximal box, giving the layer $\ell=0$, two boxes in layer $\ell=-1$ and one box in layer $\ell=-2$. The boxes associated to singletons do not  count.}
\label{f.SepPhenomBoxes}
\end{figure}

For each cluster $C\subseteq V$, the edges between nearest neighbour points inside $C$ can be represented by the shorter arc in $\p B(R_2)$ between them, and the union of these arcs is a connected arc, denoted by  $[C]\subseteq \p B(R_2)$. For $C=V$, it can happen that $[C]= \p B(R_2)$.
We now define a sequence of {\bf bounding boxes} for each cluster $C$,
$$
B_i(C):=\big\{x\in \mathcal{D}_{\Theta} : \dist_\infty(x,[C]) \leq s2^{j(C)-2+i} \big\}\,,
$$ 
for $i=0,1,\dots,J(C):=j(C^p)-j(C)-1$. These boxes are not always rectangles, for two reasons: the arcs $[C]$ might contain corners of $B(R_2)$, and part of the boundary of a box might come from $\Theta$; see Figure~\ref{f.SepPhenomBoxes}. We will call $B_0(C)$ the {inner bounding box} of $C$, and
$B_J(C)=B_{J(C)}(C)$ the {outer bounding box} of $C$. It may happen, of course, that $B_0(C)=B_J(C)$.

Note that  for disjoint clusters we always have $B_J(C_1)\cap B_J(C_2)=\emptyset$. Moreover, 
\begin{align}
\dist_\infty\big(B_J(C_1), B_J(C_2)\big) &\geq s ( 2^{\max\{j(C^p_1),j(C^p_2)\}-1} - 2^{j(C^p_1)-3} -  2^{j(C^p_2)-3} )\nonumber\\
& \geq  s 2^{j(C^p_i)-2}\qquad\textrm{for }i=1,2.\label{e.BCdist}
\end{align}
Similarly, if $C_1\subsetneqq C_2$, then 
\begin{equation}\label{e.BCinnerdist}
\dist_\infty\big(\p B_J(C_1),\p B_J(C_2)\big) \geq s2^{j(C_2)-3}\,.
\end{equation}
Furthermore, it is clear that the $\ell^\infty$-diameter of a bounding box satisfies
\begin{equation}\label{e.BCdiam}
s2^{j(C)-1+i} \leq \diam_\infty (B_i(C)) \leq Ns2^{j(C)+i}\,.
\end{equation}
These estimates mean that the boundaries of different bounding boxes are well-separated, with distances at least on the order of the sizes of the boxes themselves.


Let $\CC$ be the collection of all the bounding boxes $B_i(C)$ of all clusters with $j(C)\geq 1$ and $Ns2^{j(C)+i} \leq R_2/4$. The union of them is contained in $B(R_2)\setminus B(3R_2/4)$. In some special cases, when the set $V$ is very evenly distributed on the boundary, $\CC$ could be empty.

The boxes $B \in \CC$ can be classified into layers, as follows. We let $M=M(N,R_2,s):=\max\{m : Ns2^{m} \leq R_2/4 \}$, and then, for $B=B_i(C)$, we let $\ell(B):=j(C)+i-M \leq 0$, the {\bf layer} of $B$. Hence, the maximal boxes in $\CC$ are in the 0th layer, and $\CC$ naturally decomposes into a forest of rooted trees, with the 0th layer boxes as roots, and the inner bounding boxes of the smallest non-singleton clusters as the leaves, in the $(1-M)$th layer. See again Figure~\ref{f.SepPhenomBoxes}. 

Conditioned on the event $\mathcal{A}_\Theta(R_1, R_2)$, for each interface of $\Gamma$ started from $V$ there is a first time to exit the $\ell$th layer (the union of $\ell$th layer boxes of $\CC$). The set $V^\ell$ of these exit points naturally decomposes as $\bigcup_{k=1}^{k_\ell} V^\ell_k$ according to which box contains them on its boundary. For $\ell=-M$ we let $\bigcup_{k=1}^{k_{-M}} V^{-M}_k$ be the decomposition into the connected components of $G_1$. When $\CC=\emptyset$, we set $M=0$. 

We now define the {\bf relative qualities}
\begin{equation}\label{e.RelQual}
\RelQ(\ell):=\min_{1\leq k\leq k_\ell} \min_{x,y \in V^\ell_k} \frac{\dist(x,y)}{Ns2^{\ell+M-1}}
\end{equation}
for $\ell=-M,\dots,-1,0$, with $\RelQ(\ell):=0$ if not all $N$ interfaces manage to exit the $\ell$th layer. This is basically the usual quality (as in \cite[Appendix]{\SchrammSteif} or our Section~\ref{s.coupling}), just measured on the scale of the boxes in the $\ell$th layer. We extend the definition for $1\leq \ell < L:=\lfloor\log_2 (R_2/R_1)\rfloor$ by $\RelQ(\ell):=\Qual(\Gamma(2^{-\ell}R_2))$, where $\Gamma(R)$ is the $N$-tuple of interfaces $\Gamma$ restricted to $A(R,R_2)$, with quality measured at radius $R$. Note that $\RelQ(0)$ is comparable to the original quality, i.e., it is measured on the scale of $R_2$. Finally, we set $\RelQ(L):=\Qual(\Gamma)$, so that our goal becomes to prove 
$$\PB{\RelQ(L) >\frac{1}{N} \md \RelQ(L) > 0 } > c(N) >0\,.$$

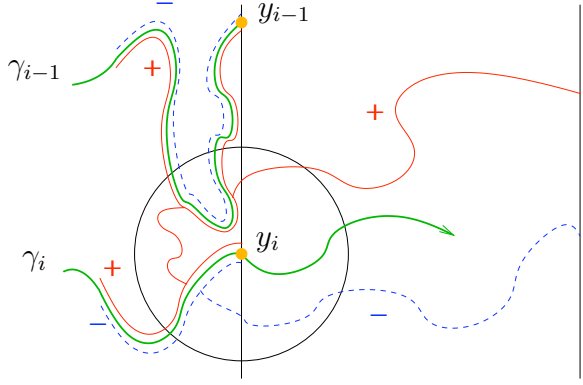
\begin{figure}[htbp]
\SetLabels
(-0.03*.8)$\gamma_{i-1}$\\
(-0.03*.3)$\gamma_i$\\
(.43*.95)$y_{i-1}$\\
(.4*.35)$y_i$\\
\endSetLabels
\centerline{
\AffixLabels{
\includegraphics[height=2 in]{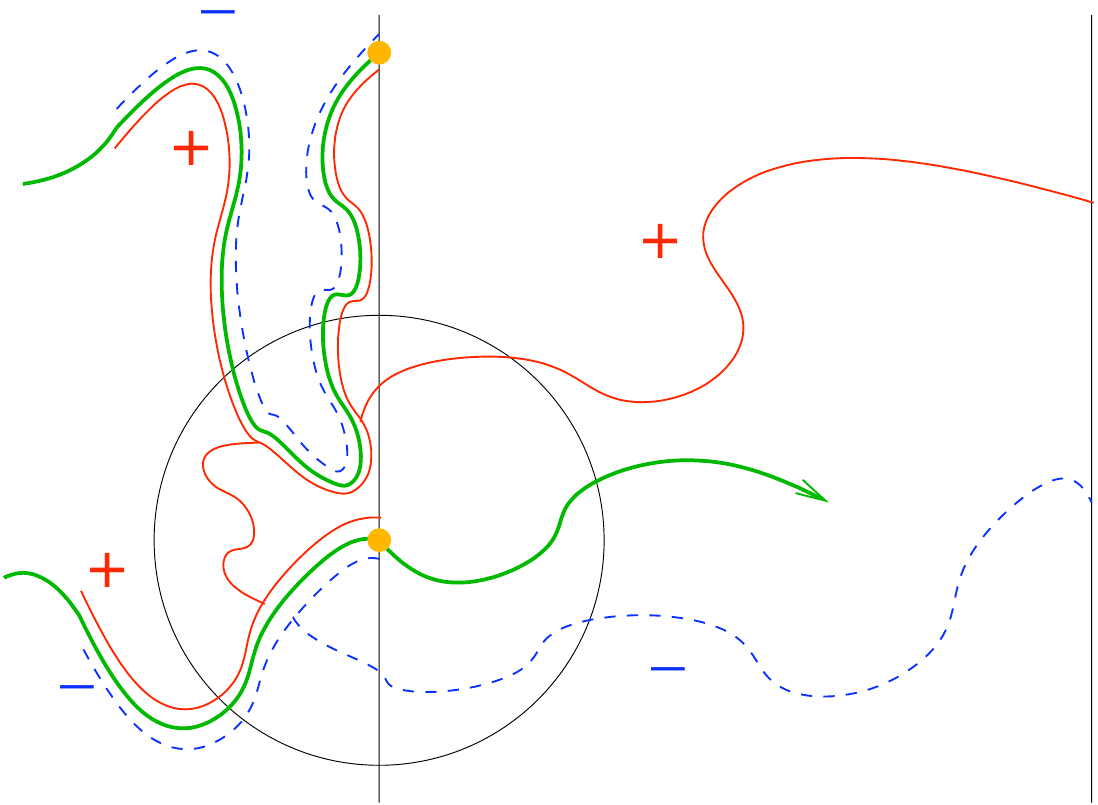}
}}
\caption{Good quality implies a large probability of successfully gluing arms and continuing the interface nicely.}
\label{f.gluing}
\end{figure}

We will use two lemmas, which are direct adaptations of widely used basic results regarding separation of interfaces to our situation. For the first one, see, e.g., \cite[Lemma A.3]{\SchrammSteif}; the second one comes from \cite[Lemma A.2]{\SchrammSteif} or \cite[Lemma 2]{\KestenScaling}. Both are easily proved by RSW and FKG techniques. That the RSW technology works nicely in our situation is ensured by the inequalities (\ref{e.BCdist}, \ref{e.BCinnerdist}, \ref{e.BCdiam}). Note that a good quality implies that the interfaces can be continued nicely with RSW and FKG gluing techniques, despite the fact that we require only that the endpoints have a large distance, while a piece of a neighbouring interface $\gamma_{i\pm 1}$ still can come close to an endpoint $y_i$. The reason that this works is that the neighbouring interface has the right colour on its side closer to $y_i$. See Figure~\ref{f.gluing}.

\begin{lemma}\label{l.moveon}
For any $\rho>0$ there is a $c_1(\rho)=c_1(N,\rho)>0$ with 
\begin{equation}
\PB{\RelQ(\ell+1) > \frac{1}{N} } > c_1(\rho)\, \Pb{\RelQ(\ell)>\rho}\label{e.moveon}
\end{equation}
for any $\ell$.
\end{lemma}

Two immediate consequences are that there is a $c_2=c_2(N)>0$ with
\begin{equation}\label{e.goon}
\PB{ \RelQ(\ell+1) > \frac{1}{N} \md \RelQ(\ell)>\frac{1}{N} } > c_2 \,,
 \end{equation}
and, since $\RelQ(-M)=2/N$, 
\begin{equation}\label{e.goend}
\PB{ \RelQ(\ell) >\frac{1}{N}} > {c_2}^{\ell+M}\,.
\end{equation}
 
\begin{lemma}\label{l.nobadsep}
For $\delta>0$, if $\rho=\rho(N,\delta)>0$ is small enough, then the probability of the event $\{$the total set
of interfaces in layer $\ell$ (i.e., from the boundary of layer $\ell-1$ to the boundary of layer $\ell$) is non-empty and has relative quality at most  $\rho \}$ is less than $\delta$.
\end{lemma}

Now, for a small $\rho>0$ whose value will be fixed later,
\begin{align*}
\Pb{\RelQ(L) > 0} & < \Pb{\RelQ(L-1) > 0}\\ 
&= \Pb{ \RelQ(L-1) > \rho} + \Pb{\rho \geq \RelQ(L-1) > 0} \\
& < \Pb{\RelQ(L-1) > \rho} + \delta\,\Pb{ \RelQ(L-2) > 0}, \text{ by Lemma~\ref{l.nobadsep}}\\
& < \sum_{\ell=0}^{L+M-1} \delta^\ell\, \Pb{\RelQ(L-1-\ell) > \rho} + \delta^{L+M-1}, \text{ by iterating previous step}\\
& < \sum_{\ell=0}^{L+M-1} \delta^\ell\, \PB{\RelQ(L-\ell)>\frac{1}{N}} c_1(\rho)^{-1} + \delta^{L+M-1}, \text{ by~(\ref{e.moveon})}\\
& < \sum_{\ell=0}^{L+M-1} \delta^\ell\, \PB{\RelQ(L)>\frac{1}{N}} c_2^{-\ell} c_1(\rho)^{-1} + \delta^{L+M-1}, \text{ by (\ref{e.goon})}\\
& < \PB{\RelQ(L)>\frac{1}{N}} \left( \frac{1}{c_1(\rho)}\, \sum_{\ell=0}^{L+M-1} \left(\frac{\delta}{c_2}\right)^\ell + \left(\frac{\delta}{c_2}\right)^{L+M-1} \right) \,, 
\end{align*}
by applying (\ref{e.goend}) to the second term. Hence, if we choose $\rho$ so small that the $\delta$ in Lemma~\ref{l.nobadsep} becomes smaller than $c_2/2$, then the last line becomes at most
$$
\PB{\RelQ(L)>\frac{1}{N}}\, K(N)\,,
$$
for some constant $K(N)$, and we have proved Proposition~\ref{p.separation}.\qed
\medskip

To deduce Lemma~\ref{l.separation} from the proposition, one might worry that conditioning on the remote ``future'' $r$ will manifest itself only gradually, so will not produce a good separation already at $R_1$. But there is no such problem:

\proofof{Lemma~\ref{l.separation}} 
We set now $N=4$. Let us denote by $\ev{Q}$ the event $\{\Qual(\Gamma(R_1)) >1/4\}$,  by $\ev{A}$ the event $\mathcal{A}_\Theta (r,R_2)$,  by $\ev{B}$ the event  $\mathcal{A}_\Theta (R_1,R_2)$, and by $\ev{C}$ the 4-arm event $\A_4(r,R_1)$. Obviously, $\Ps{\ev{A}} \leq \Ps{\ev{B}}\Ps{\ev{C}}$, hence 
$$\Ps{\ev{A} \md \ev{B}} \leq \alpha_4(r,R_1)\,.$$
Also, by the usual RSW and FKG gluing techniques (similarly to Lemma~\ref{l.moveon}), 
$$
\Ps{\ev{A} \md \ev{Q}} \asymp \alpha_4(r,R_1)\,.
$$
Therefore,
$$
\Ps{\ev{Q} \md \ev{A}} = \frac{ \Ps{\ev{Q}, \ev{A}} }{ \Ps{\ev{A}} } \geq \frac{ \Ps{\ev{A} \md \ev{Q}} \Ps{\ev{Q}} }{\Ps{\ev{B}}\Ps{\ev{C}}} \asymp \Ps{\ev{Q}\md \ev{B}}\,,
$$
which is at least a constant $c>0$ by Proposition~\ref{p.separation}, and we have proved the first statement of the lemma: $\Ps{\ev{Q}\md \ev{A}}>c'>0$.

We still have to add $U_\Theta=\tau \in \{0,1\}$ to the conditioning. The faces $\Theta$ and the four interfaces $\Gamma(R_1)$ together induce a configuration $\Theta'$ of four faces around $B(R_1)$. Besides $\ev{A}$, condition also on any $\Theta'$ that satisfies $\Qual(\Theta')=\Qual(\Gamma(R_1))>1/4$, which is just the event $\ev{Q}$. By Proposition~\ref{p.separation}, we have $\Pb{\Qual(\Gamma(r)) > 1/4 \md \ev{A}, \Theta'} > c'' > 0$, hence
\begin{align*}
\PB{\Qual(\Gamma(r))>\frac 1 4,\, \ev{Q} \md \ev{A}} &= \sum_{\Qual(\Theta')>1/4} \PB{\Qual(\Gamma(r))>\frac 1 4 \md \ev{A}, \Theta'} \cdot \Pb{\Theta' \md \ev{A}}  \\
& > c''\, \Pb{\ev{Q} \md \ev{A}} = c'' c' > 0\,.
\end{align*}
Let us denote by $\ev{Q}^*$ the event $\{\Qual(\Gamma(r))>1/4\}$ and by $\ev{A}^\tau$ the event $\ev{A} \cap \{U_\Theta=\tau\}$. So, we have just shown that  $\Pb{\ev{Q}^*,\, \ev{Q} \md \ev{A}} > c'c''$. On the other hand,
with the usual gluing arguments, we have $\Pb{ \ev{A}^\tau \md \ev{A},\,\ev{Q}^*} > c''' >0$. Altogether, $\Pb{\ev{Q},\,\ev{A}^\tau \md \ev{A}} > c'c''c'''$. Therefore,
$$
\Pb{\ev{Q} \md \ev{A}^\tau} =  \frac{ \Pb{\ev{Q},\, \ev{A}^\tau \md \ev{A}} }{ \Pb{\ev{A}^\tau \md \ev{A}} } > c'c''c''' > 0\,,
$$
and Lemma~\ref{l.separation} is proved.
\qed
\medskip

Let us conclude with a discussion of what even stronger versions one might imagine. Two overly naive questions closely related to each other: (a) Could the lower bound $c(N)$ on the success probability approach 1 as $R_2/R_1$ tends to $\infty$? (b) Do the interfaces become well-separated quickly, say, by radius $R_2/2$, and remain so all along till $R_1$? The answer to both, of course, is ``no'', because there are many scales, and even from a small quality at some scale we still have some chance to fulfill the conditioning. 
Another direction is to ask how $c(N)$ depends on $N$. From the proof strategy it appears that a larger $N$ does not make it harder for neighboring interfaces to be about $1/N$-separated. However, a larger $N$ also means more randomness, hence the typical quality at $R_1$ should not be on the order of $1/N$, only at most $1/(N\log^b N)$, so $c(N)$ is not expected to be independent of $N$. What seems likely is that if we look at the interface endpoints at $R_1$ locally, around some point of $\p B(R_1)$, rescaled by $N$, then we get a translation-invariant point process on $\R$ in the $N\to\infty$ limit, which should have some repulsion similar to the eigenvalue processes in random matrix theory. See \cite{\DubedatCommutation} for the closely related question of growing SLE$_\kappa$ curves conditioned on not to intersect each other.

Finally, there could be a ``super-strong'' version that we do not know how to prove, though we expect it to be true: 

\begin{conjecture}[Super-strong Separation Lemma] Given an arbitrary starting configuration $\Theta$ of faces around $B(R_2)$ (with any number of faces), condition on at least four of the interfaces starting from the endpoints of the faces to reach $B(R_1)$. Under this conditioning, the set of interfaces reaching $R_1$ have positive quality with positive probability.
\end{conjecture}

If we condition on a $\Theta$ with $N$ faces, then we can apply the Strong Separation Lemma for 4 interfaces, and lose a factor of $N \choose 4$ in the probability. But for a general $\Theta$, the number of faces could be even on the order of $R_2/\eta$, so the loss would be huge.

\subsection{Brownian intersections}

We will now sketch a proof of the following basic result:

\begin{proposition}[Lawler's Separation Lemma \cite{\LawlerStrict}]\label{p.LawlerSep}
Let $B(r)$ and $B(R)$ be the Euclidean disks around the origin, with $0< 2r < R$, and let $\Gamma=\{\gamma_1(t),\gamma_2(t)\}$ be two independent planar Brownian motions, started at arbitrary but different points $\gamma_i(0) \in\p B(R)$, stopped when reaching $\p B(r)$, at stopping times $\tau_i(r)$. Then there is some absolute constant $c>0$ such that 
$$
\PB{ \Qual(\Gamma(r)) > 1/4 \md  \gamma_1[0,\tau_1(r)]\cap \gamma_2[0,\tau_2(r)]=\emptyset } > c\,,
$$
where $\Gamma(s)$ is the pair of Brownian motions, each stopped at $\tau_i(s)$, $r\leq s\leq R$, and $\Qual(\Gamma(s))$ is the quality at $\p B(s)$, defined by 
$$
\Qual(\Gamma(s))=\frac{1}{s} \min_{i=1,2} \dist\Big(\gamma_i(\tau_i(s)),\ \gamma_{3-i}[0,\tau_{3-i}(s)]\Big)\,.
$$ 
\end{proposition}

Lemma 4.2 in \cite{\LawlerStrict} is more general than this, involving packets of Brownian motions, but the ideas are exactly the same, so we chose the notationally simplest possible setting. Also, the Brownian motions could be started outside $B(R)$, and we could condition on their trajectories up to their first hitting of $\p B(R)$; the proof would be exactly the same.

As opposed to percolation interfaces, where the closeness of one side of a neighbouring interface does not  present any danger to fulfilling the conditioning, see Figure~\ref{f.gluing}, a Brownian motion could hit another trajectory anywhere, causing failure. That is why the notion of quality is strengthened as above.

\proof
The following two lemmas are the analogs of Lemmas~\ref{l.moveon} and~\ref{l.nobadsep}. For simplicity, we state them for ordinary qualities instead of the relative ones introduced in~(\ref{e.RelQual}).

\begin{lemma}\label{l.BMmoveon} 
For any $\rho>0$ there is an $\eps(\rho)>0$ with 
$$
\PB{\Qual(\Gamma(s/2)) > \frac{1}{4} } > \eps(\rho)\, \Pb{\Qual(\Gamma(s)) >\rho}\,,
$$
whenever $r\leq s/2 < s \leq R$.
\end{lemma}

\begin{figure}[htbp]
\SetLabels
(.66*.72)$\gamma_1$\\
(0.61*0.02)$\gamma_2$\\
(.04*.5)$s$\\
(.18*.5)$s/2$\\
(.40*.65)$U_1$\\
(.34*.22)$U_2$\\
(.28*.66)$V_1$\\
(.24*.30)$V_2$\\
(.86*.85)$\gamma_1$\\
(0.88*0.2)$\gamma_2$\\
(.845*.56)$\gamma^*$\\
(.84*.4)$\beta$\\
(0.91*0.44)$y$\\
(0.76*0.39)$B_z(s/4)$\\
\endSetLabels
\centerline{
\AffixLabels{
\includegraphics[height=2.5in]{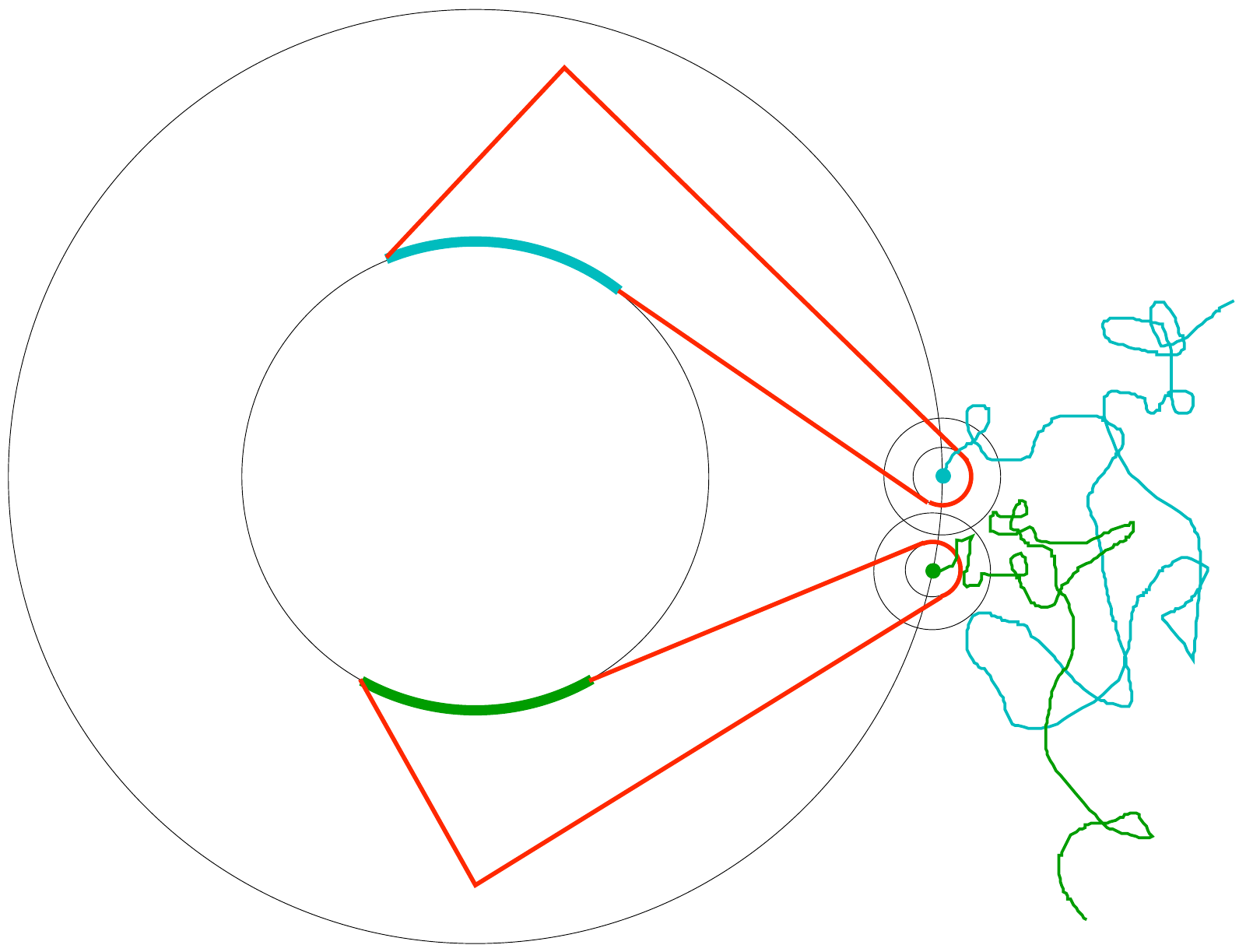}
\hskip .4 in
\includegraphics[height=2.5in]{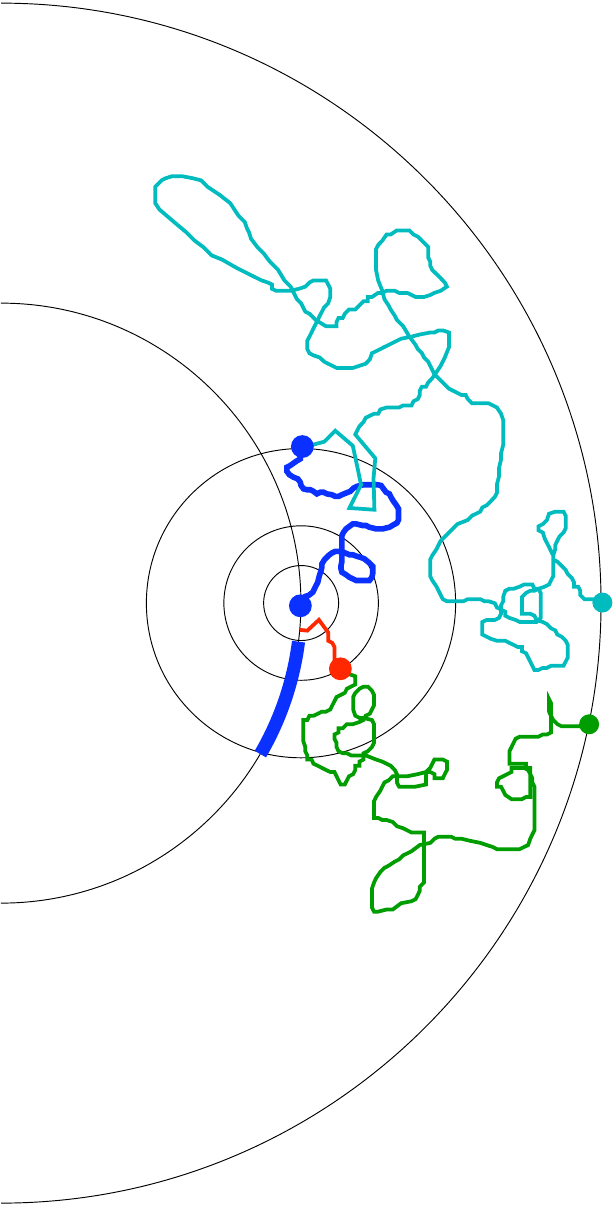}}
}
\caption{Proving Lemma~\ref{l.BMmoveon} and Lemma~\ref{l.BMnobadsep}.}
\label{f.BM}
\end{figure}

\proof 
Given any pair $\Gamma(s)$ satisfying $\Qual(\Gamma(s)) >\rho$, it is straightforward to construct two disjoint domains in the plane, $U_1$ and $U_2$, with the following properties, as depicted on the left part of Figure~\ref{f.BM}:
the interior of each $U_i$ is disjoint from $\gamma_{3-i}[0,\tau_{3-i}(s)]$ and from $B(s/2)$, 
but contains the $s\rho/2$-ball around $\gamma_i(\tau_i(s))$; $\p U_i$ contains an arc of $\p B(s/2)$, called $V_i$;
the harmonic measure of $V_i$ for Brownian motion in $U_i$ started from $\gamma_i(\tau_i(s))$ is at least $h(\rho)>0$, independently of $s$; and
finally, $V_1$ and $V_2$ are well-separated in the sense that if each $\gamma_i$ exits $U_i$ at $V_i$, then we necessarily have $\Qual(\Gamma(s/2))>1/4$.
This proves the lemma with $\eps(\rho)=h(\rho)^2$.\qed
\medskip

\begin{lemma}\label{l.BMnobadsep} 
For any $\delta>0$, if $\rho=\rho(\delta)>0$ is small enough, then, for arbitrary but different starting points $\gamma_i(0) \in \p B(s)$,
$$
\PB{ \gamma_1[0,\tau_1(s/2)]\cap \gamma_2[0,\tau_2(s/2)]=\emptyset,\  \Qual(\Gamma(s/2)) \leq \rho} < \delta\,.
$$
\end{lemma}

\proof 
Let $\gamma^*$ be the subpath of $\gamma_1[0,\tau_1(s/2)]$ from the last time it entered the ball of radius $s/4$ around $z:=\gamma_1(\tau_1(s/2))$ till its end $\tau_1(s/2)$. If the event $\{\Qual(\Gamma(s/2)) \leq \rho\}$ occurs, then there is a first time that $\gamma_2$ enters the ball of radius $s\sqrt{\rho}$ around $z$; let us denote the entering point by $y$. The connected component of $B_z(s/4) \setminus (B(s/2) \cup \gamma^*)$ that contains $y$ will be denoted by $D$, and the part of $\p B(s/2) \cap \p D$ that has distance between $s \rho$ and $s/4$ from $z$ will be denoted by $\beta$. See the right hand picture of Figure~\ref{f.BM}. It is now an easy exercise about harmonic measure to show that, regardless of what $\gamma^*$ and $y$ are, the probability that a Brownian motion started at $y$ exits $D$ not through $\gamma^* \cup \beta$ is close to zero if $\rho$ is small, which implies the lemma. 
\qed
\medskip

Based on these two lemmas, the rest of the proof of Proposition~\ref{p.LawlerSep} is identical to that of Proposition~\ref{p.separation}.
\qed

\bibliographystyle{../../halpha}
\addcontentsline{toc}{section}{Bibliography}

\bibliography{Pivotal.ref}

\ \\
\ \\
{\bf Christophe Garban}\\
ENS Lyon, CNRS\\
\url{http://perso.ens-lyon.fr/christophe.garban/}\\
Partially supported by ANR grant BLAN06-3-134462.\\
\\
{\bf G\'abor Pete}\\
Institute of Mathematics, Technical University of Budapest\\
\url{http://www.math.bme.hu/~gabor}\\
Supported by an NSERC Discovery Grant at the University of Toronto, and an EU Marie Curie International Incoming Fellowship at the Technical University of Budapest\\
\\
{\bf Oded Schramm} (December 10, 1961 -- September 1, 2008)\\
Microsoft Research\\ 
\url{http://research.microsoft.com/en-us/um/people/schramm/}\\

\end{document}